\documentclass[11pt]{amsart}
\usepackage{amsmath, amsfonts, amsthm, amssymb}
\usepackage[all,cmtip,line]{xy}
\usepackage{array}
\usepackage{enumerate}
\usepackage{srcltx} 

\newtheorem{thm}{Theorem}[section]
\newtheorem{lem}{Lemma}[section]
\newtheorem{prop}[lem]{Proposition}
\newtheorem{cor}[lem]{Corollary}

\newtheorem{defn}[lem]{Definition}
\newtheorem{rem}[lem]{Remark}
\newtheorem{conj}[lem]{Conjecture/Question}
\newtheorem{claim}[lem]{Claim}
\numberwithin{equation}{section}

\newcommand \eps{\varepsilon}

\newlength{\originalbase}
\title{Interior $W^{2,p}$ estimate for small perturbations to the complex Monge-Ampere equation}

\author{Jingrui Cheng, Yulun Xu}

\date{Nov 2022}

\begin{document}

\maketitle

\begin{abstract}
Let $w_0$ be a bounded,  $C^3$,  strictly plurisubharmonic function defined on $B_1\subset \mathbb{C}^n$.  Then $w_0$ has a neighborhood in $L^{\infty}(B_1)$ with the following property: for any continuous,  plurisubharmonic function $u$ in this neighborhood solving $1-\eps \le MA(u)\le 1+\eps$,  one has $u\in W^{2,p}(B_{\frac{1}{2}})$,  as long as $\eps>0$ is small enough depending only on $n$ and $p$.  This partially generalizes Caffarelli's interior $W^{2,p}$ estimates for real Monge-Ampere to the complex version.
\end{abstract}

\section{Introduction}
Monge-Ampere equations are second-order partial differential equations whose leading term is the determinant of the Hessian of a real unknown function.  The Hessian is required to be positive or at least nonnegative,  so that the equations are elliptic or degenerate elliptic.  Monge-Ampere equations can be divided into real or complex,  depending on whether one is considering real Hessian or complex Hessian.  In the real case,  the Hessian is $u_{ij}$,  so that the positivity of the Hessian is a convexity condition.  In the complex case,  the Hessian is $u_{i\bar{j}}$,  and its positivity is a plurisubharmonicity condition.

For both real and complex Monge-Ampere,  the existence and regularity theory with smooth data has been well established.  In the real case,  it is proved by Caffarelli-Nirenberg-Spruck \cite{CNS} on smooth strictly convex domains on $\mathbb{R}^n$.  
In the complex case,  the foundations of an existence and regularity theory were laid out by Yau \cite{Yau} in the setting of a compact K\"ahler manifold,  and by Caffarelli-Kohn-Nirenberg-Spruck \cite{CKNS},  in the setting of a smooth pseudo-convex domain.

Another important aspect about the Monge-Ampere equations is their apriori estimates,  starting with interior ones.  In general,  the results known for the real case is much stronger than the complex case,  due to the fact that the solution being convex gives much more stringent constraint than being plurisubharmonic.
For example,  the interior gradient estimate for real Monge-Ampere equation is more or less a trivial matter (if we know the solution is bounded),  since the underlying solution considered is convex.
This is not the case for complex Monge-Ampere,  and the boundedness of plurisubharmonic function only gives the gradient being in $L^2$.

Arguably the most important estimate of Monge-Ampere is to get second derivative estimates.  If we get such estimates in $L^{\infty}$,  then the equation becomes unformly elliptic and the standard theory can apply.  On the other hand,  the bad news is that for both real and complex Monge-Ampere equations,  there are no purely interior $C^2$ estimates.  Indeed,  having a convex solution to $\det u_{ij}=1$ in a domain doesn't imply $u\in C^2$ in the interior,  due to a counterexample by Pogorelov \cite{P2} (a counterexample for the complex version is given by He in \cite{He}).  In general,  one needs to impose some boundary conditions (say,  $u=0$ on the boundary),  in order to conclude that $D^2u$ is bounded in the interior (Pogorelov's estimate \cite{P}).  Based on that,  Caffarelli's proved the following interior $W^{2,p}$ estimate when the right hand side is a small perturbation of a constant:
\begin{thm}\label{t1.1}
Let $\Omega\subset \mathbb{R}^n$ be a convex domain such that $B_1\subset \Omega\subset B_n$ and $u$ is a weak solution to $\det u_{ij}=f$ with $|f-1|<\eps$ and $u=0$ on $\partial \Omega$.  Then for any $1<p<\infty$,  if $\eps>0$ is small enough depending only on $p$ and $n$,  then $||u||_{W^{2,p}(B_{\frac{1}{2}})}$ can be bounded by a constant depending only on $p$ and $n$.
\end{thm}
In the above theorem,  the weak solution is defined using the measure of the image of the gradient mapping.  One could also replace the boundary condition by a $strict\,\,convexity$ assumption,  meaning that the supporting plane of a convex function touches the function only at one point.  That is,  we have:
\begin{thm}\label{t1.2}
Let $\Omega\subset \mathbb{R}^n$ be a bounded convex domain,  and $u$ is a weak solution to $\det u_{ij}=f$ with $|f-1|<\eps$ which is strictly convex.  Then for any compact subdomain $\Omega'$,  one has $||u||_{W^{2,p}(\Omega')}\le C$,  where $C$ depends on $p$,  $n$,  $dist(\Omega',\partial \Omega)$,  the modulus of strict convexity of $u$,  as long as $\eps$ is small enough depending only on $p$ and $n$.
\end{thm}
Theorem \ref{t1.2} actually follows from Theorem \ref{t1.1}.  Indeed,  for any $x_0\in\Omega'$,  we can take $l_{x_0}$ to be a linear function touching $u$ from below,  then $S_c:=\{u(x)\le l_{x_0}(x)+c\}$ will be contained in $\Omega$ for $c$ small enough,  due to the strict convexity.  Then one can normalize $S_c$ to be in the situation of Theorem \ref{t1.1}.\\

The goal of this paper is to generalize (partially) Theorem \ref{t1.1} and Theorem \ref{t1.2} to the complex Monge-Ampere equations.  More precisely,  we show that
\begin{thm}\label{main theorem baby}
Let $\Omega\subset \mathbb{C}^n$ be a bounded domain with $B_{1-\gamma_0}\subset \Omega\subset B_{1+\gamma_0}$ for some $\gamma_0>0$.  Let $u\in C^2(\Omega)\cap PSH(\Omega)\cap C(\bar{\Omega})$ be such that $1-\eps \le \det u_{i\bar{j}}\le 1+\eps$ in $\Omega$ and $u=0$ on $\partial \Omega$.  
Given $1<p< \infty$,  if $\gamma_0$ is small enough depending only on $n$,  and $\eps$ small enough depending only on $n$ and $p$,  then 
\begin{equation*}
||u||_{W^{2,p}(B_{\frac{1}{2}})}\le C,  \,\,\,\,\,||\sum_i\frac{1}{u_{i\bar{i}}}||_{L^p(B_{\frac{1}{2}})}\le C,
\end{equation*}
where the constant $C$ depends only on $n$ and $p$.
\end{thm}
This theorem should be understood as the analogue of Theorem \ref{t1.1}.
More generally,  we have:
\begin{thm}\label{main theorem}
Let $w_0$ be a smooth function in the unit ball such that for some $C_0>1$:
\begin{equation*}
      \frac{1}{C_0}I\le (w_0)_{z_i\bar{z}_j}\le C_0I,\,\,\,|D^3w_0|\le C_0\text{ in $B_1$.}
\end{equation*}
Then there exists $\delta_0>0$ small enough,  depending only on $C_0$ and $n$,  such that for all $u\in PSH(B_1)\cap C(B_1)$ with $|u-w_0|\le \delta_0$ on $B_1$,  solving $1-\eps\le MA(u)\le 1+\eps$,  we have $u\in W^{2,p}(B_{\frac{1}{2}})$ and $\sum_i\frac{1}{u_{i\bar{i}}}\in L^p(B_{\frac{1}{2}})$,  as long as $\eps$ is small enough depending only on $n$ and $p$.
\end{thm}
 In the above,  $MA(u)$ is the complex Monge-Ampere operator defined for continuous plurisubharmonic functions,  in the Bedford-Taylor sense (see \cite{BT}),  so that $MA(u)=\det u_{i\bar{j}}$ when $u\in C^2$.

If one compares Theorem \ref{main theorem baby} to Theorem \ref{t1.1},  or Theorem \ref{main theorem} to Theorem \ref{t1.2},  the biggest difference is that we have to assume our solution is close to a smooth plurisubharmonic function.  The reason we have to make this assumption is related to whether one has Pogorelov type estimates for complex Monge-Ampere equations,  which has been open until now.  Note that the interior $C^2$ estimates for the complex Monge-Ampere equation with zero boundary values were studied by F.  Schulz in \cite{Schulz},  using the integral approach of N.  M.  Ivochikina \cite{Iv} for real Monge-Ampere equations. 
However,  the proof in \cite{Schulz} is not complete,  which was first pointed out by Blocki \cite{Blocki}.
We will comment more on the technical aspect later,  but for now,  let us present some direct consequences of our main theorem,  which seems new and interesting.

First we observe that Theorem \ref{main theorem} would give us the following result in the manifold setting:
\begin{cor}\label{c1.2}
 Let $(M,\omega_0)$ be a compact K\"ahler manifold.  
Let $\varphi\in PSH(M,\omega_0)\cap C(M)$ be the solution to:
\begin{equation*}
(\omega_0+\sqrt{-1}\partial\bar{\partial}\varphi)^n=f\omega_0^n,\,\,\,\omega_0+\sqrt{-1}\partial\bar{\partial}\varphi>0,
\end{equation*}
where $|f-1|<\eps$ and $\int_Mf\omega_0^n=\int_M\omega_0^n$.
Let $1<p<\infty$,  then we have that $\varphi\in W^{2,p}(M)$ as long as $\eps$ is small enough depending only on $p$,  $n$ and the background metric $\omega_0$.
\end{cor}
Similar results would also hold for the setting of bounded domains.  In other words,  we have:

\begin{cor}\label{c1.2N}
Let $\Omega\subset \mathbb{C}^n$ be a bounded domain.  Let $u_0\in C^3(\Omega)\cap C(\bar{\Omega})\cap PSH(\Omega)$ be the solution to $\det(u_0)_{i\bar{j}}=f_0>0$ in $\Omega$ and $u_0|_{\partial \Omega}=\varphi_0$.  Let $u\in C(\bar{\Omega})\cap PSH(\Omega)$ be the solution to $MA(u)=f$ and $u|_{\partial \Omega}=\varphi$ such that $|f-f_0|<\eps$,  $|\varphi-\varphi_0|<\eps$.  Let $\Omega'$ be a compact subdomain of $\Omega$,  then we have $u\in W^{2,p}(\Omega')$,  as long as $\eps$ is small enough depending only on $\Omega',\,\Omega$,  the $C^3$ bound and complex Hessian lower bound of $u_0$ in a neighborhood of $\Omega'$,  
$p$ and $n$.
\end{cor}
One more application of Theorem \ref{main theorem} is the following Liouville theorem for entire scalar flat metric on $\mathbb{C}^n$,  which is a generalization of a result by Yu Wang \cite{Wang}:
\begin{cor}\label{c1.3}
Let $u$ be a $C^2$ plurisubharmonic function on $\mathbb{C}^n$.  Denote $\omega_u=\sqrt{-1}\partial\bar{\partial}u$ and assume that $\omega_u$ is scalar flat,  namely
\begin{equation*}
\sum_{i,j=1}^nu^{i\bar{j}}\partial_{i\bar{j}}\big(\log \det u_{a\bar{b}}\big)=0.
\end{equation*}
Then there exists $\eps_n>0$ small enough depending only on $n$,  such that we can deduce $u$ is quadratic,  provided that:
\begin{enumerate}
\item $\overline{\lim}_{r\rightarrow \infty}\frac{\sup_{B_r}|u(z)-|z|^2|}{r^2}\le \eps_n$,
\item $1-\eps_n\le \det u_{a\bar{b}}\le 1+\eps_n$.
\end{enumerate}
\end{cor}
The result by Yu Wang \cite{Wang} is a special case of the above Corollary with $\det u_{a\bar{b}}=1$.\\

Finally we observe the following $C^{2,\alpha}$ estimate for complex Monge-Ampere.  For $K>0$,  $0<\alpha<1$ and $C_1>1$,   we define the following class of functions:
\begin{equation*}
\mathcal{F}(K,\alpha,C_1)=\{\text{$f$ is defined on $B_1$}:||f||_{\alpha,B_1}\le K,\,\,\,\frac{1}{C_1}\le f\le C_1\text{ on $B_1$}\}.
\end{equation*}
\begin{cor}\label{c1.4}
Let $w_0$ be as in Theorem \ref{main theorem}.  Then there exists $\delta_0>0$,  depending only on $C_0$,  $n$,  $K$,  $\alpha$,  $C_1$,  such that for any $u\in PSH(B_1)\cap C(B_1)$ with $|u-w_0|\le \delta_0$ on $B_1$ and solving $MA(u)=f$ for some $f\in\mathcal{F}(K,\alpha,C_1)$,  we have $u\in C^{2,\alpha}(B_{\frac{1}{2}})$.
\end{cor}

Next we would like to explain the ideas of proof for Theorem \ref{main theorem baby}.  The heart of the idea is from Caffarelli's paper \cite{Ca1}  which we explain first.  Since the solution $u$ is strictly convex,  we may consider sections of $u$ of the form $S_c(x_0):=\{(u-l)(x)\le u(x_0)+c\}$ which is strictly contained in $\Omega$,  where $l(x)$ is the supporting linear function of $u$ at $x_0$ and $c>0$ is called the ``height" of the section.  Now we solve $\det w_{ij}=1$ on this open set,  equaling $u$ on the boundary.  From Pogorelov estimate,  we know that $w$ is smooth in the interior.  By doing Taylor expansion for $w$,  we find that the sections of $w$ will be close to ellipsoids. 
On the other hand,  since $f$ is close to 1,  we also have $u$ is very close to $w$ by maximum principle,  hence the sections of $u$ are close to ellipsoids as well.  
If the shape of the ellipsoids are comparable to a ball for heights going to 0,  then the second derivatives are under control at that point.  The whole point of $W^{2,p}$ estimate is then to estimate the measure of the set where the shape of such ellipsoids loses control. (which is reflected by the opening of the paraboloid touching $u$ from below) For this purpose,  we will need a version of Vitali's covering lemma,  but adapted to sections.  To establish the covering lemma for sections,  a  crucial property we need is the following engulfing property: 

If $S_{\mu_1}(x_1)\cap S_{\mu_2}(x_2)\neq \emptyset$,  with $\mu_1\le \mu_2$,  then $S_{\mu_1}(x_1)\subset S_{10\mu_2}(x_2)$.  
 
This property would be a result of compactness.  Indeed,  if $u$ were the standard solution $\frac{1}{2}|x|^2$,  we would have $S_{\mu_1}(x_1)=B_{\sqrt{\mu_1}}(x_1)$,  $S_{\mu_2}(x_2)=B_{\sqrt{\mu_2}}(x_2)$ and the engulfing property indeed holds.  We can still expect this property if $u$ is close to a quadratic polynomial.

We follow similar lines of argument in the proof of Theorem \ref{main theorem baby}.  The first hurdle we face is to take sections with $u$.  Unlike the convex function,  given $x_0\in \Omega$,  it is not clear whether one can find a pluriharmonic function $h$,  for which $\{u-h<u(x_0)+c\}$ is compactly contained in $\Omega$ for $c>0$ small.  Even though this is not clear in general,  we show that,  however,  it is indeed possible if $u$ is close to a smooth plurisubharmonic function whose complex Hessian has a lower bound.

Next we need to solve the Dirichlet problem $\det w_{i\bar{j}}=1$,  $w=u$ on the boundary of a small section of $u$,  similar to what we did in the real case.  The problem we are facing now,  is that we do not know if $w$ is smooth,  since Pogorelov's estimate is not known for the complex Monge-Ampere equations.  However,  if $u$ is close to a smooth,  strictly plurisubharmonic function,  then the section defined by $u$ will be close to an ellipsoid from the very beginning.  
This will allow us to use Savin's perturbation result to conclude that $w$ is indeed smooth in the interior.  In this paper,  we use an induction process to construct sections $S_{\mu}(x_0)$ for $\mu>0$ and small,  which takes the form $\{(u-h_{\mu,x_0})(z)\le u(x_0)+\mu\}$,  where $h_{\mu,x_0}(z)$ is pluriharmonic.  Moreover,  we also show that $S_{\mu}(x_0)$ remains close to an ellipsoid in the induction process.

A drawback with our construction is that it is highly non-canonical,  since it relies on solving Dirichlet problem on a sequence of smaller and smaller sections ``centered at" $x_0$.  This construction of $S_{\mu}(x_0)$ does not commute with linear transformations we use to normalize the ellipsoids.  We explain this matter in greater detail in Subsection 4.2,  under Proposition \ref{p3.15}.

The fact that $S_{\mu}(x_0)$ is non-canonical makes it apparently very hard to relate $S_{\mu}(x_0)$ to $S_{\mu}(x_1)$,  even if $x_0$ and $x_1$ are very close.  This would make it seeming impossible to prove the engulfing property for the $S_{\mu}(x_0)$ we constructed.
(let us also comment that the  $h_{\mu,x_0}(z)$ above is unbounded in general as $\mu\rightarrow 0$,  which is totally different from the real case.) However,  one surprising thing we observed is that,  if a section in the form $\{u-h\le u(x_0)+\mu\}$ happens to be close to an ellipsoid,  then the shape of this ellipsoid is ``unique" in a quatitative sense.
This key observation allows us to show the engulfing property.\\

The above discussion shows the importance to understand whether we have Pogorelov estimates for complex Monge-Ampere equations.  In particular,  we are motivated to make the following definition:
\begin{defn}
Let $\Omega\subset \mathbb{C}^n$ be a bounded domain and $\Omega'\subset \Omega$ be compactly contained in $\Omega$.   We say that $(\Omega,\Omega')$ has the Pogorelov property if:
\begin{enumerate}
\item There exists $u\in PSH(\Omega)\cap C(\bar{\Omega})$,  solving $MA(u)=1$ in $\Omega$ in the sense of Bedford-Taylor,  and $u=0$ on $\partial \Omega$,
\item The solution $u$ is $C^2$ (hence $C^{\infty}$) on $\Omega'$.
\end{enumerate}
\end{defn}
If one carefully checks the argument of the present paper,  what we really proved is the following result:
\begin{thm}
Let $u\in PSH(B_1)\cap C(B_1)$ solve $1-\eps\le MA(u)\le 1+\eps$ in $B_1$.  Assume that there exists finitely many Pogorelov pairs $(\Omega_i,\Omega_i')$,  $1\le i\le N$,  such that
\begin{enumerate}
\item Each $\Omega_i$ is of the form $\{u-h_i<c_i\}$ for some pluriharmonic function $h_i$ and $c_i\in\mathbb{R}$.
\item $\bar{B}_{\frac{1}{2}}\subset \cup_{i=1}^N\Omega_i'$.
\end{enumerate}
Then we have $u\in W^{2,p}(B_{\frac{1}{2}})$,  as long as $\eps$ is small enough,  depending only on $n$,  $p$,  the lower and upper Hessian bound and $C^3$ bound for $u_i$ on $\Omega_i'$.   Here $u_i$ is the solution to $MA(u_i)=1$ on $\Omega_i$ and $u_i=0$ on $\partial\Omega_i$.
\end{thm}
For now,  it seems mysterious to characterize when the assumptions (1) and (2) above hold.  We are only able to verify such assumptions when $u$ is close to a smooth function for the moment.  For example,  in the setting of Theorem \ref{main theorem baby},  $\Omega=\{u<0\}$,  and $(\Omega,B_{0.6})$ has the Pogorelov property,  as long as $\Omega$ is close to a unit ball,  thanks to Savin's $C^{2,\alpha}$ estimates for small perturbations.\\

To conclude the Introduction,  we will explain the organization of the rest of the paper.

In Section 2,  we include some definitions,  notations and some preliminary results we will use again and again in this paper.

In Section 3,  we show how to reduce Theorem \ref{main theorem} to Theorem \ref{main theorem baby}.  The later is a special case of the former by taking $w_0=|z|^2-1$.  Section 4-6 below are devoted to the proof of Theorem \ref{main theorem baby}.

In Section 4,  we construct sections $S_{\mu}(x_0)$ for $u$ which are of the form $\{(u-h_{\mu,x_0})(z)\le u(x_0)+\mu\}$,  where $h_{\mu,x_0}$ is pluriharmonic and $S_{\mu}(x_0)$ is close to an ellipsoid.  The second half of this section focuses on the engulfing property of sections.

In Section 5,  we prove some measure-theoretic lemmas which will be needed to estimate the ``bad" set where the second derivative loses control.
These lemmas are all standard results for balls,  but we have to adapt them to $S_{\mu}(x_0)$ we constructed in Section 4.  The engulfing property of $S_{\mu}(x_0)$ is crucially used in establishing these lemmas.

In Section 6,  we verify that the ``bad sets" fits in the assumptions of the measure theoretic lemmas in Section 5,  and obtain the power decay of the measure of the bad set.  Contrary to the real case,  we first obtain control for the mixed Hessian $u_{i\bar{j}}$,  then the full $W^{2,p}$ estimate follows from the classical $L^p$ estimate for Laplacian.

In Section 7,  we discuss some implications of Theorem \ref{main theorem}.  In particular,  we give detailed proofs for Corollary \ref{c1.2},  \ref{c1.2N} ,\ref{c1.3} and \ref{c1.4}.

\section{preliminaries}
The key result we will need again and again is the following lemma:

\begin{lem}\label{Savin estimate}
Let $u$ be a viscosity solution to $det (u_{i\bar{j}})=1$ in $B_1$. Suppose that $||u-w||_{L^{\infty}}\le \delta$, where $w$ is a smooth solution to $det (w_{i\bar{j}})=1$ in $B_1$.  If $\delta$ is small enough depending only on the smoothness of $w$ and $n$,  we have $||u||_{C^{4}(B_{0.99})}\le C$,  where $C$ has the same dependence as $\delta$.
\end{lem}
The small perturbation theorem of Savin is for more general fully nonlinear elliptic equations,  and applies to equations of the form $F(D^2u,x)=0$,  where $F(r,x):\mathcal{S}\times B_1\rightarrow \mathbb{R}$ is $C^2$ in the $x$ variable,  elliptic in the $r$ variable,  and uniformly elliptic only for $r$ in a neighborhood of 0,  with $F(0,x)=0$.  Then for any solution to $F(D^2u,x)=0$ with $||u||_{L^{\infty}}$ small enough,  we would get $C^{2,\alpha}$ estimate in the interior.

Lemma \ref{Savin estimate} follows from the general perturbation theorem of Savin in \cite{Savin} by writing the complex Monge-Ampere operator in the real form.  The details can be found in Yu Wang \cite{Wang}.  Also once we get $C^{2,\alpha}$ estimate the above equation,  it is straightforward to apply standard elliptic estimates to improve the estimate to $C^4$.
In this paper,  we will mostly apply Lemma \ref{Savin estimate} with $w=|z|^2-1$.

Another thing we will need is the following interpolation lemma:

We need the following interpolation estimate.
\begin{lem}\label{interpolation}
Let $u$ be a $C^4$ function in $B_{r_0}(0)$ with $|D^4u|_{L^{\infty}(B_{r_0})}\le C$,  and $|u|_{L^{\infty}(B_{r_0})}\le \mu$,  then for any $0<\lambda<r_0$,  one has
\begin{equation*}
\begin{split}
    &|u(0)|\le \mu,\,\,\,|D u(0)|\le C(\lambda^3+\frac{\mu}{\lambda}),\,\, \\
    &|D^2u(0)|\le C(\lambda^2+\frac{\mu}{\lambda^2}), \,\, |D^3 u(0)|\le C(\lambda +\frac{\mu}{\lambda^3}).
\end{split}
\end{equation*}
\end{lem}
\begin{proof}
From the Taylor expansion,  we find that,  for $x\in B_{r_0}$,  the following estimate holds:
\begin{equation*}
|u(x)-\sum_{|\alpha|\le 3}\frac{D^{\alpha}u(0)}{\alpha!}x^{\alpha}|\le C_nC|x|^4.
\end{equation*}
If we restrict to $|x|\le \lambda$,  we find that
\begin{equation*}
|\sum_{|\alpha|\le 3}\frac{D^{\alpha}u(0)}{\alpha!}x^{\alpha}|\le C_nC\lambda^4+\mu.
\end{equation*}
This is equivalent to:
\begin{equation*}
\sup_{|y|\le 1}|\sum_{|\alpha|\le 3}\frac{D^{\alpha}u(0)}{\alpha!}\lambda^{|\alpha|}y^{\alpha}|\le C_nC\lambda^4+\mu.
\end{equation*}
Hence we would have:
\begin{equation*}
|D^mu(0)|\lambda^m\le C_n'C\lambda^4+\mu,\,\,\,m=1,\,2,\,3.
\end{equation*}
\end{proof}

In the present paper,  we will frequently use differentiation with respect to complex variables.  Following the usual conventions,  we define:
\begin{equation}
\partial_{z_i}=\frac{1}{2}(\partial_{x_i}-\sqrt{-1}\partial_{y_i}),\,\,\,\partial_{\bar{z}_i}=\frac{1}{2}(\partial_{x_i}+\sqrt{-1}\partial_{y_i}),\,\,\,z_i=x_i+\sqrt{-1}y_i.
\end{equation}
So that we find the Laplacian operator can be written as:
\begin{equation*}
\Delta =4\sum_{i=1}^n\partial_{z_i\bar{z}_i}.
\end{equation*}
A notion we will encounter again and again is pluriharmonic function,  which we explain below.
\begin{defn}
Let $\Omega\subset\mathbb{C}^n$ be a domain.  Let $h\in C^2(\Omega)$.  We say that $h$ is pluriharmonic if $h_{z_i\bar{z}_j}(z)=0$ for all $1\le i,\,j\le n$.  
\end{defn}
Note that $h$ being pluriharmonic will imply $h$ being harmonic,  but not the other way.  
One can also see that if $h$ is the real part of a holomorphic function,  then $h$ is pluriharmonic.

Another definition we need is:
\begin{lem}
Let $T:\mathbb{C}^n\rightarrow \mathbb{C}^n$ be a $\mathbb{C}$-linear transformation,  we define $||T||$ to be the operator norm of $T$,  namely:
\begin{equation*}
||T||=\sup_{|z|\le 1}|T(z)|.
\end{equation*}
\end{lem}
In the proof of engulfing property of sections,  we will need to frequently consider dilation maps.  Hence we introduce the following definition to make the notations simpler.
\begin{defn}
Let $E\subset \mathbb{C}^n$ be a set and $x_0\in E$.  We will sometimes denote $E$ to be $E(x_0)$ to indicate it is a ``pointed set".  Let $c>0$,  we define:
\begin{equation*}
cE(x_0)=\{x_0+c(y-x_0):y\in E(x_0)\}.
\end{equation*}
Namely $cE(x_0)$ is the image of the dilation map centered at $x_0$ by a factor $c$.
\end{defn}

\section{Reduction of Theorem \ref{main theorem} to Theorem \ref{main theorem baby}}

In this section,  we will show how to use Theorem \ref{main theorem baby} to deduce Theorem \ref{main theorem}.

To see the implication in an intuitive way,  we can take any point $x_0\in B_{\frac{1}{2}}$.  After subtracting a pluriharmonic function $h_{x_0}$,  we will see that $\{w-h_{x_0}(z)\le w_0(x_0)+\mu\}$ will be close to an ellipsoid (centered at $x_0$) when $\mu$ is small enough.  The same would be true for $\{u-h_{x_0}(z)\le \mu\}$.  After normalization the ellipsoid to a unit ball,  we are in the situation of Theorem \ref{main theorem baby} and we get $u\in W^{2,p}$ in a neighborhood of $x_0$.  We will make this idea precise in the rest of this section.

Let $w_0$ be as in Theorem \ref{main theorem}.
Take any point $x_0\in B_{\frac{1}{2}}$,  we can write down the Taylor expansion of $w_0$ at $x_0$:
\begin{equation}\label{2.1}
\begin{split}
    &w_0=w_0(x_0)+Re(\sum_il_{x_0,i}(z-x_0)_i)+\sum_{i,j}a_{x_0,i\bar{j}}(z-x_0)_i\overline{(z-x_0)}_j\\
    &+Re(\sum_{i,j}b_{x_0,ij}(z-x_0)_i(z-x_0)_j)+O(|z-x_0|^3).
\end{split}
\end{equation}
Define $h_{x_0}(z)=Re(\sum_il_{x_0,i}(z-x_0)_i)+Re(\sum_{i,j}b_{x_0,ij}(z-x_0)_i(z-x_0)_j)$,  
First we want to show that if we have another function $u_0$,  such that $|u_0-w_0|\le \delta$,  then the section $\{z:u_0-h_{x_0}(z)\le u_0(x_0)+\mu\}$ will be close to an ellipsoid if $\mu$ is small,  but much larger than $\delta$.  More precisely
\begin{lem}\label{l2.1}
Let $w_0$ be as stated in Theorem \ref{main theorem}.  Namely we assume that $w_0\in C^3(B_1)$,  and $\frac{1}{C_0}I\le (w_0)_{i\bar{j}}\le C_0I$,  $|D^3w_0|\le C_0$ on $B_{0.99}$.  
Let $\delta\ge 0$ and $u_0$ is a function on $B_1$ with $|u_0-w_0|\le \delta$ on $B_{0.95}$.  Then there exists $C_1>0$ large enough and $\mu_0>0$ small enough depending only on $C_0$,  such that for all $\mu$ with $4C_1\delta\le \mu\le \mu_0$,  we have: 
\begin{equation*}
(1-C_1\gamma)E_{\mu}(x_0)\subset \{z\in B_{\frac{1}{2C^2_0}}(x_0): (u_0-h_{x_0})(z)\le u_0(x_0)+\mu\}\subset (1+C_1\gamma)E_{\mu}(x_0).
\end{equation*}
Moreover,  $(u_0-h_{x_0})(z)=u_0(x_0)+\mu$ on $\partial \{z\in B_{\frac{1}{2C^2_0}}(x_0): (u_0-h_{x_0})(z)\le u_0(x_0)+\mu\}$.
\end{lem}
Here $\gamma=\frac{\delta}{\mu}+\mu^{\frac{1}{2}}$ and $E_{\mu}(x_0)=\{z\in\mathbb{C}^n:\sum_{i,j=1}^na_{x_0,ij}(z-x_0)_i\overline{(z-x_0)_j}\le \mu\}$.
\begin{proof}
Using (\ref{2.1}),  we see that on $B_1$:

\begin{equation}\label{2.2}
-C_0|z-x_0|^3\le w_0-w_0(x_0)-h_{x_0}(z)- \sum_{i,j}a_{x_0,ij}(z-x_0)_i\overline{(z-x_0)_j}\le C_0|z-x_0|^3.
\end{equation}
Since $|u_0-w_0|\le 2\delta_0$ on $B_{0.95}$,  we see that for any $x_0\in B_{0.95}$:
\begin{equation*}
\begin{split}
-2\delta-C_0|z-x_0|^3&\le u_0-u_0(
x_0)-h_{x_0}(z)-\sum_{i,j}a_{x_0,ij}(z-x_0)_i\overline{(z-x_0)_j}\\
&\le C_0|z-x_0|^3+2\delta.
\end{split}
\end{equation*}
Let $z\in B_{\frac{1}{2C^2_0}}(x_0)$ and $(u_0-h_{x_0})(z)\le u_0(x_0)+\mu$,  we get
\begin{equation*}
-2\delta-C_0|z-x_0|^3\le \mu-\sum_{i,j}a_{x_0,ij}(z-x_0)_i\overline{(z-x_0)_j}.
\end{equation*}
Since $a_{x_0,ij}\ge \frac{1}{C_0}I$,  we get
\begin{equation*}
\frac{1}{C_0}|z-x_0|^2\le \sum_{i,j}a_{x_0,ij}(z-x_0)_i\overline{(z-x_0)_j}\le \mu+2\delta+C_0\cdot \frac{1}{2C_0^2}|z-x_0|^2.
\end{equation*}
So that
\begin{equation*}
|z-x_0|\le \sqrt{2C_0(\mu+2\delta)}\le \sqrt{3C_0\mu}.
\end{equation*}
So that
\begin{equation*}
\begin{split}
&\sum_{i,j}a_{x_0,ij}(z-x_0)_i\overline{(z-x_0)_j}\le u_0-u_0(x_0)-h_{x_0}(z)+2\delta+C_0|z-x_0|^3\\
&\le \mu+2\delta+C_0(3C_0\mu)^{\frac{3}{2}}\le (1+C_1\gamma)\mu.
\end{split}
\end{equation*}
This proves the inclusion 
\begin{equation*}
\{z\in B_{\frac{1}{2C_0^2}}(x_0):(u_0-h_{x_0})(z)\le u_0(x_0)+\mu\}\subset (1+C_1\gamma)E_{\mu}(x_0).
\end{equation*}
Now we prove the other inclusion.  Let $z\in (1-C_1\gamma)E_{\mu}(x_0)$,  which implies
\begin{equation*}
\sum_{i,j}a_{x_0,ij}(z-x_0)_i\overline{(z-x_0)_j}\le (1-C_1\gamma)^2\mu\le (1-C_1\gamma)\mu.
\end{equation*}
So that
\begin{equation*}
|z-x_0|^2\le C_0\sum_{i,j}a_{x_0,ij}(z-x_0)_i\overline{(z-x_0)_j}\le C_0(1-C_1\gamma)\mu<\frac{1}{2C_0^2},
\end{equation*}
as long as $\mu\le \mu_0$ with $\mu_0$ small enough.  Moreover
\begin{equation*}
\begin{split}
&u_0-u_0(x_0)-h_{x_0}(z)\le \sum_{i,j}a_{x_0,ij}(z-x_0)_i\overline{(z-x_0)_j}+C_0|z-x_0|^3+2\delta\\
&\le (1-C_1\gamma)\mu+C_0(3C_0\mu)^{\frac{3}{2}}+2\delta\le \mu.
\end{split}
\end{equation*}
The last inequality would hold if we take $C_1$ to be large enough depending on $C_0$.  This proves the inclusion:
\begin{equation*}
(1-C_1\gamma)E_{\mu}(x_0)\subset \{z\in B_{\frac{1}{2C_0^2}}(x_0):(u_0-h_{x_0})(z)\le u_0(x_0)+\mu\}.
\end{equation*}
\end{proof}
Now we are ready to verify the implication from Theorem \ref{main theorem baby} to \ref{main theorem}.   

Let $u$ and $w_0$ be as stated in Theorem \ref{main theorem}.
Let $\mu>0$ and $x_0\in B_{0.8}$.  Let $T_{\mu,x_0}$ be a $\mathbb{C}$-affine transformation such that $T_{\mu,x_0}(B_{\sqrt{\mu}}(0))=E_{\mu}(x_0)$.  Define 
\begin{equation}\label{2.3N}
u_{\mu,x_0}(\zeta)=\frac{1}{\mu|\det T_{\mu,x_0}|^{\frac{2}{n}}}(u-h_{x_0}-\mu)(T_{\mu,x_0}(\sqrt{\mu}\zeta)).
\end{equation}
Since $E_{\mu}(x_0)$ is defined in terms of $a_{x_0,ij}$,  with $\frac{1}{C_0}\le a_{x_0,ij}\le C_0I$,  it is easy to see that:
\begin{equation*}
||T_{\mu,x_0}||\le C_2,\,\,\,||T_{\mu,x_0}^{-1}||\le C_2,\,\,\,\frac{1}{C_2}\le |\det T_{\mu,x_0}|^2\le C_2.
\end{equation*}
Here $C_2$ is a large enough constant depending only on $C_0$ and $n$.
Define $\Omega_{\mu}=T_{\mu,x_0}^{-1}(\{z\in B_{\frac{1}{2C_0^2}}:(u-h_{x_0})(z)\le u(x_0)+\mu\})$.
Then by straightforward calculation and Lemma \ref{l2.1},  we can see the following:
\begin{lem}\label{L3.2N}
There is $\mu_0>0$ small enough depending only on $C_0$ such that for all $4C_1\delta_0\le \mu\le \mu_0$ (with $C_1>0$ being the constant given by Lemma \ref{l2.1}),  we have
\begin{enumerate}
\item $B_{1-C_1\gamma}\subset \Omega_{\mu}\subset B_{1+C_1\gamma}$,  with $\gamma=\frac{\delta_0}{\mu}+\mu^{\frac{1}{2}}$.
\item $\det(u_{\mu,x_0})_{\zeta_i\bar{\zeta}_j}=f(T_{\mu,x_0}(\sqrt{\mu}\zeta))$ in $\Omega_{\mu}$,  $u_{\mu,x_0}=0$ on $\partial\Omega_{\mu}$.
\end{enumerate}
\end{lem}
The renormalized function $u_{\mu,x_0}$ fits in the assumptions for Theorem \ref{main theorem baby} after suitably choosing the parameters,  and Theorem \ref{main theorem} follows as a direct consequence:
\begin{cor}\label{c2.3}
Theorem \ref{main theorem} holds,  if we assume Theorem \ref{main theorem baby} and $u\in C^2(B_1)$.
\end{cor}
\begin{proof}
We wish to apply Theorem \ref{main theorem baby} to each $u_{\mu,x_0}$.  In order to do so,  we just need:
\begin{equation*}
C_1\gamma=C_1(\frac{\delta_0}{\mu}+\mu^{\frac{1}{2}})\le \gamma_0(n),\,\,\,|f(T_{\mu,x_0}(\sqrt{\mu}\zeta))-1|\le \eps(n,p).
\end{equation*}
Here $\gamma_0(n)$ and $\eps(n,p)$ are the constants given by Theorem \ref{main theorem baby}.

So we could just take $\mu$ so that $2C_1\mu^{\frac{1}{2}}\le \frac{1}{2}\gamma_0(n)$ and also $\mu\le \mu_0$ (given by Lemma \ref{L3.2N}).  With this $\mu$,  we can take $\delta_0$ so that $C_1\frac{\delta_0}{\mu}\le \frac{1}{2}\gamma_0(n)$ and also that $4C_1\delta_0\le \mu$.  We fix this choice from now on. 

Since we assumed that Theorem \ref{main theorem baby} holds,  we conclude that:
\begin{equation*}
||u_{\mu,x_0}||_{W^{2,p}}(B_{\frac{1}{2}})\le C,
\end{equation*}
where $C$ is a constant depending only on $n$ and $p$.  Then using (\ref{2.3N}) we may go back to $u$ and obtain that 
\begin{equation*}
||u||_{W^{2,p}(E_{\frac{1}{2}\mu}(x_0))}\le C'.
\end{equation*}
Here $C'$ depends on $C_0$,  $p$ and $n$.  Note that $\mu$ is already chosen which depends only on $C_0$ and $n$,  and $h_{x_0}$ is defined using $w_0$,  which can be bounded in terms of $C_0$ and $n$ as well.

Note that $E_{\frac{1}{2}\mu}(x_0)$ contains $B_{r_0}(x_0)$ for some $r_0>0$ small enough (depending only on $C_0$) for any $x_0\in B_{0.8}$.  The result of Theorem \ref{main theorem} would follow right away.
\end{proof}
Next we can use an approximation argument to remove the assumption that $u\in C^2(B_1)$.
\begin{proof}
(of Theorem \ref{main theorem},  without assuming $u\in C^2(B_1)$)

First,  we can find $f_k\in C^{\infty}(B_{0.9})$ such that $f_k\rightarrow f$ in $L^2(B_{0.9})$ and $|f_k-1|\le \eps$ (since $|f-1|\le \eps$,  one can see that the standard smoothing by convolution will preserve this property).   We can also find a sequence of $g_k\in C^{\infty}(\partial B_{0.9})$,  such that $g_k\rightarrow u$ uniformly on $\partial B_{0.9}$ (since $u$ is assumed to be continuous).  

Let $v_k$ be the solution to the Dirichlet problem:
\begin{equation*}
\det(v_k)_{i\bar{j}}=f_k\,\,\,\text{ in $B_{0.9}$},\,\,\,v_k=g_k\,\,\text{ on $\partial B_{0.9}$}.
\end{equation*}
From Caffarelli-Kohn-Nirenberg-Spruck \cite{CKNS},  we know that $v_k\in C^{\infty}(\bar{B}_{0.9})$.  Also from the following Lemma \ref{local stability},  we know that $v_k\rightarrow u$ uniformly on $\bar{B}_{0.9}$.  Hence for large enough $k$,  $v_k$ will fullfil the assumption of Theorem \ref{main theorem},  and each $v_k$ is smooth.  Hence we may use Corollary \ref{c2.3} to conclude that 
\begin{equation*}
||v_k||_{W^{2,p}(B_{\frac{1}{2}})}\le C,
\end{equation*}
where $C$ depends only on $C_0$,  $n$ and $p$.  In particular,  $C$ is uniform in $k$.
Passing to the limit,  we see that $u\in W^{2,p}(B_{\frac{1}{2}})$,  with the same bound $C$.
\end{proof}
In the above,  we used the following stability estimate from Dinew-Kolodziej \cite{DK} to deduce the uniform convergence of the approximation sequence:
\begin{lem}(\cite{DK})\label{local stability}
Let $\omega_E$ be the Euclidean K\"ahler form on $\mathbb{C}^n$.  Let $q>1$. Consider $u,v\in PSH(\Omega)\cap C(\bar{\Omega})$.  Assume that
\begin{equation*}
    MA(u) =f, \,\,\, MA(v)=g .
\end{equation*}
for some $f,g \in L^q (\Omega, dV)$. Then
\begin{equation*}
    \sup_{\Omega}(v-u) \le \sup_{\partial \Omega}(v-u) +c(q,n,diam(\Omega))||f-g||^{\frac{1}{n}}_{L^q(\Omega)}.
\end{equation*}
\end{lem}

\section{construction of sections}
From now on we will focus on the proof of Theorem \ref{main theorem baby}.  Our first step is to construct sections of $u$ which are close to ellipsoids via an induction process.  Next,  we prove some fine properties of the sections which ensures that they are good differentiation basis.  
\subsection{Inductive construction of sections}
Let us summarize our construction into the following proposition:
\begin{prop}\label{p3.1}
Let $\Omega$ and $u$ be as stated in Theorem \ref{main theorem baby},  with $\gamma_0$ small enough depending only on  $n$.  Let $0<\sigma<1$ be given.  Then there exists $\eps>0$ depending only on $\sigma$ and $n$,  such that if $|f-1|\le \eps$,  the following hold:
\begin{enumerate}
\item There exists $\mu_0>0$ small enough depending only on $n$ and $\sigma$,  such that for all $x_0\in B_{0,8}$ and all $\mu\le \mu_0$,  there exists a degree 2 pluriharmonic polynomial $h_{\mu,x_0}(z)$ with $h_{\mu,x_0}(x_0)=0$,  such that
\begin{equation*}
\begin{split}
&(1-0.1\sigma)E_{\mu}(x_0)\subset S_{\mu}(x_0):=\{z\in \Omega: (u-h_{\mu,x_0})(z)\le u(x_0)+\mu\}\\
&\subset (1+0.1\sigma)E_{\mu}(x_0).
\end{split}
\end{equation*}
In the above,  $E_{\mu}(x_0)=\{z\in\mathbb{C}^n:\sum_{i,j=1}^na_{\mu,x_0,ij}(z-x_0)_i\overline{(z-x_0)_j}\le \mu\}$,  with $a_{\mu,x_0,ij}$ being positive Hermitian and $\det a_{\mu,x_0,ij}=1$.
\item There is a function $c(\sigma):\sigma\in (0,1)\rightarrow \mathbb{R}_{>0}$,  such that for any $x_0\in B_{0.8}$ and any $0<\mu_1\le \mu_2\le \frac{\mu_0}{1+c(\sigma)}$,  one has $S_{\mu_1}(x_0)\subset S_{(1+c(\sigma))\mu_2}(x_0)$.  Moreover,  $0<c(\sigma)\le C_{2,n}\sigma^{\frac{1}{2}}$ for some dimensional constant $C_{2,n}$.  
\item There is a dimensional constant $C_{3,n}>0$ such that for all $0<\mu\le \mu_0$ and any $x_0\in B_{0.8}$,  there exists a $\mathbb{C}$-linear transformation $T_{\mu,x_0}$,  such that $|\det T_{\mu,x_0}|=1$,  $T_{\mu,x_0}(B_{\sqrt{\mu}}(0))=E_{\mu}(x_0)$,  $T_{\mu_0,x_0}=id$.  
Moreover,  for any $0<\mu_1<\mu_2\le \mu_0$ and any $x_0\in B_{0.8}$:
\begin{equation*}
||T_{\mu_1,x_0}\circ T_{\mu_2,x_0}^{-1}||\le C_{3,n}(\frac{\mu_2}{\mu_1})^{\frac{C_{3,n}\sigma^{\frac{1}{2}}}{-\log(0.1\sigma)}},\,\,\,||T_{\mu_2,x_0}\circ T_{\mu_1,x_0}^{-1}||\le C_{3,n}(\frac{\mu_2}{\mu_1})^{\frac{C_{3,n}\sigma^{\frac{1}{2}}}{-\log(0.1\sigma)}}.
\end{equation*}
\end{enumerate}
\end{prop}
\begin{rem}
We will make a choice of $\sigma$ later on,  depending on the value of $p$ in the $W^{2,p}$ estimate.  (The larger $p$ is,  the smaller $\sigma$ needs to be.) So that the choice of $\eps$ eventually depends only on $p$ and $n$.
\end{rem}

We fix some $0<\sigma<1$,  and describe the construction of $S_{\mu}(x_0)$.

First we solve the following Dirichlet problem on $\Omega$.
\begin{equation}\label{3.1}
\begin{split}
   & \det((v_0)_{i\bar{j}})=1  \text{ in } \Omega \\
   & v_0=0 \text{ on } \partial \Omega.
\end{split}
\end{equation}
To start the process,  we need that $v_0$ is smooth in the interior.  This is guaranteed by the fact that $\Omega$ is close to $B_1$.  More precisely,  we have:
\begin{lem}\label{dirichlet problem estimate}
Let $\Omega\subset \mathbb{C}^n$ be a bounded domain and $B_{1-\gamma}(0)\subset \Omega\subset B_{1+\gamma}(0)$ for some $0<\gamma<1$.  Let $v_0$ be the solution to the Dirichlet problem in (\ref{3.1}),  then 
\begin{equation*}
|z|^2-1-3\gamma \le v_0\le |z|^2-1+3\gamma.
\end{equation*}
Moreover,  there exists $\gamma_n>0$ small enough,  such that if $\gamma\le \gamma_n$,  we have $v_0\in C^{4}(\bar{B}_{0.9})$ with $||v_0-(|z|^2-1)||_{C^4,  B_{0.9}}\le C$.  Here $C$ depends only on $n$.
\end{lem}
\begin{proof}
From the assumption,  we see that $1-\gamma\le |z|\le 1+\gamma$ on $\partial \Omega$,  hence
\begin{equation*}
|z|^2-1-3\gamma \le v_0\le |z|^2-1+3\gamma,\,\,\,\text{ on $\partial \Omega$.}
\end{equation*}
Note that both $|z|^2-1-3\gamma$ and $|z|^2-1+3\gamma$ satisfy $\det u_{i\bar{j}}=1$. Hence from maximum principle,  we see that 
\begin{equation*}
|z|^2-1-3\gamma\le v_0\le |z|^2-1+3\gamma\text{ in $\Omega$.}
\end{equation*}
If $\gamma$ is small enough,  then we may use Savin's estimate (Lemma \ref{Savin estimate}) to see that $v_0$ is bounded in $C^{2,\alpha}$ on $B_{0.95}$ by a dimensional constant $C$.  Then one can differentiate the $\det(v_0)_{i\bar{j}}=1$ and use classical elliptic estimates to conclude that $|v_0|_{C^4,B_{0.9}}\le C'$.
\end{proof}
As a consequence,  we get that $v_0$ is actually $C^3$ close to $|z|^2-1$,  hence convex if $\gamma$ is small enough.  More precisely:
\begin{cor}\label{c3.4}
Let $v_0$ be as in Lemma \ref{dirichlet problem estimate}.  Then for any $0<\gamma\le \gamma_n$ with $\gamma_n$ small enough,  we have:
\begin{equation*}
|D^m(v_0-(|z|^2-1))|_{B_{0.9}}\le C_n\gamma^{1-\frac{m}{4}},\,\,\,m=1,\,2,\,3.
\end{equation*}
\end{cor}
\begin{proof}
This follows from Lemma \ref{dirichlet problem estimate} and the interpolation estimates.
\end{proof}
In order to define sections for $u$,  we need to show that $u$ and $v_0$ are sufficiently close.  This is guaranteed by the following lemma: 
\begin{lem}\label{l5}
Assume that $\det u_{i\bar{j}}=f$ in $\Omega$ and $u|_{\partial\Omega}=0$.  Let $v_0$ be the solution to the Dirichlet problem (\ref{3.1}).  Assume that $1- \eps \le f\le 1+\eps$.  Then we have $(1+\eps)^{\frac{1}{n}}v_0\le u\le (1-\eps)^{\frac{1}{n}}v_0$.  In particular
\begin{equation*}
|v_0-u|\le 4\eps \text{ in $\Omega$.}
\end{equation*}
\end{lem}

\begin{proof}
Since $\det\big((1+\eps)^{\frac{1}{n}}(v_0)_{i\bar{j}}\big)=1+\eps \ge f =\det u_{i\bar{j}} \ge 1-\eps =\det\big((1-\eps)^{\frac{1}{n}} (v_0)_{i\bar{j}}\big)$ and those three functions all have the same boundary value,
we can use the maximum principle to conclude that 
\begin{equation*}
(1+\eps)^{\frac{1}{n}}v_0\le u\le (1-\eps)^{\frac{1}{n}}v_0.
\end{equation*}
So that
\begin{equation*}
u-v_0\le ((1-\eps)^{\frac{1}{n}}-1)v_0\le \frac{2}{n}\eps|v_0|\le \frac{4}{n}\eps\le 4\eps.
\end{equation*}
In the above,  we used Lemma \ref{dirichlet problem estimate} that 
$|v_0|\le 2$ (if $\gamma$ is small enough).  The lower estimate for $u-v_0$ is completely similar.
\end{proof}
To define sections for the first step,   we need an analogue of Lemma \ref{l2.1}:
\begin{lem}\label{l3.6}
Let $v\in C^3(B_1)$ and $|v_{i\bar{j}}-\delta_{ij}|\le c$,  $|D^3v|\le c$ on $B_{0.9}$ where $0<c<1$.  Let $\delta\ge 0$ and $u_0$ is a function on $B_1$ with $|u_0-v|\le \delta$ on $B_{0.95}$.  Then for small enough $c>0$ (depending only on $n$) and for any $x_0\in B_{0.8}$,  there is a degree 2 pluriharmonic polynomail $h_{v,x_0}(z)$,  such that for all $\mu$ with $4\delta\le \mu\le 0.9(0.9-|x_0|)^2$ we have:
\begin{equation*}
(1-\gamma)E_{\mu}(x_0)\subset \{z\in B_1:(u_0-h_{v,x_0})(z)\le u_0(x_0)+\mu\}\subset (1+\gamma)E_{\mu}(x_0),
\end{equation*}
and $E_{\mu}(x_0)\subset B_1$.  
Here $\gamma=\frac{2\delta}{\mu}+(3c)^{\frac{3}{2}}\mu^{\frac{1}{2}},$ and $E_{\mu}(x_0)=\{z:\sum_{i,j}v_{i\bar{j}}(x_0)(z-x_0)_i\overline{(z-x_0)_j}\le \mu\}$.
\end{lem}
\begin{proof}
The proof is very similar to Lemma \ref{l2.1}.  
First we can write down the Taylor expansion of $v$ at $x_0$: 
\begin{equation*}
\begin{split}
&v(z)=v_0(x_0)+Re\big(\sum_il_{0,x_0,i}(z-x_0)_i\big)+\sum_{i,j}v_{i\bar{j}}(x_0)(z-x_0)_i\overline{(z-x_0)_j}\\
&+Re\big(\sum_{i,j}b_{0,x_0,ij}(z-x_0)_i(z-x_0)_j\big)+O(|z-x_0|^3).
\end{split}
\end{equation*}
Define $h_{v,x_0}(z)=Re\big(\sum_il_{0,x_0,i}(z-x_0)_i+\sum_{i,j}b_{0,x_0,ij}(z-x_0)_i(z-x_0)_j\big)$,  and use the bound for $D^3v$,  $u_0-v$,  we get
\begin{equation*}
-2\delta-c|z-x_0|^3\le (u_0-h_{v,x_0})(z)-u_0(x_0)-\sum_{i,j}v_{i\bar{j}}(x_0)(z-x_0)_i\overline{(z-x_0)_j}\le c|z-x_0|^3+2\delta.
\end{equation*} 
Let $z\in B_1$ with $(u_0-h_{v,x_0})(z)\le u_0(x_0)+\mu$,  we get
\begin{equation*}
\sum_{i,j}v_{i\bar{j}}(x_0)(z-x_0)_i\overline{(z-x_0)_j}\le \mu+2\delta+c|z-x_0|^3.
\end{equation*}
By choosing $c$ small,  we may assume that $v_{i\bar{j}}(x_0)\ge \frac{1}{2}I$,  so that
\begin{equation*}
\frac{1}{2}|z-x_0|^2\le \mu+2\delta+2c|z-x_0|^2.
\end{equation*}
Hence if $4c<\frac{1}{2}$,  we get 
\begin{equation*}
|z-x_0|^2\le 2(\mu+2\delta)\le 3\mu.
\end{equation*}
Hence
\begin{equation*}
\sum_{i,j}v_{i\bar{j}}(x_0)(z-x_0)_i\overline{(z-x_0)_j}\le \mu(1+\frac{2\delta}{\mu}+3^{\frac{3}{2}}c\mu^{\frac{1}{2}}).
\end{equation*}
This proves the inclusion 
\begin{equation*}
\{z\in B_1:(u_0-h_{v,x_0})(z)\le u_0(x_0)+\mu\}\subset (1+C_1\gamma)E_{\mu}(x_0).
\end{equation*}
Next we show that $E_{\mu}(x_0)\subset B_1$.  From $|v_{i\bar{j}}-\delta_{ij}|\le c$,  we see that $E_{\mu}(x_0)\subset B_{(\frac{\mu}{1-c})^{\frac{1}{2}}}(x_0)$.  Hence we just need to make sure $(\frac{\mu}{1-c})^{\frac{1}{2}}\le 0.9-|x_0|$.  If $c$ is chosen small enough,  this is indeed true.

Now we prove the inclusion that
\begin{equation*}
(1-\gamma)E_{\mu}(x_0)\subset \{z\in B_1:(u_0-h_{v,x_0})(z)\le u_0(x_0)+\mu\}.
\end{equation*}
 Assume that $z\in (1-\gamma)E_{\mu}(x_0)$,  so that 
\begin{equation*}
\sum_{i,j}v_{i\bar{j}}(x_0)(z-x_0)_i\overline{(z-x_0)_j}\le (1-\gamma)^2\mu.
\end{equation*}
Hence $|z-x_0|\le \sqrt{2\mu}$,  using $v_{i\bar{j}}(x_0)\ge \frac{1}{2}I$.  
Therefore
\begin{equation*}
\begin{split}
&(u_0-h_{v,x_0})(z)-u_0(x_0)\le \sum_{i,j}v_{i\bar{j}}(x_0)(z-x_0)_i\overline{(z-x_0)_j}+c|z-x_0|^3+2\delta\\
&\le (1-\gamma)\mu+(2\mu)^{\frac{3}{2}}+2\delta\le \mu.
\end{split}
\end{equation*}
\end{proof}
Now we choose $\mu_0>0$ so that:
\begin{equation}\label{3.2}
3^{\frac{3}{2}}\mu_0^{\frac{1}{2}}=\min(\frac{1}{20}\sigma,\frac{1}{20}\gamma_n),  \mu_0<0.9\cdot 0.1^2.
\end{equation}
where $\gamma_n$ is the constant given by Lemma \ref{dirichlet problem estimate}.  There is no loss of generality to assume that $\sigma\le \gamma_n$ and we will assume this throughout this section.

We can apply Lemma \ref{l3.6} to $v_0$ and $u$ and construct sections $S_{\mu}(x_0)$ for all $\mu_0^2< \mu\le \mu_0$. 
\begin{cor}\label{c3.7}
Let $u$ and $\Omega$ be as stated in Theorem \ref{main theorem baby}.  Assume that $\gamma_n$ is small enough,  $\mu_0$ is chosen according to (\ref{3.2}),  and $16\eps\le \mu_0^2$,  $\frac{8\eps}{\mu_0}\le \frac{1}{20}\sigma$.   For $\mu_0^2< \mu\le \mu_0$ and $x_0\in B_{0.8}$ we define:
\begin{equation*}
S_{\mu}(x_0)=\{z\in B_1:(u-h_{v_0,x_0})(z)\le u_0(x_0)+\mu\}.
\end{equation*}
Then we have 
\begin{equation*}
(1-0.1\sigma)E_{\mu}(x_0)\subset S_{\mu}(x_0)\subset (1+0.1\sigma)E_{\mu}(x_0),
\end{equation*}
 for any $\mu_0^2< \mu\le \mu_0.$
\end{cor}
\begin{proof}
This follows directly from Lemma \ref{l3.6},  by choosing $u_0$ to be $u$,  $v$ to be $v_0$ in that lemma.  Choose $\delta=4\eps$,  and assume that $\gamma_n$ stated in Theorem \ref{main theorem baby} is small enough so as to make $|(v_0)_{i\bar{j}}-\delta_{ij}|$ and $|D^3v_0|$ small enough on $B_{0.9}$.  After these choice,  we may use Lemma \ref{l3.6} to conclude that for and $x_0\in B_{0.8}$ and $4\cdot 4\eps\le \mu\le 0.9\cdot 0.1^2$,  we can conclude:
\begin{equation*}
(1-\gamma)E_{\mu}(x_0)\subset S_{\mu}(x_0)\subset (1+\gamma)E_{\mu}(x_0).
\end{equation*}
Note that the range $16\eps \le \mu\le 0.9\cdot 0.1^2$ contains the range $\frac{1}{2}\mu_0\le \mu\le \mu_0$.  Moreover,  from the assumptions on the parameters,  we get
\begin{equation*}
\gamma\le \frac{2\cdot4\eps}{\mu_0}+3^{\frac{3}{2}}\mu_0^{\frac{1}{2}}\le \frac{\sigma}{20}+\frac{\sigma}{20}\le \frac{\sigma}{10}.
\end{equation*}
\end{proof}
Now we need to define $S_{\mu}(x_0)$ for $\mu\le  \mu_0^2$.

Let $\tilde{T}_{1,x_0}$ be a $\mathbb{C}$-affine map such that $\tilde{T}_{1,x_0}(B_{\sqrt{\mu_0}}(0))=E_{\mu_0}(x_0)$. We hope to estimate how far $\tilde{T}_{1,x_0}$ is away from identity map,  in terms of how the ellipsoid $E_{\mu_0}(x_0)$ is close to a ball.  For that we need the following lemma: 
\begin{lem}\label{l3.8}
Let $E\subset \mathbb{C}^n$ be an ellipsoid,  given by:
\begin{equation*}
E=\{z:\sum_{i,j=1}^na_{i\bar{j}}(z-x_0)_i\overline{(z-x_0)_j}\le r^2\},
\end{equation*}
with $a_{i\bar{j}}$ being positive Hermitian matrix,  $\det a_{i\bar{j}}=1$.  Then there is a $\mathbb{C}$-affine transform $T$ such that $T(B_r(0))=E$,  $\det T=1$,  and $||T-I||\le \max_{1\le i\le n}|\lambda_i^{-\frac{1}{2}}-1|$,  $||T^{-1}-1||\le \max_{1\le i\le n}|\lambda_i^{\frac{1}{2}}-1|$,  where $\lambda_i$ are eigenvalues of $a_{i\bar{j}}$.
\end{lem}
\begin{proof}
First we consider when $a_{i\bar{j}}$ is diagnal,  so that $E=\{z:\sum_{i=1}^n\lambda_i|z_i|^2\le r^2\}$,  with $\Pi_i\lambda_i=1$.  Then we define $T(w)=x_0+(\frac{w_1}{\sqrt{\lambda_1}},\cdots,\frac{w_n}{\sqrt{\lambda_n}})$.  In the general case,  we can take a unitary transformation $U$,  so that $E$ becomes diagnal under the new coordinate,  then the desired $\mathbb{C}$-affine map is given by $T=U^{-1}T'U$,  where $T'$ is the dilation map along coordinate axis.  Hence the result would follow from the diagnal case.
\end{proof}
As a consequence,  we see that:
\begin{lem}\label{l3.9}
Let $u$ and $\Omega$ be as stated in Theorem \ref{main theorem baby}.  Assume that $\gamma$ is small enough depending on $n$.  Let $\tilde{T}_{1,x_0}$ be the $\mathbb{C}$-affine map given by Lemma \ref{l3.8},  applied to $E_{\mu_0}(x_0)$.  Then for some $C_n'>0$,  we have
\begin{equation*}
||\tilde{T}_{1,x_0}-I||\le C_n'\gamma^{\frac{1}{2}},\,\,\,||\tilde{T}_{1,x_0}^{-1}-I||\le C_n'\gamma^{\frac{1}{2}}.
\end{equation*}
\end{lem}
\begin{proof}
Note that $E_{\mu}(x_0)=\{z:\sum_{i,j}(v_0)_{i\bar{j}}(x_0)(z-x_0)_i\overline{(z-x_0)_j}\le \mu\}$.  It follows from Corollary \ref{c3.4} that $|(v_0)_{i\bar{j}}-\delta_{ij}|\le C_n\gamma^{\frac{1}{2}}$,  for all $\gamma$ small enough.  So that $1-C_n\gamma^{\frac{1}{2}}\le \lambda_i\le 1+C_n\gamma^{\frac{1}{2}}$.  Therefore
\begin{equation*}
|\lambda_i^{\frac{1}{2}}-1|\le|(1+C_n\gamma^{\frac{1}{2}})^{\frac{1}{2}}-1|\le C_n'\gamma^{\frac{1}{2}}, 
\end{equation*}
which implies that $||\tilde{T}^{-1}_{1,x_0}-I||\le C_n'\gamma^{\frac{1}{2}}$.  The estimate for $||\tilde{T}_{1,x_0}-I||$ is similar.
\end{proof}

To define $S_{\mu}(x_0)$ for $\mu<\mu_0^2$,  we need to rescale the ellipsoid $E_{\mu_0}(x_0)$ to be a unit ball.  Define the change of coordinate:
\begin{equation*}
z^{(1)}=\frac{1}{\sqrt{\mu_0}}\tilde{T}_{1,x_0}^{-1}(z).
\end{equation*}
Then $\frac{1}{\sqrt{\mu_0}}\tilde{T}^{-1}_{1,x_0}(E_{\mu_0}(x_0))=B_1(0)$.  Define $\Omega_{x_0,1}=\frac{1}{\sqrt{\mu_0}}\tilde{T}_{1,x_0}^{-1}(S_{\mu_0}(x_0))$.  According to Corollary \ref{c3.7},  we have that:
\begin{equation*}
B_{1-0.1\sigma}(0)\subset \Omega_{x_0,1}\subset B_{1+0.1\sigma}(0),\,\,\,\Omega_{x_0,1}\text{ is pseudoconvex.}
\end{equation*}
Define $v_{x_0,1}$ be the solution to the following Dirichlet problem:
\begin{equation}\label{3.3}
\begin{split}
&\det(v_{x_0,1})_{i\bar{j}}=1,\,\,\,\text{ in $\Omega_{x_0,1}$},\\
&v_{x_0,1}=0\,\,\,\,\text{ on $\partial \Omega_{x_0,1}$}.
\end{split}
\end{equation}
We can normalize $u$ on $\Omega_{x_0,1}$ to be:
\begin{equation}\label{3.4}
u_{x_0,1}=\frac{1}{\mu_0}\big(u-u(x_0)-h_{v_0,x_0}-\mu_0\big)(\tilde{T}_{1,x_0}(\sqrt{\mu_0}z^{(1)})),\,\,\,z^{(1)}\in \Omega_{x_0,1}.
\end{equation}
Then 
\begin{equation*}
\begin{split}
&\det (u_{x_0,1})_{i\bar{j}}=f_{x_0,1}\text{ in $\Omega_{x_0,1}$}, f_{x_0,1}(z^{(1)})=f(\tilde{T}_{1,x_0}(\sqrt{\mu_0}z^{(1)})),\\
&u_{x_0,1}=0\text{ on $\partial \Omega_{x_0,1}$}.
\end{split}
\end{equation*}
We observe that the following holds for $u_{x_0,1}$ and $v_{x_0,1}$:
\begin{lem}\label{l3.10}
Let $v_{x_0,1}$ and $u_{x_0,1}$ be defined by (\ref{3.3}) and (\ref{3.4}) respectively.  Assume that $\sigma\le \gamma_n$,  where $\gamma_n$ is given by Lemma \ref{dirichlet problem estimate}.  Let $\mu_0$ be defined by \ref{3.2}.  Then the following hold:
\begin{enumerate}
\item $v_{x_0,1}\in C^4(B_{0.95})$,  and $|D^m(v_{x_0,1}-(|z^{(1)}|^2-1))|_{B_{0.9}}\le C_n\sigma^{1-\frac{m}{4}}$,  $m=1,2,3$,  where $C_n$ is the same $C_n$ in Corollary \ref{c3.4}.
\item $|v_{x_0,1}-u_{x_0,1}|\le 4\eps$ in $\Omega_{x_0,1}$.
\item Define 
\begin{equation*}
\tilde{h}_{x_0,1}(z^{(1)})=Re\big(\sum_i(v_{x_0,1})_i(0)z^{(1)}_i+\sum_{i,j}(v_{x_0,1})_{ij}(0)z^{(1)}_iz^{(1)}_j\big),
\end{equation*}
then for $\mu_0^2\le \mu\le \mu_0$,  we have:
\begin{equation*}
\begin{split}
&(1-0.1\sigma)\tilde{E}_{\mu}(0)\subset \{z^{(1)}\in B_1:(u_{x_0,1}-\tilde{h}_{x_0,1}-u_{x_0,1}(0))(z^{(1)})\le \mu\}\\
&\subset (1+0.1\sigma)\tilde{E}_{\mu}(0),
\end{split}
\end{equation*} 
where $\tilde{E}_{\mu}(0)=\{z^{(1)}:\sum_{i,j}(v_{x_0,1})_{i\bar{j}}(0)z^{(1)}_iz^{(1)}_{\bar{j}}\le \mu\}$.
\item Let $\tilde{T}_{2,x_0}$ be the $\mathbb{C}$-affine transform given by Lemma \ref{l3.8} normalizing $\tilde{E}_{\mu_0}(x_0)$,  then one has:
\begin{equation*}
||\tilde{T}_{2,x_0}-I||\le C_n'\sigma^{\frac{1}{2}},\,\,\,\,||\tilde{T}_{2,x_0}^{-1}-I||\le C_n'\sigma^{\frac{1}{2}},
\end{equation*}
where $C_n'$ is the constant given by Lemma \ref{l3.9}.
\end{enumerate}
\end{lem}
\begin{proof}
To prove (1),  we just note that since $v_{x_0,1}$ solves (\ref{3.3}),  with $B_{1-0.1\sigma}\subset \Omega_{x_0,1}\subset B_{1+0.1\sigma}$,  and $\sigma\le \gamma_n$,  then Lemma \ref{dirichlet problem estimate} and Corollary \ref{c3.4} can be applied to show that (1) holds,  with $\gamma$ there replaced by $\sigma$.  Item (2) above follows from Lemma \ref{l5},  since $u_{x_0,1}$ solves $\det (u_{x_0,1})_{i\bar{j}}=f_{x_0,1}$ with $|f_{x_0,1}-1|\le \eps$ (since $f$ satisfies the same).  Item (3) essentially follows from Lemma \ref{l3.6},  applied to $v=v_{x_0,1}$,  $u=u_{x_0,1}$,  $c=1$,  $\delta=4\eps$,   $x_0=0$.  Because of our choice of $\mu_0$ and $\eps$,  we would have $\gamma\le 0.1\sigma$.  Also the range for $\mu$ in Lemma \ref{l3.6} is $4\cdot 4\eps\le \mu\le 0.9^3$ (with $x_0=0$),  which contains the range $\mu_0^2\le \mu\le \mu_0$.
The proof of item (4) follows from Lemma \ref{l3.9},  because of item (3).
\end{proof}
We define 
\begin{equation*}
\tilde{S}_{1,\mu}(0)=\{z^{(1)}\in B_{0.9}\subset \Omega_{x_0,1}:(u_{x_0,1}-\tilde{h}_{x_0,1}-u_{x_0,1}(0))(z^{(1)})\le \mu\},\,\,\,\mu_0^2< \mu\le \mu_0.
\end{equation*}

We can now define $S_{\mu}(x_0)$ for $\mu_0^3< \mu\le \mu_0^2$ by transforming back to $z$ variable.  In other words,  for $\mu_0^3< \mu\le \mu_0^2$,  we define
\begin{equation*}
\begin{split}
S_{\mu}(x_0)&=\tilde{T}_{1,x_0}\big(\sqrt{\mu_0}\tilde{S}_{1,\mu_0^{-1}\mu}(0)\big)\\
&=\{z\in 0.9E_{\mu_0}(x_0)\subset S_{\mu_0}(x_0):u(z)-h_{v_0,x_0}(z)-\mu_0\tilde{h}_{x_0,1}(\frac{1}{\sqrt{\mu_0}}\tilde{T}_{1,x_0}^{-1}(z))\le u(x_0)+\mu\}.
\end{split}
\end{equation*}

Next we use an induction process to define $S_{\mu}(x_0)$,  for all $0<\mu\le \mu_0$.  Assume that for some $k_0\ge 2$,  we have defined $S_{\mu}(x_0)$ for all $\mu_0^{k_0}< \mu\le \mu_0$,  we wish to define $S_{\mu}(x_0)$ for $\mu_0^{k_0+1}< \mu\le \mu_0^{k_0}$.

We make the following $\mathbf{induction}\,\,\,\mathbf{ hypothesis}$,  stated with $k_0$:
\begin{enumerate}
\item $S_{\mu}(x_0)=\{z\in S_{\mu_0^{k-2}}(x_0):u(z)-h_{x_0,k-2}(z)\le u(x_0)+\mu\}$,  for all $\mu_0^k< \mu\le \mu_0^{k-1}$,  $2\le k\le k_0$.  Here $h_{x_0,k-2}$ is a pluriharmonic polynomial of degree 2. 
\item There exists a family of ellipsoids $E_{\mu}(x_0)$,  centered at $x_0$,  such that $(1-0.1\sigma)E_{\mu}(x_0)\subset S_{\mu}(x_0)\subset (1+0.1\sigma)E_{\mu}(x_0)$.  Moreover,  
\begin{equation*}
E_{\mu}(x_0)=\{z:\sum_{i,j}a_{i\bar{j},k-2}(z-x_0)_i\overline{(z-x_0)_j}\le \mu\},
\end{equation*}  
for all $\mu_0^k< \mu\le \mu_0^{k-1}$,  $2\le k\le k_0$,  and $\det a_{i\bar{j},k-2}=1$.
\item There exists a sequence of $\mathbb{C}$-affine coordinate change: $z^{(k)}=\frac{1}{\sqrt{\mu_0}}\tilde{T}_{k,x_0}^{-1}(z^{(k-1)})$ for $k\ge 2$,  $z^{(1)}=\frac{1}{\sqrt{\mu_0}}\tilde{T}_{x_0,1}^{-1}(z-x_0)$ with $\det \tilde{T}_{k,x_0}=1$.  $\tilde{T}_{k,x_0}$ maps $B_{\sqrt{\mu_0}}(0)$ to be the image of $E_{\mu_0^k}(x_0)$ under $z^{(k-1)}$.  
Moreover $||\tilde{T}_{k,x_0}-I||\le C'_n\sigma^{\frac{1}{2}}$,  $||\tilde{T}^{-1}_{k,x_0}-I||\le C_n'\sigma^{\frac{1}{2}}$ for all $2\le k\le k_0-1$.
\item Denote $\Omega_{x_0,k}$ to be the image of $S_{\mu_0^k}(x_0)$ under the coordinate $z^{(k)}$,  then we have $B_{1-0.1\sigma}\subset \Omega_{x_0,k}\subset B_{1+0.1\sigma},$ for all $1\le k\le k_0-1$.
\item The image of $E_{\mu}(x_0)$ under coordinate $z^{(k-1)}$ is $B_{\sqrt{\mu_0^{1-k}\mu}}(0)$,  for $\mu_0^k<\mu\le \mu_0^{k-1}$,  $2\le k\le k_0$.
\end{enumerate}
First we observe that the above inductive hypothesis indeed hold for $k_0=2$.  Indeed,  the items (1) and (3) follow from Corollary \ref{c3.7}.  Item (2) follows from item (3) of Lemma \ref{l3.10}.

Now we will construct $S_{\mu}(x_0)$ for $\mu_0^{k+1}\le \mu<\mu_0^{k_0}$,  and verify that the above inductive hypothesis continues to hold with $k_0$ replaced by $k_0+1$.  That is,  we prove:
\begin{lem}
Assume that above $\mathbf{induction}\,\,\,\mathbf{ hypothesis}$ holds with some $k_0\ge 2$,  then we can construct $S_{\mu}(x_0)$ for $\mu_0^{k_0+1}< \mu\le\mu_0^{k_0}$ which satisfies the $\mathbf{induction}\,\,\,\mathbf{ hypothesis}$ with $k_0$ replaced by $k_0+1$.
In particular,  the $\mathbf{induction}\,\,\,\mathbf{ hypothesis}$ holds for all $k_0\ge 2$.
\end{lem}
\begin{proof}
We solve the following Dirichlet problem on $\Omega_{x_0,k_0-1}$:
\begin{equation}\label{3.5New}
\det (v_{x_0,k_0-1})_{i\bar{j}}=1,\,\,\,\,\text{ in $\Omega_{x_0,k_0-1}$},\,\,\,\,\,\,v_{x_0,k_0-1}=0 \text{ on $\partial \Omega_{x_0,k_0-1}$}.
\end{equation}
 Define $T_{k,x_0}=\tilde{T}_{1,x_0}\circ \tilde{T}_{2,x_0}\cdots \tilde{T}_{k,x_0}$ so that the change of coordinate between $z^{(k)}$ and $z$ is given  by $z=x_0+T_{k,x_0}(\mu_0^{\frac{k}{2}}z^{(k)})$.
\begin{equation}\label{3.6}
u_{x_0,k_0-1}(z^{(k_0-1)})=\frac{1}{\mu_0^{k_0-1}}(u-h_{x_0,k_0-2}-u(x_0)-\mu_0^{k_0-1})(x_0+T_{k_0-1,x_0}(\mu_0^{\frac{k_0-1}{2}}z^{(k_0-1)})).
\end{equation}
Then $u_{x_0,k_0-1}$ solves:
\begin{equation*}
\begin{split}
&\det (u_{x_0,k_0-1})_{i\bar{j}}=f_{x_0,k_0-1},\,\,\text{ in $\Omega_{x_0,k_0-1}$, }\,\,f_{x_0,k_0-1}=f(x_0+T_{k_0-1,x_0}(\mu_0^{\frac{k_0-1}{2}}z^{(k_0-1)})),\\
&u_{x_0,k_0-1}=0,\,\,\,\,\text{ on $\partial \Omega_{x_0,k_0-1}$}.
\end{split}
\end{equation*}
Using the same argument as in Lemma \ref{l3.10},  (1),  we have that $v_{x_0,k_0-1}\in C^4(B_{0.95})$,  and $|D^m(v_{x_0,k_0-1}-|z^{(k_0-1)}|^2-1)|_{B_{0.9}}\le C_n\sigma^{1-\frac{m}{4}}$ for $m=1,\,2,\,3.$ This follows from our inductive hypothesis that $B_{1-0.1\sigma}\subset \Omega_{x_0,k_0-1}\subset B_{1+0.1\sigma}$,  and an application of Lemma \ref{dirichlet problem estimate} and Corollary \ref{c3.4}.  
Also we would have 
\begin{equation*}
|v_{x_0,k_0-1}-u_{x_0,k_0-1}|\le 4\eps\text{ on $\Omega_{x_0,k_0-1}$},
\end{equation*}
following the same argument as Lemma \ref{l3.10}.

Then we may consider the Taylor expansion of $v_{x_0,k_0-1}$ at $z^{(k_0-1)}=0$:
\begin{equation*}
\begin{split}
&v_{x_0,k_0-1}(z^{(k_0-1)})=v_{x_0,k_0-1}(0)+Re\big(\sum_il_iz^{(k_0-1)}_i\big)+\sum_{i,j}a_{i\bar{j}}z_i^{(k_0-1)}\bar{z}_j^{(k_0-1)}\\
&+Re\big(\sum_{i,j}b_{ij}z_i^{(k_0-1)}z_j^{(k_0-1)}\big)+O(|z^{(k_0-1)}|^3).
\end{split}
\end{equation*}
Define 
\begin{equation}\label{3.7NN}
\tilde{h}_{x_0,k_0-1}(z^{(k_0-1)})=Re(\sum_il_iz_i^{(k_0-1)})+Re\big(\sum_{i,j}b_{ij}z_i^{(k_0-1)}z_j^{(k_0-1)}\big).
\end{equation}
Then the argument for part (3) of Lemma \ref{l3.10} shows that:
\begin{equation}\label{3.8}
\begin{split}
&(1-0.1\sigma)\tilde{E}_{\mu}(0)\\
&\subset \tilde{S}_{k_0-1,\mu}:=\{z^{(k_0-1)}\in \Omega_{x_0,k_0-1}:(u_{x_0,k_0-1}-\tilde{h}_{x_0,k_0-1}-u_{x_0,k_0-1}(0))(z^{(k_0-1)})\le \mu\}\\
&\subset (1+0.1\sigma)\tilde{E}_{\mu}(0),
\end{split}
\end{equation}
for any $\mu_0^2< \mu\le  \mu_0$,  where $\tilde{E}_{\mu}(0)=\{z^{(k_0-1)}:\sum_{i,j}(v_{x_0,k_0-1})_{i\bar{j}}(0)z_i^{(k_0-1)}\bar{z}_j^{(k_0-1)}\le \mu\}$.
Now we define $\tilde{T}_{k_0,x_0}$ be the $\mathbb{C}$-linear transforma given by Lemma \ref{l3.8} normalizing $\tilde{E}_{\mu_0}(0)$ above.
Then for $\mu_0^{k_0+1}< \mu\le \mu_0^{k_0}$,  we define
\begin{equation*}
S_{\mu}(x_0)=x_0+\tilde{T}_{k_0,x_0}(\mu_0^{\frac{k_0-1}{2}}\tilde{S}_{k_0-1,\mu_0^{-(k_0-1)}\mu}).
\end{equation*}
Similarly,  we define 
\begin{equation*}
E_{\mu}(x_0)=x_0+\tilde{T}_{k_0,x_0}(\mu_0^{\frac{k_0-1}{2}}\tilde{E}_{\mu_0^{-(k_0-1)}\mu}).
\end{equation*}
Using (\ref{3.6}),  one find that,  for $\mu_0^{k_0+1}\le \mu<\mu_0^{k_0}$: 
\begin{equation*}
S_{\mu}(x_0)=\{z\in 0.9E_{\mu_0^{k_0-1}}(x_0):u(z)-h_{x_0,k_0-2}(z)-\mu_0^{k_0-1}\tilde{h}_{x_0,k_0-1}(\mu_0^{-\frac{k_0-1}{2}}T_{k_0-1,x_0}^{-1}(z-x_0))\le \mu\}.
\end{equation*}
Define
\begin{equation*}
h_{x_0,k_0-1}(z)=h_{x_0,k_0-2}(z)+\mu_0^{k_0-1}\tilde{h}_{x_0,k_0-1}(\mu_0^{-\frac{k_0-1}{2}}T_{k_0-1,x_0}^{-1}(z-x_0)).
\end{equation*}
In view of (\ref{3.7NN}),  as well as the induction hypothesis for $h_{x_0,k_0-2}(z)$,  we see that $h_{x_0,k_0-1}(z)$ is a pluriharmonic polynomial of degree 2,  and $h_{x_0,k_0-1}(x_0)=0$.  This proves the induction hypothesis,  part (1),  for $k=k_0+1$.

Part (2) simply follows from (\ref{3.8}) and the fact that $\det(v_{x_0,k_0-1})_{i\bar{j}}(0)=1$,  as well as $\det T_{k,x_0}=1$.

Since $|(v_{x_0,k_0-1})_{i\bar{j}}(0)-\delta_{ij}|\le C_n\sigma^{\frac{1}{2}}$,  the argument in Lemma \ref{l3.10},  part (4) shows that
\begin{equation*}
||\tilde{T}_{k_0,x_0}-I||\le C_n'\sigma^{\frac{1}{2}},\,\,\,\,||\tilde{T}_{k_0,x_0}^{-1}-I||\le C_n'\sigma^{\frac{1}{2}}.
\end{equation*}
Using $\tilde{T}_{k_0,x_0}$,  we can define a change of coordinates $z^{(k_0)}=\frac{1}{\sqrt{\mu_0}}\tilde{T}_{k_0,x_0}^{-1}(z^{(k_0-1)})$.
This proves part (2) of the Inductive hypothesis with $k=k_0$.  
Part (3) follows from (\ref{3.8}) and our choice of $\tilde{T}_{k_0,x_0}$ above.
\end{proof}

Now we constructed $S_{\mu}(x_0)$ for $0<\mu\le \mu_0$,  let us verify Proposition \ref{p3.1}.

Part (1) has already been proved by the above argument.  Indeed,  we define $h_{\mu,x_0}(z)=h_{x_0,k-1}(z)$ for $\mu_0^k< \mu\le \mu_0^{k-1}$.  
Indeed,  part (1) of Proposition \ref{p3.1} follows from part (1) and (2) of $\mathbf{induction}\,\,\,\mathbf{ hypothesis}$.
Now we verify part (2) of Proposition \ref{p3.1}.  First,  we make the following observation out of the above inductive process.
\begin{lem}\label{l3.12}
Let $0<\mu_1<\mu_2$,  then the following hold:
\begin{enumerate}
\item If there is some $k\ge 2$ such that $\mu_0^k< \mu_1<\mu_2\le \mu_0^{k-1}$,  then $S_{\mu_1}(x_0)\subset S_{\mu_2}(x_0)$,  $E_{\mu_1}(x_0)\subset E_{\mu_2}(x_0)$,
\item If there is some $k\ge 1$ such that $\frac{1}{2}\mu_0^k< \mu_2$ and $\mu_1\le \mu_0^{k+1}$,  then $S_{\mu_1}(x_0)\subset S_{\mu_2}(x_0)$,  $E_{\mu_1}(x_0)\subset E_{\mu_2}(x_0)$.
\end{enumerate}
\end{lem}
\begin{proof}
Part (1) above is obvious,  due to part (1) and (2) of $\mathbf{induction}\,\,\,\mathbf{ hypothesis}$.  

Part (2) requires more work,  due to that $S_{\mu}(x_0)$ and $E_{\mu}(x_0)$ are discontinuous in $\mu$ for $\mu\rightarrow \mu_0^k+$ and $\mu\rightarrow \mu_0^k-$.  
First,  without loss of generality we may assume that $\mu_0^{k+2}<\mu_1\le \mu_0^{k+1}$. 

If $\mu_0^k<\mu_2\le \mu_0^{k-1}$ (with $k\ge 2$ in this case),  
we know from part (5) of $\mathbf{induction}\,\,\,\mathbf{ hypothesis}$ that under $z^{(k-1)}$,  $E_{\mu_2}(x_0)$ is given by $B_{\sqrt{\mu_0^{1-k}\mu_2}}(0)$.  On the other hand,  the image of $E_{\mu_1}(x_0)$ under $z^{(k+1)}$ is given by   $B_{\sqrt{\mu_0^{-1-k}\mu_1}}(0)$.  Hence if we recall the transition formula given by part (3) of $\mathbf{induction}\,\,\,\mathbf{ hypothesis}$,  we see that:
\begin{equation}\label{3.9}
\begin{split}
&\text{$E_{\mu_1}(x_0)$ under $z^{(k-1)}$}=\mu_0\tilde{T}_{k+1,x_0}\circ \tilde{T}_{k,x_0}(B_{\sqrt{\mu_0^{-1-k}\mu_1}}(0))\\
&\subset (1+C_n'\sigma^{\frac{1}{2}})^2B_{\sqrt{\mu_0^{1-k}\mu_1}}(0).
\end{split}
\end{equation}
In the first inclusion above,  we used that $||\tilde{T}_{i,x_0}||\le 1+C_n'\sigma^{\frac{1}{2}}$.
Therefore
\begin{equation*}
\text{image of $S_{\mu_1}(x_1)$ under $z^{(k-1)}$}\subset B_{(1+0.1\sigma)(1+C_n'\sigma^{\frac{1}{2}})^2\sqrt{\mu_0^{1-k}\mu_1}}(0)
\end{equation*}
On the other hand,
\begin{equation*}
\text{image of $S_{\mu_2}(x_2)$ under $z^{(k-1)}$}\supset (1-0.1\sigma)B_{\sqrt{\mu_2\mu_0^{1-k}}}(0)
\end{equation*}

We will be able to show $S_{\mu_1}(x_0)\subset S_{\mu_2}(x_0)$ if we can ensure:
\begin{equation*}
(1+C_n'\sigma^{\frac{1}{2}})^4(1+0.1\sigma)^2\mu_1\le (1-0.1\sigma)^2\mu_2.
\end{equation*}
This can be guaranteed if we take $\mu_0$ small enough so that $(1+C_n'\sigma^{\frac{1}{2}})^4\frac{(1+0.1\sigma)^2}{(1-0.1\sigma)^2}\le \mu_0^{-1}$,  since $\mu_2\ge \mu_0\mu_1$.

The other case is when $\frac{1}{2}\mu_0^k<\mu_2\le \mu_0^k$.  The calculation in this case is similar to the case when $\mu_0^k<\mu_2\le \mu_0^{k-1}$,  except that we need to use the coordinate $z^{(k)}$,  and we may conclude:
\begin{equation*}
\text{image of $S_{\mu_2}(x_0)$ under $z^{(k)}$}\supset (1-0.1\sigma)B_{\sqrt{\mu_2\mu_0^{-k}}}(0),
\end{equation*}
On the other hand,  
\begin{equation*}
\text{image of $S_{\mu_1}(x_0)$ under $z^{(k+1)}$}\subset (1+0.1\sigma)B_{\sqrt{\mu_1\mu_0^{-1-k}}}(0).
\end{equation*}
Hence,  using the transition between $z^{(k)}$ and $z^{(k+1)}$:
\begin{equation*}
\text{image of $S_{\mu_1}(x_0)$ under $z^{(k)}$}\subset \sqrt{\mu_0}(1+C_n'\sigma^{\frac{1}{2}})(1+0.1\sigma)B_{\sqrt{\mu_1\mu_0^{-1-k}}}(0).
\end{equation*}
We will have the inclusion as long as we can make sure:
\begin{equation*}
(1-0.1\sigma)^2\mu_2\ge (1+C_n'\sigma^{\frac{1}{2}})^2(1+0.1\sigma)^2\mu_1.
\end{equation*}
We will still have this since $\mu_2\ge \frac{1}{2}\mu_1\mu_0$ and we can take $\mu_0$ small enough.
\end{proof}

With the help of the previous lemma,  we are ready to prove the almost monotonicity of sections claimed in part (2):
\begin{cor}
Let $c(\sigma)=\frac{(1+0.1\sigma)^2}{(1-0.1\sigma)^2}(1+C_n'\sigma^{\frac{1}{2}})^2-1$ with $C_n'$ given by Lemma \ref{l3.10}.  

 Then for all $0<\mu_1\le \mu_2\le \frac{\mu_0}{1+c(\sigma)}$ and any $x_0\in B_{0.8}$,  
\begin{equation*}
S_{\mu_1}(x_0)\subset S_{(1+c(\sigma))\mu_2}(x_0).
\end{equation*}
\end{cor}
\begin{proof}
Denote $c(\sigma)=\frac{(1+0.1\sigma)^2}{(1-0.1\sigma)^2}(1+C_n'\sigma^{\frac{1}{2}})^2-
1$.
Let $k\ge 1$ be such that $\mu_0^{k+1}<\mu_1\le \mu_0^k$ and $\mu_2>\mu_1$.  
There are several cases to consider:

Case 1: $\mu_0^{k+1}<\mu_1<\mu_2\le \frac{\mu_0^k}{1+c(\sigma)}$.
Then from Lemma \ref{l3.12},  we know that $S_{\mu_1}(x_0)\subset S_{(1+c(\sigma))\mu_2}(x_0)$.

Case 2: $\frac{\mu_0^k}{1+c(\sigma)}<\mu_2\le \frac{\mu_0^{k-1}}{1+c(\sigma)}$ $(k\ge 2$ for this case).

First,  under the coordinate $z^{(k-1)}$,  we have the following inclusions:
\begin{equation}\label{4.12NNN}
\begin{split}
&\text{image of $S_{(1+c(\sigma))\mu_2}(x_0)$ under $z^{(k-1)}$}\supset B_{(1-0.1\sigma)\sqrt{(1+c(\sigma))\mu_2\mu_0^{1-k}}}(0)\\
&=B_{(1+0.1\sigma)(1+C_n'\sigma^{\frac{1}{2}})\sqrt{\mu_2\mu_0^{1-k}}}(0).
\end{split}
\end{equation}
On the other hand,  if we consider the image of $S_{\mu_1}(x_0)$ under $z^{(k)}$,  we have
\begin{equation*}
\text{image of $S_{\mu_1}(x_0)$ under $z^{(k)}$}\subset B_{(1+0.1\sigma)\sqrt{\mu_1\mu_0^{-k}}}(0)
\end{equation*}
Then we use the transition between $z^{(k)}$ and $z^{(k-1)}$: $z^{(k)}=\frac{1}{\sqrt{\mu_0}}\tilde{T}_{k,x_0}(z^{(k-1)})$ to get:
\begin{equation}\label{4.13NNN}
\text{Image of $S_{\mu_1}(x_0)$ under $z^{(k-1)}$}\subset \sqrt{\mu_0}\tilde{T}_{k,x_0}^{-1}(B_{(1+0.1\sigma)\sqrt{\mu_1\mu_0^{-k}}}(0))\subset B_{(1+0.1\sigma)(1+C_n'\sigma^{\frac{1}{2}})\sqrt{\mu_1\mu_0^{1-k}}}(0).
\end{equation}
Combining (\ref{4.12NNN}) and (\ref{4.13NNN}),  we get that:
\begin{equation*}
S_{\mu_1}(x_0)\subset S_{(1+c(\sigma))\mu_2}(x_0).
\end{equation*}

Case 3: $\mu_2\ge \frac{\mu_0^{k-1}}{1+c(\sigma)}$ ($k\ge 2$ for this case).  

Without loss of generality,  we may assume $\sigma$ small enough so that $1+c(\sigma)<2$,  then the conclusion would follow from Lemma \ref{l3.12}.  Since $\mu_1\le \mu_0^k$ but $\mu_2\ge \frac{1}{2}\mu_0^{k-1}$.
\end{proof}

Now we verify part (3) of Proposition \ref{p3.1}.
\begin{lem}
Define $T_{\mu,x_0}=T_{k,x_0}:=\tilde{T}_{1,x_0}\circ \tilde{T}_{2,x_0}\cdots \tilde{T}_{k,x_0}$ for $\mu_0^{k+1}<\mu\le \mu_0^k$,  $k\ge 1$.  Then for $0<\mu_1<\mu_2\le \mu_0$,  we have:
\begin{equation*}
||T_{\mu_1,x_0}^{-1}\circ T_{\mu_2,x_0}||,\,\,\,||T_{\mu_2,x_0}^{-1}\circ T_{\mu_1,x_0}||\le C_{3,n}\big(\frac{\mu_2}{\mu_1}\big)^{\frac{C_{3,n}\sigma^{\frac{1}{2}}}{\log(0.1\sigma)}}.
\end{equation*}
Here $C_{3,n}$ is some dimensional constant.
\end{lem}
\begin{proof}
First we find $k_1\ge k_2$ such that $\mu_0^{k_1+1}<\mu_1\le \mu_0^{k_1}$,  $\mu_0^{k_2+1}<\mu_2\le \mu_0^{k_2}$.
Then we have:
\begin{equation*}
T_{\mu_1,x_0}^{-1}\circ T_{\mu_2,x_0}=\tilde{T}_{k_1,x_0}^{-1}\circ \tilde{T}^{-1}_{k_1-1,x_0}\cdots\tilde{T}_{k_2+1,x_0}^{-1},\,\,\,T_{\mu_2,x_0}^{-1}\circ T_{\mu_1,x_0}=\tilde{T}_{k_2+1,x_0}\circ \tilde{T}_{k_2+1,x_0}\cdots \tilde{T}_{k_1,x_0}.
\end{equation*}
Then we may use part (3) of $\mathbf{induction}\,\,\,\mathbf{ hypothesis}$ that:
\begin{equation*}
||T_{\mu_1,x_0}^{-1}\circ T_{\mu_2,x_0}||\le \Pi_{k=k_2}^{k_1-1}||T_{k,x_0}^{-1}||\le (1+C_n'\sigma^{\frac{1}{2}})^{k_1-k_2}.
\end{equation*}
On the other hand,  we easily have the bound $\mu_0^{k_2-k_1+1}\le \frac{\mu_2}{\mu_1}$,  hence
\begin{equation*}
||T_{\mu_1,x_0}^{-1}\circ T_{\mu_2,x_0}||\le (1+C_n'\sigma^{\frac{1}{2}})^{-\frac{\log(\mu_2\mu_1^{-1})}{-\log(\mu_0)}+1}.
\end{equation*}
Recall our choice of $\mu_0$ made in (\ref{3.2}): $3^{\frac{3}{2}}\mu_0^{\frac{1}{2}}=\frac{\sigma}{20}$.  Then the claimed estimate follows easily.  The same computation works also for $T_{\mu_2,x_0}^{-1}\circ T_{\mu_1,x_0}$.
\end{proof}

As a direct consequence,  the diameter of $S_{\mu}(x_0)$ should go to zero as $\mu\rightarrow 0$.  Namely:
\begin{cor}
Assume that $\sigma$ is small enough (depending on $n$).  Then the diameter of $S_{\mu}(x_0)$ goes to zero as $\mu\rightarrow 0$.  This convergence is uniform for $x_0\in B_{0.8}$.
\end{cor}
\begin{proof}
Since $T_{\mu,x_0}=id$ for $\mu_0^2<\mu\le \mu_0$,  we see that,  for any $\mu\le \mu_0^2$,  
\begin{equation*}
||T_{\mu,x_0}^{-1}||,\,\,\,||T_{\mu,x_0}||\le C_{3,n}(\frac{\mu}{\mu_0})^{\frac{C_{3,n}\sigma^{\frac{1}{2}}}{\log(0.1\sigma)}}.
\end{equation*}
On the other hand $T_{\mu,x_0}(B_{\sqrt{\mu}}(0))=E_{\mu}(x_0)-x_0$,  we see that 
\begin{equation*}
diam\,E_{\mu}(x_0)\le 2\mu^{\frac{1}{2}}||T_{\mu,x_0}||\le 2C_{3,n}\mu^{\frac{1}{2}}\cdot (\frac{\mu}{\mu_0})^{\frac{C_{3,n}\sigma^{\frac{1}{2}}}{\log(0.1\sigma)}}.
\end{equation*}
If $\sigma$ is small enough so that $\frac{1}{2}+\frac{C_{3,n}\sigma^{\frac{1}{2}}}{\log(0.1\sigma)}>0$,  then the right hand side will go to zero as $\mu\rightarrow 0$.
The result follows from that $S_{\mu}(x_0)\subset (1+0.1\sigma)E_{\mu}(x_0)$.
\end{proof}

\subsection{Further properties of sections}
Our intention will be to use the sections given by Proposition \ref{p3.1} to replace the role of balls in the uniformly elliptic case.  The most crucial property that we need is the following ``engulfing property" of sections,  formulated below:
\begin{prop}\label{p3.15}
Assume that $x_1,\,x_2\in B_{0.8}$,  $0<\mu_1,\,\mu_2\le \mu_0$ and $\mu_1\le 4\mu_2$.  Let $\sigma>0$ be small enough (depending only on dimension).  Assume also that $S_{\mu_1}(x_1)\cap S_{\mu_2}(x_2)\neq \emptyset$,  then $S_{\mu_1}(x_1)\subset 10S_{\mu_2}(x_2)$.
\end{prop}
In Caffarelli's proof for $W^{2,p}$ estimate in the real case,  we also need this ``engulfing property",  but this property is not a problem in the real case.  The essential point is the ``invariance of sections under linear transformations".  To be more  clear,  in the real case,  the sections are simply defined as: (with $u$ being strictly convex function)
\begin{equation*}
S_{\mu}(x_0)=\{x:u(x)\le u(x_0)+\nabla u(x_0)\cdot (x-x_0)+\mu\},\,\,\,\,\mu>0.
\end{equation*}
Now we define $v(y)=\frac{1}{r^2}u(x_0'+rTy)$ where $T$ is linear with $\det T=1$.  Such a transformation would preserve the Monge-Ampere equation.  \sloppy Under the change of coordinates $x=x_0'+rTy$,   $S_{\mu}(x_0)$ will be transformed to a section of $v$ centered at $y_0:=T^{-1}(\frac{1}{r}(x_0-x_0'))$,  with height $\frac{\mu}{r^2}$.

We no longer have this property in the complex case.  Indeed,  if you do a similar change of coordinates (now with $T$ being $\mathbb{C}$-linear),  and you do the same construction described in Proposition \ref{p3.1} for the function $v$ in the variable $y$,  then transform back to $x$,  you will get different sections than the direct construction in the original $x$ coordinates.  The two definitions will differ by an addition of a pluriharmonic function.

What saves us is the following ``uniqueness" property,  which shows that an addition of a pluriharmonic function will not affect the sections we get,  as long as they are close to ellipsoids.
\begin{lem}\label{key lemma}
Let $u$ be a function defined on an open set $U\subset \mathbb{C}^n$ and let $h(z)$ be pluriharmonic function on $U$ such that $h(0)=0$.  Let $0<\gamma<1$ and $\mu>0$ be such that:
\begin{equation}\label{key lemma assumption}
\begin{split}
&B_{(1-\gamma)\sqrt{\mu}}(0)\subset \{u\le u(0) +\mu\}\subset B_{(1+\gamma)\sqrt{\mu}}(0)\subset U,\\
&(1-\gamma)E_{\mu}(0)\subset\{u\le h + u(0)+\mu\}\subset (1+\gamma)E_{\mu}(0)\subset U.
\end{split}
\end{equation}
In the above,  $E_{\mu}(0)=\{z\in\mathbb{C}^n:\sum_{i,j=1}^na_{i\bar{j}}z_i\bar{z}_j\le \mu\}$ with $a_{i\bar{j}}$ positive Hermitian.  Then we have:
\begin{equation}\label{3.11}
B_{(1-\gamma)\sqrt{\mu}}(0)\subset (1+\gamma)E_{\mu}(0),\,\,\,\,(1-\gamma)E_{\mu}(0)\subset B_{(1+\gamma)\sqrt{\mu}}(0).
\end{equation}
\end{lem}
\begin{proof}
First,  by considering a unitary transformation if necessary,  we may assume that $E_{\mu}(0)=\{z\in\mathbb{C}^n:\sum_{i=1}^n\lambda_i|z_i|^2\le \mu\}$ with $0<\lambda_1\le\lambda_2\cdots
\le\lambda_n$.  
Also it will suffice to prove one of the two inclusions in (\ref{3.11}),  say $B_{(1-\gamma)\sqrt{\mu}}(0)\subset (1+\gamma)E_{\mu}(0)$.  To prove the other inclusion,  we may consider a change of coordinates: $w_i=\sqrt{\lambda_i}z_i$,  so that $E_{\mu}(0)$ becomes $B_{\sqrt{\mu}}(0)$,  and we use $v(z):=u(z)-h(z)$ to replace $u$.  Then the second inclusion would follow from the first.

We wish to argue by contradiction and assume that $B_{(1-\gamma)\sqrt{\mu}}(0)$ is not contained in $(1+\gamma)E_{\mu}(0)$,  then we must have:
\begin{equation*}
B_{(1-\gamma)\sqrt{\mu}}(0)\cap \partial \big((1+\gamma)E_{\mu}(0)\big)\neq \emptyset.
\end{equation*}

Note that $u\ge h(z)+ u(0)+\mu$ on $(1+\gamma)\partial E_{\mu}(0)$ and $u\le  u(0)+ \mu$ in $B_{(1-\gamma)\sqrt{\mu}}$.
Therefore,  $h(z)\le 0$ on $B_{(1-\gamma)\sqrt{\mu}}\cap (1+\gamma)\partial E_{\mu}$.  We will show that this hypersurface, if nonempty,  actually bounds a nontrivial region,  so that we get $h\le 0$ in a neighborhood of 0. Since $h$ is a pluriharmonic function and $h(0)=0$, we can use the strong maximum principle to get $h\equiv0$.  Hence (\ref{key lemma assumption}) would give us $B_{(1-\gamma)\sqrt{\mu}}(0)\subset \{u\le h+\mu\}\subset (1+\gamma)E_{\mu}(0)$,  contrary to what we assume above.

First we present the argument when $n=2$.   We want to show that if $\lambda_2>\frac{(1+\gamma)^2}{(1-\gamma)^2}$,  then $h\equiv 0$,  which would contradict (\ref{key lemma assumption}).  On the other hand,  if $\lambda_2\le \frac{(1+\gamma)^2}{(1-\gamma)^2}$,  then we would have $B_{(1-\gamma)\sqrt{\mu}}\subset (1+\gamma)E_{\mu}$,  which is another contradiction.

To see that $\lambda_2> \frac{(1+\gamma)^2}{(1-\gamma)^2}$ implies $h\equiv 0$,  we fix some $z_{1,*}$ and consider the cross section between $(z_{1,*},z_2)$ and $B_{(1-\gamma)\sqrt{\mu}}\cap (1+\gamma)E_{\mu}$ (viewed as a subset in $\mathbb{C}$ for $z_2$).  They are given by:
\begin{equation}\label{3.12NN}
\begin{split}
&|z_{1,*}|^2+|z_2|^2\le (1-\gamma)^2\mu,\\
&\lambda_1|z_{1,*}|^2+\lambda_2|z_2|^2\le (1+\gamma)^2\mu.
\end{split}
\end{equation}
We want to argue that,  if $\lambda_2>\frac{(1+\gamma)^2}{(1-\gamma)^2}$,  then the boundary of the cross section will be on $(1+\gamma)\partial E_{\mu}\cap B_{(1-\gamma)\sqrt{\mu}}$,   for all $z_{1,*}$ close to zero.  Since $h\le 0$ on the boundary of cross section,  we would have $h(z_{1,*},z_2)\le 0$ in the interior of the section $(z_{1,*},z_2)$ since $h$ is pluriharmonic.  This is true for all $z_{1,*}$ close to 0.  Hence $h\ge 0$ in a neighborhood of 0 and we can conclude from strong maximum principle that $h\equiv 0$ since $h(0)=0$.

The boundary of the cross section is on $(1+\gamma)\partial E_{\mu}\cap B_{(1-\gamma)\sqrt{\mu}}$ if and only if the second inequality in (\ref{3.12NN}) implies the first (with $z_{1,*}$ fixed).  
The first inequality is equivalent to:
\begin{equation*}
|z_2|^2\le (1-\gamma)^2\mu-|z_{1,*}|^2,
\end{equation*}
whereas the second inequality is equivalent to:

\begin{equation*}
|z_2|^2\le \frac{(1+\gamma)^2}{\lambda_2}\mu-\frac{\lambda_1}{\lambda_2}|z_{1,*}|^2.
\end{equation*}
In order for the second inequality to be stronger than the first,  we need that:
\begin{equation}\label{3.13}
\frac{(1+\gamma)^2}{\lambda_2}\mu-\frac{\lambda_1}{\lambda_2}|z_{1,*}|^2\le (1-\gamma)^2\mu-|z_{1,*}|^2.
\end{equation}
If $\lambda_2>\frac{(1+\gamma)^2}{(1-\gamma)^2}$,  the inequality in (\ref{3.13}) is strict with $z_{1,*}=0$.  Hence (\ref{3.13}) will hold for $z_{1,*}$ close enough to zero. 
This proves our earlier claim about the boundary of the cross section and finishes the argument that $\lambda_2>\frac{(1+\gamma)^2}{(1-\gamma)^2}$ implies $h\equiv 0$.  If $\lambda_1\le \lambda_2\le \frac{(1+\gamma)^2}{(1-\gamma)^2}$,  then $B_{(1-\gamma)\sqrt{\mu}}\subset (1+\gamma)E_{\mu}$.  Indeed,  $B_{(1-\gamma)\sqrt{\mu}}$ is given by:
\begin{equation*}
|z_1|^2+|z_2|^2\le (1-\gamma)^2\mu,  
\end{equation*}
but then $\lambda_1|z_1|^2+\lambda_2|z_2|^2\le \frac{(1+\gamma)^2}{(1-\gamma)^2}\cdot (1-\gamma)^2\mu=(1+\gamma)^2\mu$,  so that $(z_1,z_2)\in (1+\gamma)E_{\mu}$.

Now we look at general $n$.  As in $n=2$,  $\lambda_n\le \frac{(1+\gamma)^2}{(1-\gamma)^2}$ will immediately imply that $B_{(1-\gamma)\sqrt{\mu}}\subset (1+\gamma)E_{\mu}$ (since $\lambda_n$ is the largest eigenvalue),  immediately giving what we want to prove.  We just need to show that $\lambda_n>\frac{(1+\gamma)^2}{(1-\gamma)^2}$ implies $h\equiv 0$,  and then (\ref{key lemma assumption}) will give the result. 

For this we consider the cross section between $(z_{1,*},\cdots,  z_{n-1,*},z_n)$ and $B_{(1-\gamma)\sqrt{\mu}}\cap(1+\gamma)E_{\mu}$,  and it is given by:
\begin{equation}\label{3.14N}
\begin{split}
&\sum_{i=1}^{n-1}|z_{i,*}|^2+|z_n|^2\le (1-\gamma)^2\mu,\\
&\sum_{i=1}^{n-1}\lambda_i|z_{i,*}|^2+\lambda_n|z_n|^2\le (1+\gamma)^2\mu.
\end{split}
\end{equation} 
The boundary of the cross section is on $(1+\gamma)\partial E_{\mu}\cap B_{(1-\gamma)\sqrt{\mu}}$ iff the second inequality in (\ref{3.14N}) is stronger than the first.  This would mean:
\begin{equation}\label{3.15N}
\frac{1}{\lambda_n}\big((1+\gamma)^2\mu-\sum_{i=1}^{n-1}\lambda_i|z_{i,*}|^2\big)\le (1-\gamma)^2\mu-\sum_{i=1}^{n-1}|z_{i,*}|^2.
\end{equation}
Note that if $\lambda_n>\frac{(1+\gamma)^2}{(1-\gamma)^2}$,  the inequality in (\ref{3.15N}) is strict with $z_{i,*}=0$,  $1\le i\le n-1$.  Hence (\ref{3.15N}) will continue to hold for $(z_{1,*},\cdots,z_{n-1,*})$ close to 0.  Hence we would have $h\le 0$ in a neighborhood of 0,  and strong maximum principle would give $h\equiv 0$.

\end{proof}

As a consequence,  we deduce that:
\begin{cor}\label{c3.17}
Let $u$ be a function defined on an open set $U\subset \mathbb{C}^n$ with $0\in U$.  Let $h_1(z),\,h_2(z)$ be pluriharmonic functions on $U$ such that $h_1(0)=h_2(0)=0$.  Let $0<\gamma<1$ and $\mu>0$ be such that:
\begin{equation*}
(1-\gamma)E_{p,\mu}(0) \subset \{u\le h_p+u(0)+\mu\}\subset (1+\gamma)E_{p,\mu}(0)\subset U,\,\,\,\,p=1,\,2.
\end{equation*}
In the above,  $E_{p,\mu}(0)=\{z\in\mathbb{C}^n:\sum_{i,j}a_{p,i\bar{j}}z_i\bar{z}_j\le \mu\}$ with $a_{p,i\bar{j}}$ being positive Hermitian and $\det a_{p,i\bar{j}}=1$,  $p=1,\,2$.  Let $T_p$ be a $\mathbb{C}$-linear transformation mapping $B_{\sqrt{\mu}}$ to $E_{p,\mu}$,  then we have:
\begin{equation*}
||T_1^{-1}\circ T_2||\le \frac{(1+\gamma)^2}{(1-\gamma)^2},\,\,\,\,||T_2^{-1}\circ T_1||\le \frac{(1+\gamma)^2}{(1-\gamma)^2}.
\end{equation*}
\end{cor}
\begin{proof}
We can apply a map $T_1^{-1}$ to the above picture and reduce $E_{1,\mu}$ to be $B_{\sqrt{\mu}}(0)$.  Denote $\tilde{E}_{2,\mu}=T_1^{-1}(E_{2,\mu})=T_1^{-1}\circ T_2(B_{\sqrt{\mu}})$.  Then we know that the eigenvalues of $\tilde{E}_{2,\mu}$ is between $\frac{(1+\gamma)^2}{(1-\gamma)^2}$ and $\frac{(1-\gamma)^2}{(1+\gamma)^2}$.  This implies the result.
\end{proof}
Before we move further,  we first want to explain the idea why Corollary \ref{c3.17} help with proving engulfing property.  Suppose,  say,  we have $S_{\mu}(x_1)\cap S_{\mu}(x_2)\neq \emptyset$.  Take $x_*\in S_{\mu}(x_1)\cap S_{\mu}(x_2)$.   We can normalize $S_{100\mu}(x_1)$ to be close to a unit ball and $S_{\mu}(x_1)$ will then be close to a ball with radius 0.1.  Now we can define a section $\tilde{S}^*$ with height $\frac{1}{100}$ centered at $x'_*$(image of $x_*$ under the new coordinate) using the new coordinate,  so that its shape is comparable to a ball.

Then we go back to the original coordinate,  we get a section $S^*$ centered at $x_*$ with height $\mu$ whose shape is comparable with $S_{\mu}(x_1)$.  Because of Corollary \ref{c3.17},  $S^*$ and $S_{\mu}(x_*)$ will also have similar shapes.  Therefore we see that the shapes of $S_{\mu}(x_*)$ and $S_{\mu}(x_1)$ are comparable (Lemma \ref{l3.18}).  Similarly,  $S_{\mu}(x_*)$ and $S_{\mu}(x_2)$ are also comparable.  Hence,  if we normalize $S_{\mu}(x_1)$ to be close to a unit ball,  the other section $S_{\mu}(x_2)$ will be close to an ellipsoid whose shape is not too eccentric.  Then the engulfing property would follow from the standard engulfing property for balls.

\begin{lem}\label{l3.18}
Let $x_*\in S_{\mu}(x_0)\cap B_{0.8}$ for some $x_0\in B_{0.8}$ and $0<\mu\le \mu_0$.  Let $T_{\mu,x_*}$ and $T_{\mu,x_0}$ be the $\mathbb{C}$-linear transfomation given by Proposition \ref{p3.1},  part (3).  Then for $\sigma$ chosen small enough depending only on $n$,  
\begin{equation*}
||T_{\mu,x_*}^{-1}\circ T_{\mu,x_0}||\le 1.1^3,\,\,\,\,||T_{\mu,x_0}^{-1}\circ T_{\mu,x_*}||\le 1.1^3.
\end{equation*}
\end{lem}
\begin{proof}
First we can find $k\ge 1$ such that $\mu_0^{k+1}<\mu\le \mu_0^k$.  We will work under the coordinate $z^{(k-1)}$,  defined as $z=x_0+\mu_0^{\frac{k-1}{2}} T_{k-1,x_0}(z^{(k-1)})$,  where $T_{k-1,x_0}$ is the $\mathbb{C}$-linear transformation mapping a ball to $E_{\mu_0^{k-1}}(0)$.  We have that
\begin{equation*}
\text{$S_{\mu}(x_0)$ under $z^{(k-1)}$}=\{z^{(k-1)}:u_{x_0,k-1}(z^{(k-1)})\le u_{x_0,k-1}(0)+\tilde{h}_{x_0,k-1}(z^{(k-1)})+\frac{\mu}{\mu_0^{k-1}}\},
\end{equation*}
where $u_{x_0,k-1}$ is the normalized $u$ on $\Omega_{x_0,k-1}$ ($S_{\mu_0^{k-1}}(x_0)$ under $z^{(k-1)}$ which is close to a ball.)
Also we have 
\begin{equation*}
\text{$E_{\mu}(x_0)$ under $z^{(k-1)}$}=\{z^{(k-1)}:\sum_{i,j}(v_{x_0,k-1})_{i\bar{j}}(0)z^{(k-1)}_i\bar{z}^{(k-1)}_j\le \frac{\mu}{\mu_0^{k-1}}\}.
\end{equation*}
Here $v_{x_0,k-1}$ solves the Dirichlet problem (\ref{3.5New}).
On the other hand,  denote $z_*$ to be the image of $x_*$ under $z^{(k-1)}$,  then we know that $z_*\in\text{the image of $S_{\mu}(x_0)$ under $z^{(k-1)}$}$.

We can find a section centered at $z_*$ under $z^{(k-1)}$ using Lemma \ref{l3.6},  applied to $v_{x_0,k-1}$ and $u_{x_0,k-1}$,  $\delta=4\eps$,  $x_0=z_*$,  $\mu$ replaced by $\frac{\mu}{\mu_0^{k-1}}\le \mu_0$,  which is between $\mu_0$ and $\mu_0^2$,   we have that:
\begin{equation}\label{3.16}
\begin{split}
&(1-\gamma)\tilde{E}_{\frac{\mu}{\mu_0^{k-1}}}(z_*)\subset \{z^{(k-1)}\in B_1:(u_{x_0,k-1}-h_*)(z^{(k-1)})\le u_{x_0,k-1}(z_*)+\frac{\mu}{\mu_0^{k-1}}\}\\
&\subset (1+\gamma)\tilde{E}_{\frac{\mu}{\mu_0^{k-1}}}(z_*).
\end{split}
\end{equation}
In the above,  $\gamma=\frac{8\eps\mu_0^{k-1}}{\mu}+3\mu^{\frac{1}{2}}\mu_0^{-\frac{k-1}{2}}.$ It is less that $0.1\sigma$,  because of our choice of $\mu_0$ in (\ref{3.2}) and our assumption that $\eps$ is small (we made a precise choice of $\eps$ in Corollary \ref{c3.7}.)
Also the coefficients of $\tilde{E}_{\frac{\mu}{\mu_0^{k-1}}}(z_*)$ is defined by $(v_{x_0,k-1})_{i\bar{j}}(z_*)$.  Let $T_*$ be the $\mathbb{C}$-linear transform given by Lemma \ref{l3.8} such that $T_*(B_{\sqrt{\mu\mu_0^{1-k}}}(0))=\tilde{E}_{\mu\mu_0^{1-k}}(0)$.  Then similar to Lemma \ref{l3.10},  part (4),  we would have that (using $(v_{x_0,k-1})_{i\bar{j}}$ is close to identity by $\sigma^{\frac{1}{2}}$):
\begin{equation*}
||\tilde{T}_*-I||\le C_n'\sigma^{\frac{1}{2}},\,\,\,||\tilde{T}_*^{-1}-I||\le C_n'\sigma^{\frac{1}{2}}.
\end{equation*}
We can tranform the picture (\ref{3.16}) back to the $z$ variable,  and obtain that:
\begin{itemize}
\item  There is an ellipsoid $E_{\mu}^*$ centered at $x_*$,  such that:
\begin{equation*}
(1-0.1\sigma)E_{\mu}^*(x_*)\subset S^*_{\mu}(x_*)\subset (1+0.1\sigma)E_{\mu}^*(x_*)\subset S_{\mu_0^{k-1}}(x_0),  
\end{equation*}
where $S_{\mu}^*(x_*)=\{z\in S_{\mu_0^{k-1}}(x_0):(u-h_*)(z)\le u(x_*)+\mu\}$,  with $h_*$ being a quadratic pluriharmonic polynomial.
\item With $T_*=T_{k-1,x_0}\circ \tilde{T}_*$,  then we have $T_*(B_{\sqrt{\mu}}(0))=E_{\mu}^*(z_*)$,  and $||T_{k-1,x_0}^{-1}\circ T_*-I||\le 0.1,\,\,\,||T_*^{-1}\circ T_{k-1,x_0}-I||\le 0.1$,  if $\sigma$ is small enough depending only on $n$.
\end{itemize}
On the other hand,  with $S_{\mu}(x_*)$,  $E_{\mu}(x_*)$ being given by Proposition \ref{p3.1},  we also have $(1-0.1\sigma)E_{\mu}(x_*)\subset S_{\mu}(x_*)\subset (1+0.1\sigma)E_{\mu}(x_*)$.  Also we have a $\mathbb{C}$-linear transform $T_{\mu,x_*}$ such that $T_{\mu,x_*}(B_{\sqrt{\mu}}(0))=E_{\mu}(x_*)$.  
Using Corollary \ref{c3.17},  we have
\begin{equation*}
||T_*^{-1}\circ T_{\mu,x*}||\le \frac{(1+0.1\sigma)^2}{(1-0.1\sigma)^2},\,\,\,||T_{\mu,x_*}^{-1}\circ T_*||\le \frac{(1+0.1\sigma)^2}{(1-0.1\sigma)^2}
\end{equation*}
Hence we get that
\begin{equation*}
||T_{k-1,x_0}^{-1}\circ T_{\mu,x_*}||,\,\,\,||T_{\mu,x_*}^{-1}\circ T_{k-1,x_0}||\le 1.1^2,\,\,\text{ if $\sigma$ is small enough}.
\end{equation*}
Finally,  we note that $T_{\mu,x_0}=T_{k-1,x_0}\circ \tilde{T}_{k,x_0}$,  and $||\tilde{T}_{k,x_0}||,\,\,||\tilde{T}_{k,x_0}^{-1}||\le 1+C_n'\sigma^{\frac{1}{2}}$,  so the result would follow if $\sigma$ is small enough depending only on $n$. 
\end{proof}

Then the engulfing property would follow from Lemma \ref{l3.18}.  
\begin{proof}
(of Proposition \ref{p3.15})

First we want to reduce to when $\mu_1$ and $\mu_2$ are comparable.  Indeed,  if $\mu_1< \mu_2$,  we know from Proposition \ref{p3.1} that $S_{\mu_1}(x_1)\subset S_{(1+c(\sigma))\mu_2}(x_1)$ and $c(\sigma)\rightarrow 0$ as $\sigma\rightarrow 0$.  Hence we may assume that $\sigma$ is small enough so that $S_{\mu_1}(x_1)\subset S_{\mu_1'}(x_1)$,  for some  $\mu_2<\mu_1'\le 4\mu_2$.  Hence it will suffice to prove the following statement:

For any $0<\mu_2\le \mu_1\le 4\mu_2$,  if $S_{\mu_1}(x_1)\cap S_{\mu_2}(x_2)\neq \emptyset$,  we have $S_{\mu_1}(x_1)\subset 10S_{\mu_2}(x_2)$.  

First,  we choose $x_*\in S_{\mu_1}(x_1)\cap S_{\mu_2}(x_2)$.   By choosing $\sigma$ small enough,  we know from Lemma \ref{l3.18} that
\begin{equation*}
||T_{\mu_p,x_*}^{-1}\circ T_{\mu_p,x_p}||\le 1.1^3,\,\,\,||T_{\mu_p,x_p}^{-1}\circ T_{\mu_p,x_*}||\le 1.1^3,\,\,\,p=1,\,2.
\end{equation*}
On the other hand,  we know from Lemma \ref{l3.9} that,  with $\sigma$ small enough,  we have
\begin{equation*}
||T_{\mu_1,x_*}^{-1}\circ T_{\mu_2,x_*}||\le1.1,\,\,\,||T_{\mu_2,x_*}^{-1}\circ T_{\mu_1,x_*}||\le 1.1.
\end{equation*}
This follows from the observation that $\mu_1$ and $\mu_2$ must belong to the same level or adjacent levels (either $\mu_0^{k+1}<\mu_2\le \mu_1\le \mu_0^{k}$,  or $\mu_0^{k+1}<\mu_2\le \mu_0^k<\mu_1\le \mu_0^{k-1}$).
Hence we get:
\begin{equation}\label{3.17}
||T_{\mu_1,x_1}^{-1}\circ T_{\mu_2,x_2}||,\,\,\,||T_{\mu_2,x_2}^{-1}\circ T_{\mu_1,x_1}||\le 1.1^6\cdot 1.1\le 2.
\end{equation}
Now we show the containment of sections:
\begin{equation}\label{3.18}
\begin{split}
&S_{\mu_1}(x_1)\subset (1+0.1\sigma)E_{\mu_1}(x_1)= T_{\mu_1,x_1}(B_{(1+0.1\sigma)\sqrt{\mu_1}}(x_{1,*}))\\
&=T_{\mu_2,x_2}\circ T_{\mu_2,x_2}^{-1}\circ T_{\mu_1,x_1}(B_{(1+0.1\sigma)\sqrt{\mu_1}}(x_{1,*}))\subset T_{\mu_2,x_2}(B_{2(1+0.1\sigma)\sqrt{\mu_1}}(x_{1,*}')).
\end{split}
\end{equation}
In the above,  we denote $x_{1,*}$ and $x_{1,*}'$ so that $T_{\mu_1,x_1}(x_{1,*})=x_1$,  $T_{\mu_2,x_2}^{-1}\circ T_{\mu_1,x_1}(x_{1,*})=x_{1,*}'$.
We also used the bound (\ref{3.17}).
On the other hand,  
\begin{equation*}
S_{\mu_2}(x_2)\subset T_{\mu_2,x_2}(B_{(1+0.1\sigma)\sqrt{\mu_2}}(x_2^*)),  
\end{equation*}
where $x_{2,*}=T_{\mu_2,x_2}^{-1}(x_2)$.  It follows that $B_{(1+0.1\sigma)\sqrt{\mu_2}}(x_2^*)\cap B_{2(1+0.1\sigma)\sqrt{\mu_1}}(x_{1,*}')\neq \emptyset$.  
Since $\mu_1\le 4\mu_2$,  it follows that 
\begin{equation}\label{3.19}
B_{2(1+0.1\sigma)\sqrt{\mu_1}}(x_{1,*}')\subset B_{(1+0.1\sigma)9\sqrt{\mu_2}}(x_2^*).
\end{equation}
Hence it follows from (\ref{3.18}) that:
\begin{equation*}
\begin{split}
&S_{\mu_1}(x_1)\subset T_{\mu_2,x_2}(B_{(1+0.1\sigma)9\sqrt{\mu_2}}(x_2^*))\subset T_{\mu_2,x_2}(B_{9.1\sqrt{\mu_2}}(x_2^*))\\
&=9.1E_{\mu_2}(x_2)\subset 10 S_{\mu_2}(x_2).
\end{split}
\end{equation*}
\end{proof}

Finally let us include the following inclusion result for future reference.
\begin{lem}
Let $\sigma$ be small enough depending on $n$.  Then for any $0<\mu\le \frac{\mu_0}{121}$,  we have
\begin{equation*}
10S_{\mu}(x_0)\subset S_{121\mu}(x_0)\subset 12S_{\mu}(x_0).
\end{equation*}
\end{lem}
\begin{proof}
Let $k\ge 1$ be such that $\mu_0^{k+1}<\mu\le \mu_0^k$.  There are two cases to consider:

If $121\mu\le \mu_0^k$,  then $E_{121\mu}(x_0)$ and $E_{\mu}(x_0)$ are defined by the same coefficients,  hence
\begin{equation*}
10S_{\mu}(x_0)\subset 10(1+0.1\sigma)E_{\mu}(x_0).
\end{equation*}
On the other hand,
\begin{equation*}
(1-0.1\sigma)E_{121\mu}(x_0)=(1-0.1\sigma)11E_{\mu}(x_0)\subset S_{121\mu}(x_0).
\end{equation*}
Hence $10S_{\mu}(x_0)\subset S_{121\mu}(x_0)$,  as long as $10(1+0.1\sigma)\le 11(1-0.1\sigma)$.  
On the other hand,
\begin{equation*}
S_{121\mu}(x_0)\subset(1+0.1\sigma)E_{121\mu}(x_0)=11(1+0.1\sigma)E_{\mu}(x_0)\subset 12(1-0.1\sigma)E_{\mu}(x_0)\subset 12 S_{\mu}(x_0).
\end{equation*}

If $\mu_0^k<121\mu\le \mu_0^{k-1}$,  then under $z^{(k-1)}$,  $E_{121\mu}(x_0)$ becomes a ball with radius $11\sqrt{\mu\mu_0^{1-k}}$,  hence $S_{121\mu}(x_0)$ contains the ball $B_{11(1-0.1\sigma)\sqrt{\mu\mu_0^{1-k}}}(0)$ and is contained in $B_{11(1+0.1\sigma)\sqrt{\mu\mu_0^{1-k}}}$ under $z^{(k-1)}$.
On the other hand,  $z^{(k)}$ and $z^{(k-1)}$ differs by a coordinate change $\tilde{T}_{k,x_0}$,  for which we have the bound $||\tilde{T}_{k,x_0}||\le 1+C_n'\sigma^{\frac{1}{2}}$,  $||\tilde{T}_{k,x_0}^{-1}||\le 1+C_n'\sigma^{\frac{1}{2}}$.  Therefore,  $E_{\mu}(x_0)$ is contained in $B_{(1+C_n'\sigma^{\frac{1}{2}})\sqrt{\mu\mu_0^{1-k}}}(0)$ ,  and contains $B_{\frac{\sqrt{\mu\mu_0^{1-k}}}{1+C_n'\sigma^{\frac{1}{2}}}}(0)$ under $z^{(k-1)}$.  Therefore,  working under $z^{(k-1)}$,  one has
\begin{equation*}
\begin{split}
&\text{ image of $10S_{\mu}(x_0)$ in $z^{(k-1)}$}\subset \text{image of $10(1+0.1\sigma)E_{\mu}(x_0)$ under $z^{(k-1)}$}\\
&\subset B_{10(1+0.1\sigma)(1+C_n'\sigma^{\frac{1}{2}})\sqrt{\mu\mu_0^{1-k}}}(0)\subset (1-0.1\sigma)B_{11\sqrt{\mu\mu_0^{1-k}}}(0)\subset \text{image of $S_{121\mu}(x_0)$ under $z^{(k-1)}$}.
\end{split}
\end{equation*}
The above inequality holds as long as $10(1+0.1\sigma)(1+C_n'\sigma^{\frac{1}{2}})\le 11(1-0.1\sigma).$
On the other hand,  still working under $z^{(k-1)}$:
\begin{equation*}
\begin{split}
&\text{Image of $S_{121\mu}(x_0)$ in $z^{(k-1)}$}\subset 11(1+0.1\sigma)B_{\sqrt{\mu
\mu_0^{1-k}}}(0)\\
&\subset 11(1+0.1\sigma)(1+C_n'\sigma^{\frac{1}{2}})\text{$\times E_{\mu}(x_0)$ in $z^{(k-1)}$}\subset 12(1-0.1\sigma)\text{$\times E_{\mu}(x_0)$ in $z^{(k-1)}$}\\
&\subset 12 \text{$\times S_{\mu}(x_0)$ in $z^{(k-1)}$}.
\end{split}
\end{equation*}
In the above,  we need to require that $11(1+0.1\sigma)(1+C_n'\sigma^{\frac{1}{2}})\le 12(1-0.1\sigma)$.
\end{proof}
As a consequence,  we also have the following version of engulfing property:
\begin{cor}
Under the assumptions of Proposition \ref{p3.15},  we have:
\begin{equation*}
S_{\mu_1}(x_0)\subset S_{121\mu_2}(x_0).
\end{equation*}
\end{cor}

\section{Some measure theoretic lemmas}
In this section we will prove a covering lemma which will be used in the $W^{2,p}$ estimate.
In the following,  $m$ always denotes the standard Lebesgue measure.

\begin{lem}\label{covering lemma}
Let $S_{\mu_{\alpha}}(x_{\alpha}) \subset \mathbb{R}^d$ be a family of sets.  Assume that $0<\mu_{\alpha}\le \mu_0$ for all $\alpha$ and $\cup_{\alpha}S_{\mu_{\alpha}}(x_{\alpha})$ is bounded. Assume that the volume of $S_{\mu_{\alpha}}(x_{\alpha})$ is comparable to that of a standard ball with radius $\sqrt{\mu_{\alpha}}$. I.e. there exists a uniform constant $C$ such that 
\begin{equation*}
    \frac{1}{C}m(B_{\sqrt{\mu_{\alpha}}}(0))\le m(S_{\mu_{\alpha}}(x_{\alpha})) \le Cm(B_{\sqrt{\mu_{\alpha}}}(0))
\end{equation*}
Assume that $S_{\mu_{\alpha}}(x_{\alpha})$ satisfies the following engulfing property:

For any $S_{\mu_{\alpha_1}}(x_{\alpha_1})$ and $S_{\mu_{\alpha_2}}(x_{\alpha_2})$ with $S_{\mu_{\alpha_1}}(x_{\alpha_1}) \cap S_{\mu_{\alpha_2}}(x_{\alpha_2}) \neq \emptyset$, if $\sqrt{\mu_{\alpha_1}}\le 2\sqrt{\mu_{\alpha_2}}$, then $S_{\mu_{\alpha_1}}(x_{\alpha_1}) \subset 10S_{\mu_{\alpha_2}}(x_{\alpha_2})$.

Let $X$ be a measurable set with $X \subset \cup_{\alpha}S_{\mu_{\alpha}}(x_{\alpha})$, then one can choose a sequence (finite or infinite) $S_{\mu_i}(x_i)$, such that:

(1) $S_{\mu_i}(x_i)$ are all disjoint.

(2)$X \subset \cup_i 10S_{\mu_i}(x_i)$. 
\end{lem}

\begin{proof}
The proof for this lemma is very similar to the standard Vitali's covering lemma in measure theory.  

First,  we choose a set $S_{\mu_1}(x_1)$ with $\sqrt{\mu_1}>\frac{1}{2}\sup_{\alpha}\sqrt{\mu_{\alpha}}$.  Then we consider all $S_{\mu_{\alpha}}(x_{\alpha})$ which does not intersect $S_{\mu_1}(x_1)$,  and you choose $S_{\mu_2}(x_2)$ so that $\sqrt{\mu_2}>\frac{1}{2}\sup\{\sqrt{\mu_{\alpha}}:S_{\mu_{\alpha}}(x_{\alpha})\cap S_{\mu_1}(x_1)=\emptyset\}$.  Then you consider all $S_{\mu_{\alpha}}(x_{\alpha})$ which does not intersect with $S_{\mu_1}(x_1)$ or $S_{\mu_2}(x_2)$ and you choose $S_{\mu_3}(x_3)$ among those so that $\sqrt{\mu_3}>\frac{1}{2}$sup of $\sqrt{\mu_{\alpha}}$ among all $S_{\mu_{\alpha}}(x_{\alpha})$ which don't intersect with $S_{\mu_1}(x_1)$ or $S_{\mu_2}(x_2)$.  We continue this process.  

This process may stop in finite steps.  If this happens,  then we get a finite sequence $S_{\mu_1}(x_1),\,S_{\mu_2}(x_2),\cdots,\,S_{\mu_N}(x_N)$ such that they are mutually disjoint,  and all $S_{\mu_{\alpha}}(x_{\alpha})$ must intersect with one of them.  Let $i_0$ be the first index such that $S_{\mu_{i}}(x_{i})\cap S_{\mu_{\alpha}}(x_{\alpha})\neq \emptyset$,  then $S_{\mu_{i_0-1}}(x_{i_0-1})\cap S_{\mu_{\alpha}}(x_{\alpha})=\emptyset$,  and due to our inductive choice,  $\sqrt{\mu_{\alpha}}\le 2\sqrt{\mu_{i_0}}$.  Using the engulfing property,  we get:
\begin{equation*}
S_{\mu_{\alpha}}(x_{\alpha})\subset 10S_{\mu_{i_0}}(x_{i_0}).
\end{equation*}
Therefore,  $X\subset \cup_{\alpha}S_{\mu_{\alpha}}(x_{\alpha})\subset \cup_{i=1}^N10S_{\mu_i}(x_i)$.  

The other possibility is that we find an infinite sequence of $\{S_{\mu_i}(x_i)\}_{i=1}^{\infty}$.  They are mutually disjoint because of our construction.  Then we must have that $\mu_i\rightarrow 0$,  since $\cup_{\alpha}S_{\mu_{\alpha}}(x_{\alpha})$ is bounded and the volume of $S_{\mu_{\alpha}}(x_{\alpha})$ is comparable to that of a standard ball with radius $\sqrt{\mu_{\alpha}}$.  In particular,  if you define $d_i=\sup\{\mu_{\alpha}:S_{\mu_{\alpha}}(x_{\alpha})\cap S_{\mu_j}(x_{j})=\emptyset,\,\,1\le j\le i\}$,  then we have $d_i\rightarrow 0$.  It follows that any $S_{\mu_{\alpha}}(x_{\alpha})$,  there exists $S_{\mu_{i_0}}(x_{i_0})$ such that $\sqrt{\mu_{\alpha}}\le 2\sqrt{\mu_{i_0}}$ such that $S_{\mu_{\alpha}}(x_{\alpha})\cap S_{\mu_{i_0}}(x_{i_0})\neq \emptyset$.  
Hence if you use the engulfing property,  you see that $S_{\mu_{\alpha}}(x_{\alpha})\subset 10S_{\mu_{i_0}}(x_{i_0})$.
\end{proof}

Using this,  we can follow the usual proof of Lebesgue differentiation theorem to conclude that:

\begin{lem}\label{L0.10}
Let $\{S_{\mu}(x)\}_{0<\mu\le \mu_0,\,x\in B_{0.8}}$ be a family of sets such that:
\begin{enumerate}
\item $S_{\mu}(x)\subset B_1$,  for all $0<\mu\le \mu_0$,  and $x\in B_{0.8}$,
\item There is $C>0$,  such that for all $0<\mu\le \mu_0$ and all $x\in B_{0.8}$,  $\frac{1}{C}m(B_{\sqrt{\mu}}(0))\le m(S_{\mu}(x))\le Cm(B_{\sqrt{\mu}}(0))$,
\item For any $S_{\mu_1}(x_1)$ and $S_{\mu_2}(x_2)$,  if $\sqrt{\mu_1}\le 2\sqrt{\mu_2}$ and $S_{\mu_1}(x_1)\cap S_{\mu_2}(x_2)\neq \emptyset$,  then $S_{\mu_1}(x_1)\subset 10S_{\mu_2}(x_2)$,
\item $diam\,S_{\mu}(x)$ tends to 0 as $\mu\rightarrow 0$,  uniformly for $x\in B_{0.8}$. 
\end{enumerate}
 Let $f:B_{0.8}\rightarrow \mathbb{R}$ be an $L^1$ function.  Then for $m$-a.e.  $x\in B_{0.8}$,  we have:
\begin{equation*}
\lim\sup_{x\in S_{\mu_{\alpha}}(x_{\alpha}),\,\mu_{\alpha}\rightarrow 0}\frac{1}{m(S_{\mu_{\alpha}}(x_{\alpha}))}\int_{S_{\mu_{\alpha}}(x_{\alpha})}|f(y)-f(x)|dm(y)=0.
\end{equation*}
In particular,  for all measurable set $A$,  we have:
\begin{equation*}
\lim\inf_{x\in S_{\mu_{\alpha}}(x_{\alpha}),\,\mu_{\alpha}\rightarrow 0}\frac{m(S_{\mu_{\alpha}}(x_{\alpha})\cap A)}{m(S_{\mu_{\alpha}}(x_{\alpha}))}=1,\,\,\,a.e.  \,x\in A.
\end{equation*}
\end{lem}
\begin{proof}
The proof of this lemma follows the proof of the standard Lebesgue differentiation theorem.

First,  we define the maximal function: given $f\in L^1(B_{0.8})$,  
\begin{equation*}
\mathcal{M}(f)(x)=\sup_{x\in S_{\mu}(x'),\,0<\mu\le \mu_0,\,x'\in B_{0.8}}
\frac{1}{m(S_{\mu}(x'))}\int_{S_{\mu}(x')}f(y)dm(y).
\end{equation*}
As in the proof of Lebesgue differentiation theorem,  the result would follow from the following estimate:
\begin{equation}\label{4.2NNN}
m\{x\in B_{0.8}:\mathcal{M}(|f|)(x)>t\}\le C_n\frac{||f||_{L^1}}{t},\,\,\,\forall t>0.
\end{equation}
Indeed,  for $t>0$,  we can define:
\begin{equation*}
\Omega_t:=\{x\in B_{0.8}:\lim\sup_{x\in S_{\mu_{\alpha}}(x_{\alpha}),\,\mu_{\alpha}\rightarrow 0}\frac{1}{m(S_{\mu_{\alpha}}(x_{\alpha}))}\int_{S_{\mu_{\alpha}}(x_{\alpha})}|f(y)-f(x)|dm(y)>t\}.
\end{equation*}
We only need to show that $m(\Omega_t)=0$ for any $t>0$.   On the other hand,  for any $g\in C(\bar{B}_{0.8})$,  we have 
\begin{equation}\label{5.2NNN}
\begin{split}
&\frac{1}{m(S_{\mu_{\alpha}}(x_{\alpha}))}\int_{S_{\mu_{\alpha}}(x_{\alpha})}|f(y)-f(x)|dm(y)\le\mathcal{M}(|f-g|)(x)\\
&+ \frac{1}{m(S_{\mu_{\alpha}}(x_{\alpha}))}\int_{S_{\mu_{\alpha}}(x_{\alpha})}|g(y)-g(x)|dm(y)+|f-g|(x).
\end{split}
\end{equation}
Therefore,  
\begin{equation*}
\Omega_t\subset\{x:\mathcal{M}(|f-g|)>\frac{t}{3}\}\cup\{x:|f-g|(x)>\frac{t}{3}\}.
\end{equation*}
Here we implicitly used item (4) of the assumption,  so that the middle term in (\ref{5.2NNN}) tends to zero as $\mu_{\alpha}\rightarrow 0$.
Hence
\begin{equation*}
m(\Omega_t)\le m\{x:\mathcal{M}(|f-g|)>\frac{t}{3}\}+m\{x:|f-g|(x)>\frac{t}{3}\}\le C_n\frac{3||f-g||_{L^1}}{t}+\frac{||f-g||_{L^1}}{t}.
\end{equation*}
We can then choose a sequence $g_j\in C(\bar{B}_{0.8})$ such that $g_j\rightarrow f$ in $L^1$,  so that we may conclude that $m(\Omega_t)=0$.

Now it only remains to show (\ref{4.2NNN}).  Denote the set on the left hand side to be $E$,  then for any $x\in E$,  we can find $0<\mu_x\le \mu_0$ and $y_x\in B_{0.8}$ with $x\in S_{\mu_x}(y_x)$,  such that:
\begin{equation*}
tm(S_{\mu_x}(y_x))\le \int_{S_{\mu_x}(y_x)}|f(y)|dm(y).
\end{equation*} 
Hence we get a covering of $E$: $E\subset \cup_{x\in E}S_{\mu_x}(y_x)$.
Now we are in a position to apply Lemma \ref{covering lemma} to choose a countable sequence $S_{\mu_i}(x_i)$,  which is mutually disjoint,  and $E\subset \cup_i10S_{\mu_i}(x_i)$.  Hence
\begin{equation*}
\begin{split}
&m(E)\le m(\cup_i10S_{\mu_i}(x_i))\le \sum_im(10S_{\mu_i}(x_i))=10^{2n}\sum_im(S_{\mu_i}(x_i))\\
&\le 10^{2n}\sum_i\frac{1}{t}\int_{S_{\mu_i}(x_i)}|f(y)|dm(y)\le 10^{2n}\frac{||f||_{L^1}}{t}.
\end{split}
\end{equation*}
This proves (\ref{4.2NNN}).
\end{proof}

Another lemma we will need is:

\begin{lem}\label{l0.11}
Let $\{S_{\mu}(x)\}_{0<\mu\le \mu_0,\,x\in B_{0.8}}$ satisfy the assumptions of Lemma \ref{L0.10},  and we assume additionally:
\begin{equation}\label{4.2N}
10S_{\mu}(x)\subset S_{121\mu}(x)\subset 12S_{\mu}(x),\,\,\,\text{for any $x\in B_{0.8}$,  $0<\mu\le \frac{\mu_0}{121}$.}
\end{equation}
Let $X,\,Y\subset B_{0.8}$ be two measurable sets.  Let $0<\bar{\eps}<1$.  Assume that:
\begin{enumerate}
\item For any $x_0\in B_{0.8}$,  $m(S_{\mu}(x_0)\cap X)< \bar{\eps}m(S_{\mu}(x_0))$ for any $\frac{\mu_0}{4}\ge \mu\ge \frac{\mu_0}{484}$. 
\item For any $S_{\mu}(x)$ with $m(S_{\mu}(x)\cap X)\ge \bar{\eps}m(S_{\mu}(x))$ and $\mu\le \frac{\mu_0}{2}$,  one has $S_{\mu}(x)\subset Y$.
\end{enumerate}
Then
\begin{equation*}
m(X)\le 12^{2n}\bar{\eps}m(Y).
\end{equation*}
\end{lem}
\begin{proof}
For a.e.  any $x_0\in X$,  we know from the previous lemma that:
\begin{equation}\label{4.2}
\lim_{\mu \rightarrow 0}\frac{m(X\cap S_{\mu}(x))}{m(S_{\mu}(x))}=1.
\end{equation}
Hence,  if we define 
$\mu_x'=\sup\{0<\mu\le \frac{1}{4}\mu_0: m(X\cap S_{\mu}(x))\ge \bar{\eps}S_{\mu}(x)\}$,  then $0<\mu_x'\le \frac{\mu_0}{484}$ for $x\in X$ satisfying (\ref{4.2}). 
 For such $x$,  we may choose $\mu_x$,  such that $\frac{5}{6}\mu_x'\le \mu_x\le \mu_x'$,  and that $m(X\cap S_{\mu_x}(x))\ge \bar{\eps}S_{\mu_x}(x)$.  We may assume without loss of generality that (\ref{4.2}) holds for all $x\in X$,  then we get a covering of $X$: $X\subset \cup_{x\in X}S_{\mu_x}(x)$ (otherwise we get a covering of $X$ modulo a measure zero set).

Then we may use Lemma \ref{covering lemma} to obtain a countable sequence $S_{\mu_i}(x_i)$,  such that $S_{\mu_i}(x_i)$ are mutually disjoint,  with $X\subset\cup_i10S_{\mu_i}(x_i)$.  Hence
\begin{equation*}
\begin{split}
&m(X)\le \sum_im(X\cap 10S_{\mu_i}(x_i))\le \sum_im(X\cap S_{121\mu_i}(x_i))\le \sum_i\bar{\eps}m(S_{121\mu_i}(x_i))\\
&\le \bar{\eps}\sum_im(12S_{\mu_i}(x_i))\le \bar{\eps}12^{2n}\sum_im(S_{\mu_i}(x_i))\le 12^{2n}\bar{\eps}m(Y).
\end{split}
\end{equation*}
In the second inequality,  we used (\ref{4.2N}).  

In the third inequality,  we used that $\frac{\mu_0}{4}\ge 121\mu_i>\mu_{x_i}'$,  since our choice of $\mu_x$ guarantees $\mu_i\ge \frac{5}{6}\mu_{x_i}'$.  Therefore $m(X\cap S_{121\mu_i}(x_i))<\bar{\eps}m(S_{121\mu_i}(x_i))$.

In the forth inequality,  we used (\ref{4.2N}) again.

In the last inequality,  we used that $S_{\mu_i}(x_i)$ are disjoint,  and contained in $Y$,  due to assumption (2) of this lemma.
\end{proof}

\section{The $W^{2,p}$ estimate}
\begin{defn}\label{def ak}
We define $D_k$ to be the set of $z_0\in \bar{B}_{0.8}$ such that for any $0<\mu \le \mu_0$,  $S_{\mu}(z_0)\subset B(z_0, \sqrt{10^k\mu})$.  Define $A_k=\bar{B}_{0.8}-D_k$.
\end{defn}
In the above,  one should think of $D_k$ to be the ``good set" and $A_k$ the ``bad set".  

Roughly speaking,  $D_k$ is the set on which $\lambda_i(z_0)\ge 10^{-k}$.  So that heuristically,  we can conclude that $\lambda_i(z_0)\le10^{k(n-1)}$ since the equation is $\Pi_i\lambda_i=f$.  
To see this picture,  we can pretend that $S_{\mu}(z_0)\approx \{z:\sum_{i,j}u_{i\bar{j}}(z_0)(z-z_0)_i\overline{(z-z_0)}_j\le \mu\}$.  The requirement that this ellipsoid is contained in $B_{\sqrt{10^k\mu}}(z_0)$ implies that $\frac{1}{\sqrt{\lambda_i}}\le 10^{\frac{k}{2}}$.
This of course needs to be made rigorous since our solution $u$ is merely a viscosity solution (hence only continuous.)

In order to show the $W^{2,p}$ estimates,  there are roughly two steps.  
\begin{enumerate}
\item For every $k\ge 1$,  $m(A_k\cap B_{r_k}(0))\le (12^{2n}\bar{\eps})^{k-1}m(B_{0.7})$,  for $k\ge 1$,  by choosing $\sigma$ and $\eps$ chosen sufficiently small and $r_k=r_{k-1}-\frac{1}{10}2^{-k}$,  $r_0=0.7$.
\item Show that for a.e.  $x\in D_k$,  there is a paraboloid with opening $M_0^{k(n-1)}$ touching $u$ from above at $x$,  and a paraboloid with opening $M_0^{-k}$ touching $u$ from below at $x$.  This is the viscosity interpretation of $D^2u(x)\le M_0^{(k-1)n}$ and $u_{i\bar{j}}(x)\ge M_0^{-k}$.
\end{enumerate}
We will carry out steps (1) and (2) in the following two subsections.  For the convenience of argument,  we will assume that $u\in C^2(\Omega)\cap C(\bar{\Omega})$,  and obtain a $quantitative$ $W^{2,p}$ bound on $B_{\frac{1}{2}}$.  Then the general case would follow from an approximation argument.
\subsection{Power decay of the measure of bad set}
The plan is to use Lemma \ref{l0.11} with $X=A_{k+1}\cap B_{r_{k+1}}(0)$,  $Y=A_k\cap B_{r_k}(0)$.  Note that the sections $S_{\mu}(x)$ satisfy all the assumptions of that lemma.
Therefore,  we just have to show the following two things:
\begin{enumerate}
\item Choosing $M_0$ large enough depending on $\bar{\eps}$ and $n$,  we have $m(S_{\mu}(x_0)\cap A_1\cap B_{r_1}(x_0))<\bar{\eps}m(S_{\mu}(x_0))$ for any $\frac{\mu_0}{2}\ge\mu\ge \frac{\mu_0}{242}$.
\item For all $S_{\mu}(x_0)$ with $m(S_{\mu}(x_0)\cap A_{k+1}\cap B_{r_{k+1}}(0))\ge \bar{\eps}m(S_{\mu}(x_0))$ and $\mu\le \frac{\mu_0}{2}$,  one has $S_{\mu}(x_0)\subset A_k\cap B_{r_k}(0)$,  by choosing $\sigma$ and $\eps$ small enough.
\end{enumerate}
We will start with the following lemma,  which is the analogue of Lemma 6 in \cite{Ca1}.  
\begin{lem}\label{l5.2}
Let $u_0$ be a $C^2$ solution to $\det(u_0)_{i\bar{j}}=f_0$ in $\Omega$ with $B_{1-\gamma}(0)\subset \Omega\subset B_{1+\gamma}$.  
We also assume that $|f_0-1|<\eps$ on $\Omega$.  Let $v_0$ be the solution of $\det(v_0)_{i\bar{j}}=1$ on $\Omega$ and $v_0|_{\partial \Omega}=0$.  Then there exists a dimensional constant $C_{6,n}$ such that for all $\gamma$ and $\eps$ small enough (depending only on dimension),
\begin{equation*}
\frac{m(\{x\in B_{\frac{1}{2}}:\Gamma(u_0-\frac{1}{2}v_0)=u_0-\frac{1}{2}v_0\})}{m(B_{\frac{1}{2}})}\ge 1-C_{6,n}\eps^{\frac{1}{2}}-C_{6,n}\gamma^{\frac{1}{2}}.
\end{equation*}
In the above,  $\Gamma(u_0-\frac{1}{2}v_0)$ is the convex envelope of $u_0-\frac{1}{2}v_0$ in $B_{0.9}$,  defined as:
\begin{equation*}
\Gamma(u_0-\frac{1}{2}v_0)(z)=\sup\{l(z):\text{$l(z)$ is affine and $l\le u_0-\frac{1}{2}v_0$ on $B_{0.9}$}\}.
\end{equation*}
\end{lem}
\begin{proof}
The proof follows similar lines as Lemma 6 of \cite{Ca1}.  First,  we may use Lemma \ref{dirichlet problem estimate} and Corollary \ref{c3.4} to conclude that,  with $\gamma$ chosen small enough:
\begin{equation}\label{5.1}
|D^m(v_0-(|z|^2-1))|_{B_{0.9}}\le C_n\gamma^{1-\frac{m}{4}},\,\,\,m=0,\,1,\,2,\,3.
\end{equation}
 In particular,  we know that $v_0$ is strictly convex from (\ref{5.1}),  after choosing $\gamma$ small enough.
Also by maximum principle,  
\begin{equation*}
(1+3\eps)v_0\le (1+\eps)^{\frac{1}{n}}v_0\le u_0\le (1-\eps)^{\frac{1}{n}}v_0\le (1-3\eps)v_0.
\end{equation*}
So we get:
\begin{equation*}
(\frac{1}{2}+3\eps)v_0\le \Gamma(u_0-\frac{1}{2}v_0)\le (\frac{1}{2}-3\eps)v_0.
\end{equation*}
In the following,  we will simply denote $\Gamma(u_0-\frac{1}{2}v_0)$ by $\Gamma$.
We make the following claim:
\begin{claim}\label{claim5.4}
\begin{equation}
\nabla\big((\frac{1}{2}-3\eps)v_0\big)(B_{\frac{1}{2}-\sqrt{96\eps}})\subset \nabla \Gamma(B_{\frac{1}{2}}).
\end{equation}
\end{claim}
First we use Claim \ref{claim5.4} to finish the proof and then prove the Claim \ref{claim5.4} itself.

Indeed,  we have:
\begin{equation}\label{5.3}
\begin{split}
&m(\nabla \Gamma(B_{\frac{1}{2}}))\ge m\big(\nabla (\frac{1}{2}-3\eps)v_0\big)(B_{\frac{1}{2}-\sqrt{96\eps}})=\int_{B_{\frac{1}{2}-\sqrt{96\eps}}}\det\big((\frac{1}{2}-3\eps)D^2v_0\big)dm(x)\\
&\ge (\frac{1}{2}-3\eps)^{2n}(2-C_n\gamma^{\frac{1}{2}})^{2n}m(B_{\frac{1}{2}-\sqrt{96\eps}})\ge m(B_{\frac{1}{2}})(1-C_{5,n}\eps^{\frac{1}{2}}-C_{5,n}\gamma^{\frac{1}{2}}).
\end{split}
\end{equation}
In the equality of the first line above,  we used that $v_0$ is strictly convex on $B_{0.9}$,  hence $\nabla v_0$ is injective.

In the first inequality of the second line,  we used (\ref{5.1}) with $m=2$ to get that $|D^2(v_0-|z|^2)|\le C_n'\gamma^{\frac{1}{2}}$.

On the other hand,  we have,  using the Lemma \ref{l5.4} below: 
\begin{equation}\label{5.4}
\begin{split}
&m(\nabla\Gamma(B_{\frac{1}{2}}))\le \int_{B_{\frac{1}{2}}\cap\{\Gamma=u_0-\frac{1}{2}v_0\}}\big(2(1+\eps)^{\frac{1}{n}}-\det(D^2(\frac{1}{2}v_0))^{\frac{1}{2n}}\big)dm(x)\\
&\le \int_{B_{\frac{1}{2}}\cap \{\Gamma=u_0-\frac{1}{2}v_0\}}\big(2(1+\eps)^{\frac{1}{n}}-(1-C_n'\gamma^{\frac{1}{2}})\big)^{2n}dm(x)\\
&\le (1+C_{5,n}\gamma^{\frac{1}{2}}+C_{5,n}\eps)^{2n}m(B_{\frac{1}{2}}\cap\{\Gamma=u_0-\frac{1}{2}v_0\}).
\end{split}
\end{equation}
The result follows from combing (\ref{5.3}) and (\ref{5.4}).

Now it only remains to prove Claim \ref{claim5.4}.  Indeed,  let $p=\nabla (\frac{1}{2}-3\eps)v_0(x_0)$,  with $x_0\in B_{\frac{1}{2}-\sqrt{96\eps}}$.  Denote $l_{x_0}=(\frac{1}{2}-3\eps)v_0(x_0)+p\cdot(x-x_0)$.  We just need to show that:
\begin{equation*}
\{z\in B_{0.9}:\Gamma(z)<l_{x_0}(z)\}\subset \{(z\in B_{0.9}:(\frac{1}{2}+3\eps)v_0\le l_{x_0}\}\subset B_{\sqrt{96\eps}}(x_0).
\end{equation*}
That is,  the minimum of $\Gamma-l_{x_0}$ is achieved in the interior of $B_{\sqrt{96\eps}}(x_0)\subset B_{\frac{1}{2}}$,  giving $p\in \nabla \Gamma(B_{\frac{1}{2}})$.   The first inclusion above is obvious since $\Gamma\ge (\frac{1}{2}+3\eps)v_0$.  To see the second inclusion,  we need the following calculation:

First,  we note that,  for $z\in B_{0.9}$,  we have:
\begin{equation*}
v_0(z)\ge v_0(x_0)+\nabla v_0(x_0)\cdot (z-x_0)+(2-C_n\gamma^{\frac{1}{2}})|z-x_0|^2.
\end{equation*}
In the above,  we used (\ref{5.1}) again.
Therefore,
\begin{equation*}
\begin{split}
&\{z\in B_{0.9}:(\frac{1}{2}+3\eps)v_0\le l_{x_0}\}\subset \{z\in B_{0.9}:(\frac{1}{2}+3\eps)(v_0(x_0)+\nabla v_0(x_0)\cdot (z-x_0)\\
&+(2-C_n\gamma^{\frac{1}{2}})|z-x_0|^2\}\le (\frac{1}{2}-3\eps)(v_0(x_0)+\nabla v_0(x_0)\cdot (z-x_0))\}
\end{split}
\end{equation*}
So that the above implies 
\begin{equation}
(\frac{1}{2}+3\eps)(2-C_n\gamma^{\frac{1}{2}})|z-x_0|^2\le -6\eps(v_0(x_0)+\nabla v_0(x_0)\cdot (z-x_0)).
\end{equation}
In the above,  we may assume $|v_0(x_0)|\le 1.1$ and $|\nabla v_0(x_0)|\le 2$,  and we get $|z-x_0|^2\le \frac{6\eps(2+2\cdot 2)}{(\frac{1}{2}+3\eps)(2-C_n\gamma^{\frac{1}{2}})}\le 96\eps$,  and Claim \ref{claim5.4} is proved.
\end{proof}
In the above proof,  we used the following lemma,  which is the analogue of Lemma 5 in \cite{Ca1}.
\begin{lem}\label{l5.4}
Let $u_0$ and $v_0$ be as given by Lemma \ref{l5.2},  with $\gamma$ small enough so that (\ref{5.1}) holds.  Denote $\Gamma=\Gamma(u_0-\frac{1}{2}v_0)$ as defined by Lemma \ref{l5.2}.  Then we have:
\begin{equation*}
\mathcal{M}(\Gamma)\le \big(2(1+\eps)^{\frac{1}{n}}-\det(D^2(\frac{1}{2}v_0))^{\frac{1}{2n}}\big)^{2n}\chi_{\Gamma=u_0-\frac{1}{2}v_0}\text{ on $B_{0.8}$}.
\end{equation*}
Here $\mathcal{M}(\Gamma)$ is the (real) Monge-Ampere measure of a convex function,  defined as $\mathcal{M}(\Gamma)(E)=m(\partial\Gamma(E))$.
\end{lem}
\begin{proof}
Since $u_0$ is $C^2$,  then we know that $\Gamma$ is $C^{1,1}$ on $B_{0.8}$.  Hence for any Borel set $E\subset B_{0.8}$,  one would have: $\mathcal{M}(\Gamma)(E)=\int_E\det(D^2\Gamma)(x)dm(x)$.  Moreover,  it is a general fact that $\mathcal{M}(\Gamma)$ is concentrated on the contact set $\{\Gamma=u_0-\frac{1}{2}v_0\}$.  On such a set,  we would have:
\begin{equation*}
D^2u_0\ge D^2\Gamma+D^2(\frac{1}{2}v_0),\,\,\,a.e.  
\end{equation*}
Therefore
\begin{equation*}
\det(D^2\Gamma)^{\frac{1}{2n}}+\det(D^2(\frac{1}{2}v_0))^{\frac{1}{2n}}\le \det(D^2u_0)^{\frac{1}{2n}}\le \big(2^{2n}\det(u_0)_{i\bar{j}}^2\big)^{\frac{1}{2n}}\le 2(1+\eps)^{\frac{1}{n}}.
\end{equation*}
In the above,  we used the concavity of the function $A\mapsto \det^{\frac{1}{d}}(A)$,  restricted to positive definite $d\times d$ symmetric matrices.
\end{proof}
As a corollary to Lemma \ref{l5.2},  we get:
\begin{cor}\label{c5.5}
Let $u_0$ be as Lemma \ref{l5.2} and $C_n$ be from (\ref{5.1}).  Define
$E$ to be the subset in $B_{0.8}$ such that $x_0\in E$ if and only if there is a paraboloid with opening $\frac{1}{2}(1-C_n\gamma^{\frac{1}{2}})$ touching $u_0$ from below at $x_0$ in $B_{0.9}$.  Then we have:
\begin{equation*}
\frac{m(E\cap B_{\frac{1}{2}})}{m(B_{\frac{1}{2}})}\ge 1-C_{6,n}\eps^{\frac{1}{2}}-C_{6,n}\gamma^{\frac{1}{2}},
\end{equation*}
where $C_{6,n}$ is from Lemma \ref{l5.2}.
\end{cor}
\begin{proof}
We just need to show that $\Gamma(u_0-\frac{1}{2}v_0)=u_0-\frac{1}{2}v_0$ has the property described in this lemma.  Let $x_0\in B_{\frac{1}{2}}$ with $\Gamma(u_0-\frac{1}{2}v_0)(x_0)=(u_0-\frac{1}{2}v_0)(x_0)$.  Since $\Gamma$ is a convex function on $B_{0.9}$,  we may find $p\in \mathbb{C}^n$ which defines a supporting plane for $\Gamma$,  then we have:
\begin{equation*}
(u_0-\frac{1}{2}v_0)(z)\ge \Gamma(u_0-\frac{1}{2}v_0)(z)\ge (u_0-\frac{1}{2}v_0)(x_0)+p\cdot(z-x_0),\,\,\,z\in B_{0.9}.
\end{equation*}
On the other hand,  from (\ref{5.1}),  we see that:
\begin{equation*}
v_0(z)\ge v_0(x_0)+\nabla v_0(x_0)\cdot (z-x_0)+(1-C_n\gamma^{\frac{1}{2}})|z-x_0|^2,\,\,z\in B_{0.9}.
\end{equation*}
Hence on $B_{0.9}:$
\begin{equation*}
u_0(z)\ge u_0(x_0)+(p+\frac{1}{2}\nabla v_0(x_0))\cdot (z-x_0)+\frac{1}{2}(1-C_n\gamma^{\frac{1}{2}})|z-x_0|^2.
\end{equation*}
The above has equality at $z=x_0$ and the right hand side defines a paraboloid with opening $\frac{1}{2}(1-C_n\gamma^{\frac{1}{2}})$.
\end{proof}
Having a paraboloid touching below at a point $x_0$ is a very strong condition,  and it will imply the control of the shape of $S_{\mu}(x_0)$ on all scales of $\mu$,  together with a control on the associated pluriharmonic function $h_{\mu,x_0}$.  More precisely
\begin{lem}\label{l5.6}
Let $u$ be a function on $B_{0.9}$.  Let $x_0\in B_{0.8}$ be such that there is a paraboloid with opening $\kappa>0$,  touching $u$ from below at $x_0$ in $B_{0.9}$.  
 Let $0<\tilde{\gamma}<1$,  $\mu>0$,  and $A=a_{i\bar{j}}$ be positive Hermitian with $\det A=1$.  Define $E(x_0)=\{z:\sum_{i,j}a_{i\bar{j}}(z-x_0)_i\overline{(z-x_0)}_j\le \mu\}$.  Let $h(z)$ be a degree 2 pluriharmonic polynomial with $h(x_0)=0$.  
\begin{equation*}
(1-\frac{1}{2}\tilde{\gamma})E(x_0)\subset \{z\in B_{0.9}:(u-h)(z)\le u(x_0)+\mu\}\subset \subset B_{0.9}.
\end{equation*}
Then we have the following estimates for $A$ and $h(z)$:
\begin{equation*}
\begin{split}
&\kappa(1-\tilde{\gamma})^2\le \lambda_i(A)\le \kappa^{1-n}(1-\tilde{\gamma})^{2(1-n)},\,\,\,1\le i\le n.\\
&-\frac{1}{(1-\tilde{\gamma})^{2n}}\kappa^{1-n}I_{2n}\le D^2h\le \frac{n-1}{(1-\tilde{\gamma})^{2n}}\kappa^{1-n}I_{2n}.
\end{split}
\end{equation*}
\end{lem}
\begin{proof}
Denote $S(x_0)=\{z\in B_{0.9}:(u-h)(z)\le u(x_0)+\mu\}$.  First we see that 
\begin{equation*}
\sum_{i,j}a_{i\bar{j}}(z-x_0)_i\overline{(z-x_0)}_j>(1-\tilde{\gamma})^2\mu\text{ on $(1-\frac{1}{2}\tilde{\gamma})\partial E(x_0)$.}
\end{equation*}
Hence 
\begin{equation}\label{5.6}
\frac{1}{(1-\tilde{\gamma})^2}\sum_{i,j}a_{i\bar{j}}(z-x_0)_i\overline{(z-x_0)}_j+u(x_0)>(u-h)(z),\,\,\,z\in (1-\frac{1}{2}\tilde{\gamma})\partial E(x_0).
\end{equation}
On the other hand,  we use that $u$ is being touched below by a paraboloid,  we get
\begin{equation*}
u(z)\ge \kappa|z-x_0|^2+l_{x_0}(z),\,\,z\in B_{0.9}.
\end{equation*}
Here $l_{x_0}(z)$ is an affine function with $l_{x_0}(x_0)=u(x_0)$.  
Hence
\begin{equation}\label{5.7}
\frac{1}{(1-\tilde{\gamma})^2}\sum_{i,j}a_{i\bar{j}}(z-x_0)_i\overline{(z-x_0)}_j+u(x_0)>\kappa|z-x_0|^2+l_{x_0}(z)-h(z),\,\,\,z\in (1-\frac{1}{2}\tilde{\gamma})\partial E(x_0).
\end{equation}
However,  the above inequality has equality with $z=x_0$.  Hence LHS of (\ref{5.7})$-$RHS of (\ref{5.7}) has a minimum in the interior of $S(x_0)$.  Denote this point to be $x_0'$.  Taking the complex Hessian at $x_0'$,  we see that,  
\begin{equation*}
\frac{1}{(1-\tilde{\gamma})^2}a_{i\bar{j}}\ge \kappa I.
\end{equation*}
That is $\lambda_i(A)\ge \kappa(1-\tilde{\gamma})^2.$
On the other hand,  since $\Pi_i\lambda_i(A)=1$,  we see that $\lambda_i(A)\le \kappa^{1-n}(1-\tilde{\gamma})^{2(1-n)}$.

Next we take the full Hessian of (\ref{5.7}) at $x_0'$,  we get:
\begin{equation*}
\frac{1}{(1-\tilde{\gamma})^2}D^2\big(\sum_{i,j}a_{i\bar{j}}(z-x_0)_i\overline{(z-x_0)_j}\big)\ge 2\kappa I-D^2h.
\end{equation*}
So that
\begin{equation*}
D^2h\ge -\frac{1}{(1-\tilde{\gamma})^{2n}}\kappa^{1-n}I.
\end{equation*}
On the other hand,  since $h$ is harmonic,  we get that:
\begin{equation*}
D^2h\le \frac{n-1}{(1-\tilde{\gamma})^{2n}}\kappa^{1-n}.
\end{equation*}
\end{proof}
As a consequence,  we get that:
\begin{cor}\label{c5.7}
Let $u_0$ be as stated in Lemma \ref{l5.2},  with $\gamma$ and $\eps$ small enough as required by that lemma.
Define $D$ to be the subset of $B_{0.8}$ such that $x_0\in D$ if and only $S_{\mu}(x_0)\subset B_{\sqrt{M_1\mu}}(x_0)$ for any $0<\mu\le \mu_0$,  where $M_1=2(1+0.1\sigma)^2(1-C_n\gamma^{\frac{1}{2}})^{-1}(1-0.2\sigma)^{-2}$,  where $S_{\mu}(x_0)$ is the family of sections constructed in Section 2 (applied to $u_0$).
Then we have:
\begin{equation*}
\frac{m(B_{\frac{1}{2}}\cap D)}{m(B_{\frac{1}{2}})}\ge 1-C_{6,n}\gamma^{\frac{1}{2}}-C_{6,n}\eps^{\frac{1}{2}}.
\end{equation*}
\end{cor}
\begin{proof}
First,  from Corollary \ref{c5.5},  we just need to show that $E\cap B_{\frac{1}{2}}\subset D\cap B_{\frac{1}{2}}$.  Indeed,  let $x_0\in E$ such that there is a paraboloid with opening $\frac{1}{2}(1-C_n\gamma^{\frac{1}{2}})$,  with $C_n$ coming from (\ref{5.1}).  On the other hand,  our construction in Section 2 gives $S_{\mu}(x_0)$,  $0<\mu\le \mu_0$ with 
\begin{equation*}
(1-0.1\sigma)E_{\mu}(x_0)\subset S_{\mu}(x_0)\subset(1+0.1\sigma)E_{\mu}(x_0),
\end{equation*}
where $S_{\mu}(x_0)=\{z:u(z)\le h_{\mu,x_0}(z)+u(x_0)+\mu\}$ and $E_{\mu}(x_0)=\{\sum_{i,j}a_{\mu,x_0,i\bar{j}}(z-x_0)_i\overline{(z-x_0)_j}\le \mu\}.$.  Now we are in a position to use Lemma \ref{l5.6} (we could assume $\mu_0$ small enough earlier so that $S_{\mu}(x_0)\subset B_{0.9}$ with $0<\mu\le \mu_0,\,x_0\in B_{0.8}$),  with $\kappa=\frac{1}{2}(1-C_n\gamma^{\frac{1}{2}})$,  $\tilde{\gamma}=0.2\sigma$, $A=a_{\mu,x_0,i\bar{j}}$,  to get:
\begin{equation*}
\frac{1}{2}(1-C_n\gamma^{\frac{1}{2}})(1-0.2\sigma)^2I_n\le a_{\mu,x_0,i\bar{j}}\le 2^{n-1}(1-C_n\gamma^{\frac{1}{2}})^{1-n}(1-0.2\sigma)^{2(1-n)}I_n.
\end{equation*}
This would imply that $E_{\mu}(x_0)\subset B_{(1-C_n\gamma^{\frac{1}{2}})^{-\frac{1}{2}}(1-0.2\sigma)^{-1}\sqrt{2\mu}}(x_0)$.  
So that $S_{\mu}(x_0)\subset B_{\sqrt{M_1\mu}}(x_0)$ for all $0<\mu\le \mu_0$.
\end{proof}

Now we are ready to prove the first statement made in the beginning of this subsection.  More precisely,
\begin{prop}\label{p5.8}
Let $0<\bar{\eps}<1$.  If $\sigma$ and $\eps$ are small enough depending only on $n$ and $\bar{\eps}$,  we have:
\begin{equation*}
m(S_{\mu}(x_0)\cap A_1\cap B_{r_1}(x_0))<\bar{\eps}m(S_{\mu}(x_0)),\,\,\,\text{for any $\frac{\mu_0}{484}\le \mu\le \frac{\mu_0}{4}$}.
\end{equation*}
\end{prop}
\begin{proof}
We fix some $\mu$ between $\frac{\mu_0}{484}$ and $\frac{\mu_0}{4}$.  We may also assume $\mu_0\le \frac{1}{484}$ so that $\mu_0^2<\mu\le \mu_0$,  so that $E_{4\mu}(x_0)$,  $E_{\mu}(x_0)$,  and $E_{\mu_0}(x_0)$ have the same coefficients.  

Recall from Section 2 that $\tilde{T}_{1,x_0}$ is the coordinate change such that $z=x_0+\sqrt{4\mu}\tilde{T}_{1,x_0}(w)$ will make $E_{4\mu}(x_0)$ becomes $B_{1}$ under $w$.  Denote $\tilde{\Omega}$ to be the image of $S_{4\mu}(x_0)$ under $w$.  Then from $(1-0.1\sigma)E_{4\mu}(x_0)\subset S_{4\mu}(x_0)\subset (1+0.1\sigma)E_{4\mu}(x_0)$,  we see that:
\begin{equation}\label{5.8NN}
B_{1-0.1\sigma}\subset \tilde{\Omega}\subset B_{1+0.1\sigma}.
\end{equation}
Also we define $\tilde{u}=\frac{1}{4\mu}(u-h_{4\mu,x_0})(x_0+\sqrt{4\mu}\tilde{T}_{1,x_0}(w))$,  then $\tilde{u}$,  $\tilde{\Omega}$ will fullfil the assumptions we made in Lemma \ref{l5.2}.  Now we are in a position to apply Corollary \ref{c5.7} to conclude that:
\begin{equation}\label{5.8}
\frac{m(B_{\frac{1}{2}}\cap \tilde{D})}{m(B_{\frac{1}{2}})}\ge 1-C_{6.n}\sigma^{\frac{1}{2}}-C_{6,n}\eps^{\frac{1}{2}}.
\end{equation}
Here $\tilde{D}$ is the subset of $B_{0.8}$ (under the $w$ variable) such that $w_0\in \tilde{D}$ if and only if $\tilde{S}_{\tilde{\mu}}(w_0)\subset B_{\sqrt{M_1\tilde{\mu}}}(w_0)$ for any $0<\tilde{\mu}\le \mu_0$,  where $M_1=\frac{2(1+0.1\sigma)^2}{(1-0.1C_n\sigma^{\frac{1}{2}})(1-0.2\sigma)^2}$ is given by Corollary \ref{c5.7},  but with $\gamma$ replaced by $0.1\sigma$ because of (\ref{5.8NN}).  Here $\tilde{S}_{\tilde{\mu}}(w_0)$ is the section given by the construction of Section 2,  but carried out for $\tilde{u}$.  First,  since $E_{\mu}(x_0)$ and $E_{4\mu}(x_0)$ have the same coefficients,  $E_{\mu}(x_0)$ is now $B_{\frac{1}{2}}$ under $w$.  Therefore,  if we define $\tilde{\Omega}_1$ to be the image of $S_{\mu}(x_0)$ under $w$,  we would get that:
\begin{equation}\label{5.9}
\frac{m(\tilde{\Omega}_1\cap \tilde{D})}{m(\tilde{\Omega}_1)}\ge 1-C_{7,n}\sigma^{\frac{1}{2}}-C_{7,n}\eps^{\frac{1}{2}}.
\end{equation}

Next we would like to translate the set $\tilde{D}$ back to $z$ variable,  and show that:
\begin{equation}\label{5.10}
\text{the image of $\tilde{D}$ under $z$ variable}\subset D_1,  
\end{equation}
where $D_1$ is defined in Definition \ref{def ak}.  This would finish the proof.

In order to show (\ref{5.10}),  we take $w_0\in \tilde{D}$,  then we have $\tilde{S}_{\tilde{\mu}}(w_0)\subset B_{\sqrt{M_1\tilde{\mu}}}(w_0)$ for $0<\tilde{\mu}\le \mu_0$,  and $(1-0.1\sigma)\tilde{E}_{\tilde{\mu}}(w_0)\subset \tilde{S}_{\tilde{\mu}}(w_0)\subset (1+0.1\sigma)\tilde{E}_{\tilde{\mu}}(w_0)$.
Now we switch back to $z$ coordinates,  then the above inclusions become: \begin{equation*}\begin{split}&S'_{4\mu\tilde{\mu}}(z_0)\subset x_0+\sqrt{4\mu}\tilde{T}_{1,x_0}(B_{\sqrt{M_1\tilde{\mu}}}(w_0))\subset B_{(1+C_n\sigma^{\frac{1}{2}})\sqrt{4M_1\mu\tilde{\mu}}}(z_0),\\
&(1-0.1\sigma)E_{4\mu\tilde{\mu}}'(z_0)\subset S'_{4\mu\tilde{\mu}}
(z_0)\subset (1+0.1\sigma)E_{4\mu\tilde{\mu}}'(z_0).
\end{split}
\end{equation*}
In the above,  $z_0$ is the image of $w_0$ under $z$ and $S'_{4\mu\tilde{\mu}}$ is of the form $\{u-h\le u(z_0)+4\mu\tilde{\mu}\}$,  and $E'_{4\mu\tilde{\mu}}$ is an ellipsoid having the same volume as $B_{\sqrt{4\mu\tilde{\mu}}}$. 
Now we may use Corollary \ref{c3.17} to conclude that $(1-0.1\sigma)E_{4\mu\tilde{\mu}}(z_0)\subset (1+0.1\sigma)E'_{4\mu\tilde{\mu}}(z_0)$,  where $E_{4\mu\tilde{\mu}}$ is the ellipsoid given by Section 2,  but constructed directly with $z_0$ (under $z$ coordinate).
From this we see:
\begin{equation*}
S_{4\mu\tilde{\mu}}(z_0)\subset (1+0.1\sigma)E_{4\mu\tilde{\mu}}(z_0)\subset \frac{(1+0.1\sigma)^2}{1-0.1\sigma}E_{4\mu\tilde{\mu}}'(z_0)\subset B_{(1+C_n\sigma^{\frac{1}{2}})\frac{(1+0.1\sigma)^2}{(1-0.1\sigma)^2}\sqrt{4M_1\mu\tilde{\mu}}}.
\end{equation*}
This would give the control of $S_{\mu'}(z_0)$ for $\mu'\le 4\mu\mu_0$,  and $4\mu\mu_0\ge \frac{\mu_0^2}{121}$.

For $\mu'\ge \frac{\mu_0^2}{121}$,  we may assume without loss of generality that $\mu_0^3<\mu'\le \mu_0^2$,  then 
\begin{equation*}
S_{\mu'}(z_0)\subset (1+0.1\sigma)E_{\mu'}(z_0)=z_0+(1+0.1\sigma)T_{2,z_0}(B_{\sqrt{\mu'}}(0))\subset B_{(1+0.1\sigma)(1+C_n\sigma^{\frac{1}{2}})^2\sqrt{\mu'}}(z_0).
\end{equation*}
Hence we see that if we have:
\begin{equation}\label{5.12}
10\ge \max\big((1+C_n\sigma^{\frac{1}{2}})\sqrt{M_1},\,(1+0.1\sigma)(1+C_n\sigma^{\frac{1}{2}})^2\big),  
\end{equation}
we can conclude that $z_0\in D_1$,  thereby finishing the proof.  This is indeed true if we choose $\sigma$ small enough depending only on $n$.  Moreover,  we need to take $\sigma$ and $\eps$ so that in (\ref{5.9}),  we have $C_{7,n}\sigma^{\frac{1}{2}}+C_{7,n}\eps^{\frac{1}{2}}<\bar{\eps}.$
\end{proof}

Now we move on to prove Statement (2) made in the beginning of Subsection 6.1.

First,  we need to show that,  if $m(S_{\mu}(x_0)\cap A_{k+1}\cap B_{r_{k+1}}(0))\ge \bar{\eps}m(S_{\mu}(x_0))$,  then $S_{\mu}(x_0)\subset  B_{r_k}(0)$,  with $\sigma$ and $\eps$ small enough.

For this we observe that:
\begin{lem}\label{l5.9}
Let $0<\mu\le \frac{\mu_0}{4}$ and $x_0\in B_{0.8}$ such that $S_{4\mu}(x_0)\subset B_{0.8}$.   Assume that there exists $\Lambda>0$ such that $||T_{\mu_1,x_0}||\le \Lambda$ for all $4\mu\le \mu_1\le \mu_0$.  Then for $\sigma>0$,  $\eps>0$ chosen small enough depending only on $n$ and $\bar{\eps}$,  we have 
\begin{equation*}
\frac{m(S_{\mu}(x_0)\cap D_{2\Lambda})}{m(S_{\mu}(x_0))}\ge 1-\frac{\bar{\eps}}{2}.
\end{equation*}
Here $D_{2\Lambda}$ is the set of $z_0$ such that $S_{\mu'}(z_0)\subset B_{2\Lambda\sqrt{\mu'}}(z_0)$ for all $0<\mu'\le \mu_0$.
\end{lem}
\begin{proof}
We consider the section $S_{4\mu}(x_0)$ with $\mu_0^{k+1}<4\mu\le \mu_0^k$.  From the assumption we see that $||T_{4\mu,x_0}||\le\Lambda$.

Similar to the proof of Proposition \ref{p5.8},  we consider the change of coordinates $z=x_0+\sqrt{4\mu}T_{4\mu,x_0}(w)$,  so that $E_{4\mu}(x_0)$ gets transformed to be a unit ball.  We denote $\tilde{\Omega}$ to be the image of $S_{4\mu}(x_0)$ under $w$ and $\tilde{\Omega}_1$ to be the image of $S_{\mu}(x_0)$.  Similar to the proof of Proposition \ref{p5.8}.
\begin{equation*}
\frac{m(\tilde{\Omega}_1\cap\tilde{D})}{m(\tilde{\Omega}_1)}\ge 1-C_{7,n}\sigma^{\frac{1}{2}}-C_{7,n}\eps^{\frac{1}{2}}.
\end{equation*}
Here $\tilde{D}$ is the subset of $B_{0.8}$ under the $w$ variable such that $w_0\in \tilde{D}$ if and only if $\tilde{S}_{\tilde{\mu}}(w_0)\subset B_{\sqrt{M_1\tilde{\mu}}}(w_0)$ for $0<\tilde{\mu}\le \mu_0$,  where $M_1=\frac{2(1+0.1\sigma)^2}{(1-0.1C_n\sigma^{\frac{1}{2}})(1-0.2\sigma)^2}$.  We repeat the argument in the proof of Proposition \ref{p5.8} below (\ref{5.10}),  transform this containment back to $z$ variable and conclude that 
\begin{equation}\label{5.13}
S_{\mu'}(z_0)\subset B_{(1+C_n\sigma^{\frac{1}{2}})\frac{(1+0.1\sigma)^2}{(1-0.1\sigma)^2}\Lambda\sqrt{M_1\mu'}},\,\,\,0<\mu'\le 4\mu\mu_0.
\end{equation}
Now it only remains to control $S_{\mu'}(z_0)$ for $\mu'>4\mu\mu_0$.

First,  if $\mu'\ge 4\mu$,  then we have $z_0\in S_{\mu'}(x_0)$.  Then from Lemma \ref{l3.18},  we know that:
\begin{equation*}
||T_{\mu',x_0}^{-1}\circ T_{\mu',z_0}||\le 1.1^3,\,\,\,\,||T_{\mu',z_0}^{-1}\circ T_{\mu',x_0}||\le 1.1^3.
\end{equation*}
Hence
\begin{equation*}
||T_{\mu',z_0}||\le 1.1^3\Lambda,\,\,\,\mu'\ge 2\mu.
\end{equation*}
Therefore,  we have:
\begin{equation}\label{5.14}
S_{\mu'}(z_0)\subset (1+0.1\sigma)E_{\mu'}(z_0)\subset B_{1.1^3(1+0.1\sigma)\Lambda\sqrt{\mu'}}(z_0),\,\,\,\mu'\ge 2\mu.
\end{equation}
Finally,  for $2\mu\mu_0<\mu'\le 2\mu$,  we note that $T_{2\mu,x_0}$ and $T_{\mu',x_0}$ differ at most by $\tilde{T}_{k+1,x_0}$,  whose norm is bounde by $1+C_n\sigma^{\frac{1}{2}}$,  hence we have:
\begin{equation}\label{5.15}
S_{\mu'}(z_0)\subset B_{1.1^3(1+0.1\sigma)(1+C_n\sigma^{\frac{1}{2}})\Lambda\sqrt{\mu'}}(z_0),\,\,\,2\mu\mu_0<\mu'<2\mu.
\end{equation}
Hence if one combines (\ref{5.13})-(\ref{5.15}) and define 
\begin{equation}\label{5.16}
M_2=\max\big((1+C_n\sigma^{\frac{1}{2}})^2\frac{(1+0.1\sigma)^4}{(1-0.1\sigma)^4}M_1,\,1.1^6(1+0.1\sigma)^2(1+C_n\sigma^{\frac{1}{2}})^2\big),
\end{equation}
then for any $z_0$ in the image of $\tilde{D}$ in the $z$ coordinate,  $S_{\mu'}(z_0)\subset B_{\sqrt{M_2\Lambda^2\mu'}}(z_0)$,  $0<\mu'\le \mu_0$.  Also we note that the right hand side of (\ref{5.16}) is less than 4,  if $\sigma$ is small enough depending only on $n$.
\end{proof}

As a consequence,  we get that
\begin{cor}\label{c5.10}
Let $0<\mu\le \frac{\mu_0}{4}$ and $x_0\in B_{0.8}$ such that $S_{4\mu}(x_0)\subset B_{0.8}$.  Assume that for some $k_0\ge 1$,  $m(S_{\mu}(x_0)\cap A_{k_0+1}\cap B_{r_{k_0+1}}(0))\ge \bar{\eps}m(S_{\mu}(x_0))$.  Assume also that $\sigma$ and $\eps$ are chosen small enough depending on $\bar{\eps}$ and $n$.  Then we have:
\begin{equation*}
\mu\le 10^{-\frac{|\log(\mu_0)|}{2\log(1+C_n\sigma^{\frac{1}{2}})}k_0}.
\end{equation*}
\end{cor}
\begin{proof}
Let $k_1\ge 1$ be such that $\mu_0^{k_1+1}<4\mu\le \mu_0^{k_1}$.  Then $T_{4\mu,x_0}=T_{k_1,x_0}$ and we have the estimate:
\begin{equation*}
||T_{k,x_0}||,\,\,\,||T_{k,x_0}^{-1}||\le (1+C_n\sigma^{\frac{1}{2}})^{k_1},\,\,\,1\le k\le k_1.
\end{equation*}
Therefore,  we may use Lemma \ref{l5.9},  with $\Lambda=(1+C_n\sigma^{\frac{1}{2}})^{k_1}$,  to conclude that:
\begin{equation*}
\frac{m(S_{\mu}(x_0)\cap D_{2\Lambda})}{m(S_{\mu}(x_0))}\ge 1-\frac{\bar{\eps}}{2}.
\end{equation*}
From our assumption that $m(S_{\mu}(x_0)\cap A_{k_0+1}\cap B_{r_{k_0+1}}(0))\ge \bar{\eps}m(S_{\mu}(x_0))$,  we must have $D_{k_0+1}\subset D_{2\Lambda}$ (note that we either have $D_{k_0+1}\subset D_{2\Lambda}$ or $D_{2\Lambda}\subset D_{k_0+1}$ by definition),  in other words, 
\begin{equation}\label{5.17}
10^{k_0+1}\le 4(1+C_n\sigma^{\frac{1}{2}})^{2k_1}.
\end{equation}
So that
\begin{equation*}
\mu\le \frac{1}{4}\mu_0^{k_1}\le 10^{-\frac{|\log(\mu_0)|}{2\log(1+C_n\sigma^{\frac{1}{2}})}k_0}.
\end{equation*}
\end{proof}

As a further corollary,  we see that
\begin{cor}\label{c5.11}
Assume that for some $k_0\ge 1$,  $m(S_{\mu}(x_0)\cap A_{k_0+1}\cap B_{r_{k_0+1}}(0))\ge \bar{\eps}m(S_{\mu}(x_0))$ for some $x_0\in B_{0.8}$ and $\mu\le \frac{\mu_0}{4}$,  then $S_{\mu}(x_0)\subset B_{r_{k_0}}(0)$,  if $\sigma$ and $\eps$ are chosen small enough depending on $n$ and $\bar{\eps}$.
\end{cor}
\begin{proof}
By our assumption,  $S_{\mu}(x_0)\cap B_{0.7}\neq \emptyset$,  hence we know that $S_{4\mu}(x_0)\subset B_{0.8}$ (by choosing $\sigma$,  hence $\mu_0$ small enough.)
Also we denote $k_1$ so that $\mu_0^{k_1+1}<4\mu\le \mu_0^{k_1}$.  
Then we may use Corollary \ref{c5.10} to obtain that:
\begin{equation}\label{5.18}
\begin{split}
&diam\,S_{\mu}(x_0)\le diam\,E_{4\mu}(x_0)\le (4\mu)^{\frac{1}{2}}||T_{4\mu,x_0}||\le \mu_0^{k_1}(1+C_n\sigma^{\frac{1}{2}})^{k_1}\\
&\le \big(\mu_0(1+C_n\sigma^{\frac{1}{2}})\big)^{\frac{\log 10}{2\log(1+C_n\sigma^{\frac{1}{2}})}k_0}\le \frac{1}{40}\cdot 2^{-k_0}.
\end{split}
\end{equation}
The above is true if we choose $\sigma$ small enough depending on $n$ (according to our choice of $\mu_0$ made in (\ref{3.2}),  $\mu_0\le \sigma^2$.)

On the other hand,  $r_{k_0}-r_{k_0+1}=\frac{1}{10}\cdot 2^{-k_0-1}$.  Since $S_{\mu}(x_0)\cap B_{r_{k_0+1}}(0)\neq \emptyset$,  the conclusion follows from (\ref{5.18}).
\end{proof}
The remaining part of Statement (2) is to show that $S_{\mu}(x_0)\subset A_k$.  This would directly follow from the following observation:
\begin{lem}\label{l5.12}
Let $0<\mu\le \frac{\mu_0}{4}$ and $x_0\in B_{0.8}$ be such that $S_{4\mu}(x_0)\subset B_{0.8}$ and $S_{\mu}(x_0)\cap D_{k_0}\neq \emptyset$,  for some $k_0\ge 1$.  Then we have:
\begin{equation*}
m(S_{\mu}(x_0)\cap D_{k_0+1})\ge (1-\frac{\bar{\eps}}{2})m(S_{\mu}(x_0)),
\end{equation*}
provided that $\sigma$ is small enough depending only on $\bar{\eps}$ and $n$.
\end{lem}
\begin{proof}
From the assumption,  we can find $x_*\in S_{\mu}(x_0)\cap D_{k_0}$,  which means that $S_{\mu'}(x_*)\subset B_{\sqrt{10^{k_0}\mu'}}(x_*)$ for any $0<\mu'\le \mu_0$.  Since $(1-0.1\sigma)E_{\mu'}(x_*)\subset S_{\mu'}(x_*)$ and $T_{\mu',x_*}(B_{\sqrt{\mu'}}(0))=E_{\mu'}(0)$,  we see that 
\begin{equation*}
||T_{\mu',x_*}||\le \frac{\sqrt{10^{k_0}}}{1-0.1\sigma},\,\,\,\,0<\mu'\le \mu_0.
\end{equation*}
Since $x_*\in S_{\mu}(x_0)$,  we know that $x_*\in S_{\mu'}(x_0)$ for all $\mu'\ge 4\mu$.  Hence we may use Lemma \ref{l3.18} to conclude that:
\begin{equation*}
||T_{\mu',x_*}^{-1}\circ T_{\mu',x_0}||\le 1.1^3,\,\,\,||T_{\mu',x_0}^{-1}\circ T_{\mu',x_*}||\le 1.1^3.
\end{equation*}
Hence we obtain that 
\begin{equation*}
||T_{\mu',x_0}||\le \frac{1.1^3\sqrt{10^{k_0}}}{1-0.1\sigma},\,\,\,4\mu\le \mu'\le \mu_0.
\end{equation*}
Now we are in a position to use Lemma \ref{l5.9} with $\Lambda=\frac{1.1^3\sqrt{10^{k_0}}}{1-0.1\sigma}.$ Then we conclude that:
\begin{equation*}
\frac{m(S_{\mu}(x_0)\cap D_{2\Lambda})}{m(S_{\mu}(x_0))}\ge 1-\frac{\bar{\eps}}{2}.
\end{equation*}
Here $D_{2\Lambda}$ is the set of $z_0$ such that $S_{\mu'}(z_0)\subset B_{2\Lambda\sqrt{\mu'}}(z_0)$ for all $0<\mu'\le \mu_0$.
We will be done if we can ensure that $D_{2\Lambda}\subset D_{k_0+1}$,  and we just need:
\begin{equation*}
4\Lambda^2=4\cdot \frac{1.1^6\cdot 10^{k_0}}{(1-0.1\sigma)^2}\le 10^{k_0+1}.
\end{equation*}
This is clear if $\sigma$ is small enough.
\end{proof}
Now we are ready to show the Statement (2) in the beginning of Subsection 6.1.
\begin{prop}\label{p5.13}
Let $0<\mu\le \frac{\mu_0}{4}$,   $x_0\in B_{0.8}$ be such that $m(S_{\mu}(x_0)\cap A_{k_0+1}\cap B_{r_{k_0+1}})\ge \bar{\eps}m(S_{\mu}(x_0))$ for some $k_0\ge 1$,  then $S_{\mu}(x_0)\subset A_{k_0}\cap B_{r_{k_0}}(0)$,  if $\sigma$ and $\eps$ are small enough depending on $\bar{\eps}$ and $n$.
\end{prop}
\begin{proof}
Corollary \ref{c5.11} already implies $S_{\mu}(x_0)\subset B_{r_k}(0)$,  and we just have to show $S_{\mu}(x_0)\subset A_{k_0}$.

If not,  namely $S_{\mu}(x_0)\cap D_{k_0}\neq \emptyset$,  then Lemma \ref{l5.12} would give us: 
\begin{equation*}
m(S_{\mu}(x_0)\cap D_{k_0+1})\ge (1-\frac{\bar{\eps}}{2})m(S_{\mu}(x_0)).
\end{equation*}
This is in contradiction with our assumption that: $m(S_{\mu}(x_0)\cap A_{k_0+1}\cap B_{r_{k_0+1}})\ge \bar{\eps}m(S_{\mu}(x_0))$.
\end{proof}
Now we are ready to show that:
\begin{thm}\label{t5.1}
Let $A_k$ be define by Definition \ref{def ak}.  Let $0<\bar{\eps}<1$ be given.  Let $\sigma>0$ and $\eps>0$ be small enough depending on $n$ and $\eps$,  we have:
\begin{equation*}
m(A_k\cap B_{r_k}(0))\le m(B_{0.7})(12^{2n}\bar{\eps})^{k-1},\,\,\,k\ge 1.
\end{equation*}
Here $r_k=r_{k-1}-\frac{1}{10}\cdot 2^{-k}$,  $r_0=0.7$.
In particular,  $m(A_k\cap B_{0.6})\le m(B_{0.7})(12^{2n}\bar{\eps})^{k-1}.$
\end{thm}
\begin{proof}
For $k=1$,  the above estimate is trivial.

For $k\ge 1$,  we have the following estimate holds:
\begin{equation*}
m(A_{k+1}\cap B_{r_{k+1}}(0))\le 12^{2n}\bar{\eps}m(A_k\cap B_{r_k}(0)).
\end{equation*}
This follows from taking $X=A_{k+1}\cap B_{r_{k+1}}(0)$,  $Y=B_{r_k}(0)\cap B_{r_k}(0)$ in Lemma \ref{l0.11},  where the two assumptions of that lemma indeed hold,  due to Proposition \ref{p5.8} and \ref{p5.13}.
\end{proof}

\subsection{Control of the second derivatives on the good set and completion of proof}
In this section,  we show that the derivatives are controlled on the good sets $D_k$.  
We wish to emphasize that we are assuming $u\in C^2(\Omega)$,  only for the sake of convenience of argument,  but the regularity of $u$ does not go into the quantitative estimates.

We start with the following lemma:
\begin{lem}
Let $u\in C^2(B_1)$ solving $\det u_{i\bar{j}}=f$ with $|f-1|<\eps$.  Let $x_0\in D_k$ for some $k\ge 1$.  Then 
\begin{equation*}
\frac{1}{10^k}I_n\le u_{i\bar{j}}(x_0)\le 2\cdot 10^{(n-1)k}I_{n}.
\end{equation*}
\end{lem}
\begin{proof}
From Definition \ref{def ak},  we know that $S_{\mu}(x_0)\subset B_{\sqrt{10^k\mu}}(x_0)$.  Let $0<c<1$,  then we have,  on $\partial S_{\mu}(x_0)$,
\begin{equation}\label{5.19}
u-h_{\mu,x_0}=u(x_0)+\mu\ge u(x_0)+\frac{|x-x_0|^2}{10^k}>u(x_0)+\frac{c|x-x_0|^2}{10^k}.
\end{equation}
In the above,  we noted that $|x-x_0|^2\le 10^k\mu$ on $\partial S_{\mu}(x_0)$.

On the other hand,  (\ref{5.19}) achieves equality when $x=x_0$.  Hence the function \sloppy $x\mapsto (u-h_{\mu,x_0})(x)-\frac{c|x-x_0|^2}{10^k}$ achieves minimum in the interior of $S_{\mu}(x_0)$,  say $x_{\mu}$.  Then we have:
\begin{equation}
u_{i\bar{j}}(x_{\mu})\ge \frac{c}{10^k}I_n.
\end{equation}
Since $\det u_{i\bar{j}}\le 1+\eps$,  we see that
\begin{equation*}
u_{i\bar{j}}(x_{\mu})\le (1+\eps)c^{1-n}10^{k(n-1)}I_n.
\end{equation*}
Since $diam\,S_{\mu}(x_0)\rightarrow 0$ and $x_{\mu}\in S_{\mu}(x_0)$,  we see that $x_{\mu}\rightarrow x_0$ as $\mu\rightarrow 0$.  Hence we conclude that:
\begin{equation*}
\frac{c}{10^k}I_n\le u_{i\bar{j}}(x_0)\le 2\cdot c^{1-n}10^{k(n-1)}I_n.
\end{equation*}
Since $0<c<1$ is arbitrary,  we can make $c\rightarrow 1$ to conclude the result.
\end{proof}

As a consequence,  we get that:
\begin{prop}\label{p5.15}
Let $p>1$ be given.  Let $u$ be given by Theorem \ref{main theorem baby} and is $C^2(B_1)$.  Let $\gamma$ be small enough depending only on $n$,  and $\eps$ small enough depending on $p$ and $n$,  then we have:
\begin{equation*}
\int_{B_{0.6}}\big((\Delta u)^p+(tr_u\omega_E)^p\big)\le C_{p,n}.
\end{equation*}
Here $\Delta u=\sum_iu_{i\bar{i}}$,  $tr_u\omega_E=\sum_i\frac{1}{u_{i\bar{i}}}$.  
\end{prop}
\begin{proof}
Note that on $D_k$,  we have:
\begin{equation*}
\Delta u\le 2n\cdot 10^{(n-1)k},\,\,\,tr_u\omega_E\le n\cdot 10^k.
\end{equation*}
Hence $\{\Delta u>2n\cdot 10^{(n-1)k}\}$ and $\{tr_u\omega_E>n\cdot 10^k\}$ are contained in $A_k$ for $k\ge 1$.
\begin{equation*}
\begin{split}
&\int_{B_{0.6}\cap \{\Delta u>2n\cdot 10^{n-1}\}}(\Delta u)^p\le \sum_{k=1}^{\infty}\int_{B_{0.6}\cap \{2n\cdot 10^{(n-1)k}<\Delta u\le 2n\cdot 10^{(n-1)(k+1)}\}}(\Delta u)^p\\
&\le \sum_{k=1}^{\infty}\big(2n\cdot 10^{(n-1)(k+1)}\big)^pm(A_k\cap B_{0.6})\le \sum_{k=1}^{\infty}\big(2n\cdot 10^{(n-1)(k+1)}\big)^p\cdot m(B_{0.7})(12^{2n}\bar{\eps})^{k-1}\\
&=(2n\cdot 10^{2(n-1)})^pm(B_{0.7})\sum_{k=1}^{\infty}(10^{(n-1)p}12^{2n}\bar{\eps})^{k-1}.
\end{split}
\end{equation*}
In order for the above sum to be finite,  we can choose $\bar{\eps}>0$ so that $10^{(n-1)p}12^{2n}\bar{\eps}=\frac{1}{2}$ so that the above integral $\le 2(2n\cdot 10^{2(n-1)})^pm(B_{0.7})$.
In the above,  we used Theorem \ref{t5.1}.  In order for the theorem to apply,  we need to choose $\eps$ and $\sigma$ small enough depending only on $\bar{\eps}$ and $n$,  so the choice of $\eps$ and $\sigma$ eventually depend on $p$ and $n$.

The estimation for $tr_u\omega_E$ is completely similar.  
\end{proof}
From Proposition \ref{p5.15},  we can get full second order estimate by applying the $L^p$-estimate for the Laplacian.
That is,  we get:
\begin{cor}
Under the assumption of Proposition \ref{p5.15},  we have:
\begin{equation*}
||u||_{W^{2,p}(B_{\frac{1}{2}})}\le C_{p,n},
\end{equation*}
if $\gamma$ is small enough depending only on $n$,  and $\eps$ small enough depending only on $p$ and $n$.
\end{cor}
\begin{proof}
The result would follow from Proposition \ref{p5.15} and the classical $W^{2,p}$ estimates (Gilbarg-Trudinger \cite{GT},  Chapter 9): for any $u\in C^2(B_{0.6})$ and $1<p<\infty$,

\begin{equation*}
||u||_{W^{2,p}(B_{\frac{1}{2}})}\le C(||u||_{L^p(B_{0.6})}+||\Delta u||_{L^p(B_{0.6})}).
\end{equation*}
\end{proof}

\section{Some corollaries of the main theorem}
First we prove Corollary \ref{c1.2}.  
\begin{proof}(of Corollary \ref{c1.2})
First,  we wish to use the following Lemma \ref{l7.1} to conclude that $\varphi$ is  close to zero.
Indeed,  by taking $g=1$,  $\psi=0$,  $p=2$ in the following,  we find that:
\begin{equation*}
||\varphi||_{L^{\infty}}\le c\eps^{\frac{1}{2(n+4)}}.
\end{equation*}
Here $c$ depends only on the background metric.

Now we take some point $p_0\in M$ and take normal coordinates $(z_1,\cdots,  z_n)$ at $p_0$ so that $g_{i\bar{j}}(p_0)=\delta_{ij}$ and $\nabla g(p_0)=0$.  We can choose local potential $\rho(z)$,  such that $\omega_0=\sqrt{-1}\partial\bar{\partial}\rho$ near $p_0$,  say on $B_1(p_0)$ (under local coordinates $z$).  So that on this neighborhood,  the equation can be written as:
\begin{equation}\label{7.1}
\det((\rho+\varphi)_{i\bar{j}})=
f\det(g_{i\bar{j}}),\,\,\,\text{ in $B_1$.}
\end{equation}
In order to use Theorem \ref{main theorem},  we need to zoom in (\ref{7.1}) at $p_0$ at a suitable scale so that the right hand side is close to a contant.  

Denote $u=\rho+\varphi$.  
Let $0<r_0<1$,  we perform a change of variable $z=r_0w$.  Next we define 
\begin{equation*}
\tilde{u}_{r_0}(w)=\frac{1}{r_0^2}u(r_0w),\,\,\,\tilde{\rho}_{r_0}=\frac{1}{r_0^2}\rho(r_0w),\,\,\tilde{\varphi}_{r_0}(w)=\frac{1}{r_0^2}\varphi(r_0w).
\end{equation*}
Assume that for some $C_0>0$,  we have $\frac{1}{C_0}I\le \rho_{i\bar{j}}\le C_0I$,  $|D^3\rho|\le C_0$ on $B_1$,  then we see that with the same $C_0>1$:
\begin{equation*}
\frac{1}{C_0}I\le (\tilde{\rho}_{r_0})_{w_i\bar{w}_j}\le C_0I,\,\,|D^3_w\tilde{\rho}_{r_0}|\le C_0,\,\,\,\text{ for $|w|<1$}.
\end{equation*}
Also
\begin{equation*}
|\tilde{u}_{r_0}-\tilde{\rho}_{r_0}|\le \frac{1}{r_0^2}||\varphi||_{L^{\infty}}\le c\frac{1}{r_0^2}\eps^{\frac{1}{2(n+4)}}.
\end{equation*}
Also on $\{|w|<1\}$,  we have
\begin{equation*}
\det(\tilde{u}_{r_0})_{i\bar{j}}=f\det (g_{i\bar{j}})(r_0w),
\end{equation*}
and we can estimate how close the right hand side is from 1:
\begin{equation*}
\begin{split}
|f\det(g_{i\bar{j}})-1|\le |f-1|\det g_{i\bar{j}}(r_0w)+f|\det(g_{i\bar{j}})(r_0w)-1|\le C_1\eps+C_1r_0.
\end{split}
\end{equation*}
Here $C_1$ depends only on the background metric and $n$.

Hence in order to apply Theorem \ref{main theorem},  we need to make sure that:
\begin{equation*}
\begin{split}
&C_1\eps+C_1r_0\le \eps_{p,n},\\
&c\frac{1}{r_0^2}\eps^{\frac{1}{2(n+4)}}\le \delta_0.
\end{split}
\end{equation*}
Here $\eps_{p,n}$ and $\delta_0$ are determined by Theorem \ref{main theorem}.  Hence we need to choose $r_0$ first so that $C_1r_0=\frac{1}{2}\eps_{p,n}$.  Then we fix this choice,  and choose $\eps$ small enough so as to make sure $C_1\eps\le \frac{1}{2}\eps_{p,n}$ and $c\frac{1}{r_0^2}\eps^{\frac{1}{2(n+4)}}\le \delta_0$.
Then Theorem \ref{main theorem} gives $W^{2,p}$ estimate for $\tilde{u}_{r_0}$ in $B_{\frac{1}{2}}$.  Scaling back to $u$,  we get $W^{2,p}$ estimate for $u$ in $B_{\frac{1}{2}r_0}(p_0)$.
\end{proof}

In the above proof,  we used the following stability estimate due to Kolodziej \cite{KO}
\begin{lem}(\cite{KO})\label{l7.1}
Let $(M,\omega_0)$ be a compact K\"ahler manifold.  Let $\varphi\in PSH(M,\omega_0)$ and $\psi\in PSH(M,\omega_0)$ be the solution to the complex Monge-Ampere equations:
\begin{equation*}
\begin{split}
&(\omega_0+\sqrt{-1}\partial\bar{\partial}\varphi)^n=f\omega_0^n,\,\,\,\sup_M\varphi=0,\\
&(\omega_0+\sqrt{-1}\partial\bar{\partial}\psi)^n=g\omega_0^n,\,\,\,\sup_M\psi=0.
\end{split}
\end{equation*}
Assume that there is some $c_0>0,\,p>1$ such that
\begin{equation*}
||f||_{L^p}\le c_0,\,\,\,||g||_{L^p}\le c_0.
\end{equation*}
Then we have
\begin{equation*}
\sup_M|\varphi-\psi|\le c(c_0,p)||f-g||_{L^p}^{\frac{1}{n+4}}.
\end{equation*}
\end{lem}
Now we prove Corollary \ref{c1.2N}.  The argument is similar to Corollary \ref{c1.2N}.
\begin{proof}
(of Corollary \ref{c1.2N})
We first use Lemma \ref{local stability} with $q=2$ to estimate the difference between $\varphi$ and $\varphi_0$.  We have that 
\begin{equation*}
\sup_{\Omega}|u-u_0|\le \eps+c\eps^{\frac{1}{2n}}\le c'\eps^{\frac{1}{2n}}.
\end{equation*}
Here $c$ depends only on $n$ and $diam\,\Omega$.  The rest of the argument is similar to Corollary \ref{c1.2}.  In other words,  we take any point $z_0\in\Omega'$ and consider rescaling $\tilde{u}_{r_0}(w)=\frac{1}{r_0^2}u(z_0+r_0w)$,  and we similarly consider $\tilde{u}_{0,r_0}=\frac{1}{r_0^2}u_0(z_0+r_0w)$.  Then for suitable chosen $r_0$,  we will have that $\tilde{u}_{r_0}$ satisfy the hypothesis of Theorem \ref{main theorem}.
\end{proof}
Next,  we prove the Liouville theorem Corollary \ref{c1.3}.
\begin{proof}
(of Corollary \ref{c1.3}) 

Let $r>1$,  we consider:
\begin{equation*}
u_r(w)=\frac{1}{r^2}u(rw).
\end{equation*}
Then we know that for $r$ large enough,  one has
\begin{equation*}
|u_r(w)-|w|^2|\le 2\eps,\,\,\,w\in B_1.
\end{equation*}
On the other hand,  we have $1-\eps\le \det u_{a\bar{b}}\le 1+\eps$,  hence from Theorem \ref{main theorem} we obtain:
\begin{equation*}
||u_r||_{W^{2,p}(B_{\frac{1}{2}})}\le C_{p,n},
\end{equation*}
as long as we choose $\eps$ small enough depending on $n$ and $p$.  On the other hand,  $u_r$ will also satisfy the scalar flat equation,  hence we may use the following Lemma \ref{l7.2},  and choose $p=p_n$,  then we obtain that:
\begin{equation*}
||u_r||_{C^{2,\alpha}(B_{\frac{1}{4}})}\le C.
\end{equation*}
Here $C$ is a uniform constant independent of $r$.  
Rescaling back to $u$,  it gives:
\begin{equation*}
|D^2u(x)-D^2u(y)|r^{\alpha}\le C|x-y|^{\alpha},\,\,\forall x,\,y\in B_{\frac{r}{2}}.
\end{equation*}
One can fix $x$,  $y$ and let $r\rightarrow \infty$ and conclude that $D^2u$ is a constant,  hence $u$ is a quadratic polynomial.
\end{proof}

In the following,  we used the following estimate of the scalar flat equation,  which originates from Chen-Cheng \cite{cc1},  Corollary 6.2.  Note that this is from the preprint version on arxiv,  which was deleted in the published version.
\begin{lem}\label{l7.2}
(\cite{cc1})
Let $u\in C^4(B_1)\cap PSH(B_1)$ be a bounded solution to the scalar flat equation:
\begin{equation*}
\sum_{i,j=1}^nu^{i\bar{j}}\partial_{i\bar{j}}\big(\log\det u_{a\bar{b}}\big)=0.
\end{equation*}
Assume that for some $p>3n(n-1)$,  we have $\Delta u\in L^p(B_1)$,  $\sum_i\frac{1}{u_{i\bar{i}}}\in L^p(B_1)$.  Then for any $0<\alpha<1$,  
\begin{equation*}
||u||_{C^{2,\alpha}(B_{\frac{1}{2}})}\le C.
\end{equation*}
Here $C$ depends on $n$,  $p$,  $||u||_{L^{\infty}(B_1)}$,  $||\Delta u||_{L^p(B_1)}$,  $||\sum_i\frac{1}{u_{i\bar{i}}}||_{L^p(B_1)}$.
\end{lem}

Now let us prove the Schauder type estimate,  Corollary \ref{c1.4}.

\begin{proof}
(of Corollary \ref{c1.4})
First,  we just need to prove that $u\in W^{2,p}$ for $p$ large enough,  depending only on $n$ and $\alpha$.  Indeed,  $W^{2,p}$ embeds into $C^{1,1-\frac{2n}{p}}$ for $p>2n$ (keep in mind that the real dimension is $2n$.)
In order to apply Lemma \ref{improve regularity},  we just need to choose $p$ large,  so that
\begin{equation*}
1-\frac{2n}{p}>1-\frac{\alpha}{n(2+\alpha)-1}.
\end{equation*}

Now we fix this $p$ and it only remains to show $u\in W^{2,p}$.  The argument is similar to Corollary \ref{c1.2} and \ref{c1.2N}.  We need to do rescaling,  so that we are in the situation of Theorem \ref{main theorem}.  In other words,  we take $z_0\in B_{\frac{1}{2}}$,  and consider rescaling $\tilde{u}_{r_0}(w)=\frac{1}{r_0^2(f(z_0))^{\frac{1}{n}}}u(z_0+r_0w)$.  Similarly,  we define $\tilde{w}_{r_0}(w)=\frac{1}{r_0^2(f(z_0))^{\frac{1}{n}}}w_0(z_0+r_0w)$.  Then for $|w|<1$,  the closeness between $\tilde{u}_{r_0}$ and $\tilde{w}_{r_0}$ becomes:
\begin{equation*}
|\tilde{u}_{r_0}-\tilde{w}_{r_0}|\le \frac{\delta_0}{r_0^2(f(z_0))^{\frac{1}{n}}}\le \delta_0\frac{C_1^{\frac{1}{n}}}{r_0^2}.
\end{equation*}
On the other hand,  
\begin{equation*}
\det(\tilde{u}_{r_0})_{i\bar{j}}=\frac{f(z_0+r_0w)}{f(z_0)},
\end{equation*}
and we can estimate:
\begin{equation*}
|\frac{f(z_0+r_0w)}{f(z_0)}-1|\le \frac{1}{f(z_0)}r_0^{\alpha}K\le C_1^{\frac{1}{n}}Kr_0^{\alpha}.
\end{equation*}
Hence,  in order for Theorem \ref{main theorem} to apply,  we just need to guarantee:
\begin{equation*}
\begin{split}
&\delta_0\frac{C_1^{\frac{1}{n}}}{r_0^2}\le \delta_0',\\
&C_1^{\frac{1}{n}}Kr_0^{\alpha}\le \eps_{p,n}.
\end{split}
\end{equation*}
Here $\delta_0'$ is the required closeness from the solution to the smooth background,  given by Theorem \ref{main theorem},  and $\eps_{p,n}$ is the required small $\eps$ under the above choice of $p$.

Hence we should first choose $r_0$ so that $C_1^{\frac{1}{n}}Kr_0^{\alpha}=\frac{1}{2}\eps_{p,n}$.  Then with $r_0$ fixed,  we choose $\delta_0$ small enough so as to make sure $\delta_0\frac{C_1^{\frac{1}{n}}}{r_0^2}\le \delta_0'$.
Then we may use Theorem \ref{main theorem} to conclude that $u$ is in $W^{2,p}$ on $B_{\frac{1}{2}r_0}(z_0)$.
\end{proof}

\begin{lem}\label{improve regularity}(\cite{LLZ2})
Let $\Omega$ be a domain in $\mathbb{C}^n$ and $u\in PSH(\Omega)\cap C(\Omega)$ be a weak solution of the complex Monge-Ampere equation in $\Omega$ with $0<\lambda \le f \in C^{\alpha}(\Omega)$ for some constant $\lambda$ and some $\alpha\in (0,1)$. If $u\in C^{1,\beta}(\Omega)$ with $\beta\in (\beta_0,1)$, where 
\begin{equation*}
    \beta_0 =\beta_0 (n,\alpha)=1-\frac{\alpha}{n(2+\alpha)-1}.
\end{equation*}
Then $u\in C^{2,\alpha}(\Omega)$. Furthermore the $C^{2,\alpha}$ norm of $u$ in any relatively compact subset is estimable in terms of $n,\alpha, \beta, \lambda, ||u||_{C^{1,\beta}(\Omega)}, ||f||_{C^{\alpha}(\Omega)}$ and the distance of the set to $\partial \Omega$.
\end{lem}

\end{document}